%% file: main-func-estim-paper.tex
\documentclass[a4paper, 10pt]{article}

\input{customized_packages_and_macros}

\usetikzlibrary{%
  arrows,%
  shapes.misc,
  shapes.arrows,%
  chains,%
  matrix,%
  positioning,
  scopes,%
  decorations.pathmorphing,
  shadows,%
  backgrounds,
  fit,positioning,shapes.symbols,chains
}
\usetikzlibrary{shapes.geometric,calc}

\author{B\"arbel Holm%
\thanks{Department of Computational Science and Technology, 
School of Computer Science and Communication, 
KTH Royal Institute of Technology, SE-10044 Stockholm,
\texttt{barbel@kth.se}} \and 
Svetlana Matculevich%
\thanks{RICAM Linz, Johann Radon Institute, AT-4040 Linz, 
\texttt{svetlana.matculevich@ricam.oeaw.ac.at}}}
\begin{document}
\title{Fully reliable error control for evolutionary problems}
\maketitle

\begin{abstract}
This work is focused on the application of functional-type a posteriori error 
estimates and corresponding indicators to a class of time-dependent problems. 
We consider the algorithmic part of their derivation and implementation and also
discuss the numerical properties of these bounds that comply with obtained 
numerical results. 
This paper examines two different methods of approximate solution reconstruction 
for evolutionary models, i.e., a time-marching technique and a space-time approach. 
The first part of the study presents an algorithm for global minimisation of the 
majorant on each of discretisation time-cylinders (time-slabs), the effectiveness 
of this algorithm is confirmed by extensive numerical tests. 
In the second part of the publication, the application of functional error estimates 
is discussed with respect to a space-time approach. It is followed by a set of 
extensive numerical tests that demonstrates the efficiency of proposed error control 
method.

\end{abstract}

\input{sections/introduction.tex}

\input{sections/problem-statement.tex}
\input{sections/examples-parabolic-time-stepping.tex}

\input{sections/examples-parabolic-space-time.tex}

\input{sections/conclusions.tex}

%

{\small
\bibliographystyle{plain}
\bibliography{main-func-estim-paper}
}
\end{document}

%% file: customized_packages_and_macros.tex
\usepackage{amsmath}
\usepackage{graphicx}								
\usepackage{subfig}
\usepackage{amsfonts} 							
\usepackage{dsfont}
\usepackage{algorithm}              
\usepackage{algorithmic}
\usepackage{booktabs} 							
\usepackage{multirow}								
\usepackage{setspace}
\usepackage{mathrsfs}
\usepackage{amssymb}
\usepackage{tikz}
\usepackage{color}
\usepackage{xcolor}
\usepackage{mathtools}
\usepackage{xfrac}
\usepackage{upgreek}

\usepackage{framed}
\usepackage{afterpage}
\usepackage{capt-of}
\usepackage{stackengine}

\usepackage{geometry}
\usepackage{color}
\definecolor{deepblue}{rgb}{0,0,0.5}
\definecolor{deepred}{rgb}{0.6,0,0}
\definecolor{deepgreen}{rgb}{0,0.5,0}
\definecolor{gray(x11gray)}{rgb}{0.75, 0.75, 0.75}

\usepackage{tikz}	
\usetikzlibrary{arrows,chains,matrix,positioning,scopes}

\makeatletter
\tikzset{join/.code=\tikzset{after node path={%
\ifx\tikzchainprevious\pgfutil@empty\else(\tikzchainprevious)%
edge[every join]#1(\tikzchaincurrent)\fi}}}

\makeatother
\tikzset{>=stealth',every on chain/.append style={join}, every join/.style={->}}
\tikzstyle{labeled}=[execute at begin node=$\scriptstyle,   execute at end node=$]

\tikzset{
>=stealth',
help lines/.style={dashed, thick},
axis/.style={<->},
important line/.style={thick},
connection/.style={thick, dotted},
}

\newcommand {\Minimum}   {\min\limits}

\newcommand {\R}		 {\mathbf{r}}

\newcommand {\Ieff}  {I_{\rm eff}}

\newcommand {\Rd}    {{\mathds{R}}^d}

\newcommand {\Real}   {{\mathds{R}}}

\newcommand {\Rtwo}   {{\mathds{R}}^2}
\newcommand {\Rthree} {{\mathds{R}}^3}

\def \EnergyNorm#1  {{\mid\!\mid\!\mid #1 \mid\!\mid\!\mid}^2 }   

\def \dvrg       {\mathrm{div}}	

\def \traspose#1 {{#1}^{rm T}}

\newcommand {\vectorv}  {\boldsymbol {v}}

\newcommand {\vectorvarphi}  {\boldsymbol {\varphi}}
\newcommand {\vectoreta}  {\boldsymbol {\eta}}
\newcommand {\vectorp}  {\boldsymbol {p}}
\newcommand {\vectoru}  {\boldsymbol {u}}
\newcommand {\vectorw}  {\boldsymbol {w}}

\newcommand {\vectortau}  {\boldsymbol {\tau}}

\newcommand {\vectorb}  {\boldsymbol {b}}

\newcommand {\flux}     {\boldsymbol {y}}

\def \dt       {\mathrm{\:d}t}
\def \dx       {\mathrm{\:d}x}

\def\L#1{L^{#1}}

\def\Ho#1{H_0^{#1}}

\def\HD#1#2{H^{#1}_{#2}}

\def\error{{[\,e\,]\,}}
\def\mdI{\overline{\mathrm m}_{\mathrm{d}}}
\def\incrmdI#1{\overline{\mathrm m}^{({#1})}_{\mathrm{d}}}
\def\mfI{\overline{\mathrm m}_{\mathrm{eq}}}
\def\incrmfI{\overline{\mathrm m}^{(k)}_{\mathrm{eq}}}

\def\incrmaj{\overline{\mathrm M}^{(k)}} 
\def\maj#1{\overline{\mathrm M}_{#1}}

\def\ed{{e}_{\mathrm{d}}}
\def\et{{e}_{\mathrm{T}}}
\def\errorincr{{[\,e\,]\,}^{(k)}}
\def\incred#1{{e}^{({#1})}_{\mathrm{d}}}

\def\Marker{\mathbb M}

\newcommand {\CFriedrichs} {C_{{\rm F}\Omega}}

\newcommand*\rfrac[2]{{}^{#1}\!/_{#2}}
\newcommand{\minus}{\scalebox{0.5}[1.0]{\( - \)}}

\def\Pone{{\rm P}_{1}}
\def\Ptwo{{\rm P}_{2}}

\def\RTzero{{\rm RT}_{0}}
\def\RTone{{\rm RT}_{1}}

\newtheorem{theorem}{Theorem}{}{}
{}{}
{}{}
\newtheorem{example}{Example}{}{}
\newtheorem{remark}{Remark}{}{}

\def\ProofBegin{\noindent{\bf Proof:} \:}
\def\ProofEnd{{\hfill $\square$}}

\definecolor{formalshade}{rgb}{0.95,0.95,1}

\usepackage{listings} 
\newcommand\pythonstyle{
\lstset{
		language=Python,
		basicstyle=\footnotesize\ttfamily,
		otherkeywords={self},             
		keywordstyle=\bfseries\color{deepblue}\bfseries,
		commentstyle=\itshape\color{purple!40!black},
		stringstyle=\color{deepgreen},
		frame=single,                         
		showstringspaces=false,            %
		belowcaptionskip=.75\baselineskip
}}
\lstnewenvironment{python}[1][]
{
\pythonstyle
\lstset{#1}
}
{}

\newcommand\pythoninline[1]{{\pythonstyle\lstinline!#1!}}

%% file: sections/introduction.tex
\section{Introduction}
\label{sec:introduction}

Evolutionary problems are fundamental components of simulations of 
real-life processes such as heat conduction and thermal radiation models in 
thermodynamics, which are used in modelling aiming
to understand and predict the global climate, and estimation 
of forest growth, among others. 
Most of the models mentioned above are governed by time-dependent {\em partial 
differential equations} (PDEs) or systems of PDEs, which in combination with initial (IC) 
and boundary conditions (BCs) produce so-called {\em initial-boundary value 
problems} (I-BVPs). 
The current study is focused on evolutionary problems of {\em parabolic type}, their
systematic mathematical analysis is presented in  
\cite{Ladyzhenskaya1985, Ladyzhenskayaetall1967, Wloka1987, Zeidler1990A, 
Zeidler1990B}. Their numerical analysis and study of their practical application are 
exposed in \cite{Thomee2006, Lang2001} and partially in classical references on 
the finite element method (FEM) for PDEs (see, e.g., \cite{Braess2001, 
GrossmannRoosStynes2007}). 

In this work we consider functional type error estimates for evolutionary problems. 
Therefore, we let $Q := \Omega \times (0, T)$ denote the space-time cylinder, where 
$\Omega \subset \Rd$, $d \in \{1, 2, 3\}$, is a bounded domain with Lipschitz 
boundary $\partial \Omega$, and $(0, T)$ is a given time interval with $0 < T < +\infty$. 
The surface of the cylinder is divided 
into the initial-time surface $\Sigma_0 := \Omega \times \{0\}$, 
the final-time surface $\Sigma_T := \Omega \times \{T\}$, 
and the remaining part $\Sigma := \partial \Omega \times (0, T)$.
Their union is denoted by $\partial\Omega$.

In this article, we consider the general form of a linear a 
{\em linear parabolic I-BVP problem}, which reads as follows: 
\begin{alignat}{2}
\sigma \, \partial_t u + \mathcal{L} u & = f \;\qquad {\rm in} \quad \; Q , 
\label{eq:equation}\\[-2pt]
u & = u_D \;\;\quad {\rm on} \quad \Sigma, \\[-2pt]
u(x, 0) & = u_0 \quad \; \; \; {\rm on} \quad \Sigma_0.
\label{eq:initial-condition}
\end{alignat}
%
Here, depending on the area of application, $u$ might describe the temperature 
alteration in heat conduction or the concentration of a certain substance in chemical 
diffusion. The given data includes the parameter 
 $\sigma$ (e.g., the conductivity of the material in 
electromagnetics), the source term $f$, the Dirichlet BC $u_D$ (which can be generalised 
to Neumann or Robin conditions), and the initial state $u_0$. The elliptic operator 
$\mathcal{L}$ is written in the general form
$$\mathcal{L} u := 
	-\dvrg_x (A(x) \nabla_x u(x, t)) + \vectorb(x) \cdot \nabla_x u(x, t) + c(x)\, u(x, t),
	\quad (x, t) \in Q,
$$
where $A$ is a material diffusion matrix (possibly anisotropic), and $\vectorb$ 
and $c$ stand for a convection vector-field and a reaction function (possibly anisotropic), 
respectively. 

In most cases of evolutionary systems of PDEs, 
there exists a generalised solution that can be reconstructed by one of the two discretisation 
techniques described below. The first, so-called {\em incremental time-stepping 
method}, includes horizontal and vertical methods of lines (the detailed study of 
this approach can be found in \cite{Thomee2006, Braess2001, 
Johnson2009, Lang2001}). This way of treating evolutionary systems numerically 
is preferred, when the 
implementation of adaptive (in space) solvers is considered. 
In the second approach, time is treated as an additional spatial variable \cite{Hackbusch1984, 
Womble1990, VandewallePiessens1992,HortonVandewalle1995}. This approach is usually referred 
to as the {\em space-time discretization} technique. Unlike the first 
approach, this one does not suffer from the time and space separation on discretisation level
(`curse of sequentiality'), and therefore becomes favourable in parallel computing.
Regardless of the method used, an obtained approximation contains an {\em error}. 
Therefore, it is very important to construct a proper numerical tool to 
analyse obtained results and provide reliable information on the {\em approximation 
error} in them in order to avoid the risk of drawing the wrong conclusion 
obtained from numerical information.

There are two approaches for evaluating the approximation error. An 
\emph{a priori} approach is used for qualitative verification of the theoretical 
properties of a numerical method, e.g., the rate of convergence and asymptotic 
behaviour of the approximation with respect to mesh size parameters (see, 
e.g., \cite{BrennerScott1994,Ciarlet1978,StrangFix1973} and references cited 
therein). In \cite{Thomee2006}, a priori error estimates are 
presented for both semi-discrete problems resulting in a spatial one and  
for most commonly used fully discrete schemes obtained by space-time 
discretization. 
%
This work, however, is focused on the so-called {\em a posteriori} approach, where 
the error is measured after computing the approximation. Unlike a priori error 
analysis, the latter estimates exploit only the given data, e.g., domain 
characteristics, source function together with the IC and BC, and the approximation itself. 
The upper bound of the distance between the approximate and exact solution measured 
in terms of the relevant energy norm is called an {\em error estimate} or {\em majorant}. 
The quantity replicating the distribution of the true error over the domain is called an 
\emph{error indicator}. In particular, the time-marching approach 
produces approximations, which alongside with progress of simulations 
accumulate the error. This error may eventually `blow up' in time if it is not controlled. 
Therefore, appropriate error estimates are crucial for monitoring the error's possible 
dramatic growth especially for non-linear problems. Once the error in the 
approximation is controlled reliably, it is possible to detect the areas with 
excessively large local errors and calculate a much more accurate 
approximation using local refinement.

This work presents numerical properties of \emph{functional type} 
a posteriori error estimates and corresponding indicators, initially introduced in 
\cite{Repin1997, Repin1999, Repin2000} and thoroughly studied for various classes 
of problems in \cite{NeittaanmakiRepin2004, RepinDeGruyter2008, 
Malietall2014} and references therein. Unlike alternative error indicators, e.g., 
gradient averaging indicators \cite{ZienkiewiczZhu1987, ZienkiewiczZhu1988} as well 
as hierarchically based \cite{Deuflhardetall1989} and
goal-oriented estimates \cite{BeckerRannacher1996}, functional type error estimates 
are guaranteed, which means that they always bound the error from above and below.
Moreover, they do not contain mesh-dependent 
local interpolation constants (as residual estimates 
\cite{BabushkaRheinboldt1978, BabushkaRheinboldtSIAM1978}), and are valid for 
any function from the class of conforming approximations. The are not restricted by 
the Galerkin orthogonality assumption. Detailed comparison of the above-described 
approaches can be found in \cite[Section 3.4]{Malietall2014}.
%
In the framework of a posteriori error estimates studied in this work,  
\cite{Repin2002} is the original one, where the method of deriving functional error 
estimates for parabolic I-BVPs was suggested.
The first attempt of their 
numerical analysis was presented in \cite{GaevskayaRepin2005}. 
%
The current work focuses on the practical 
part of the functional estimates application to I-BVPs of parabolic type. 
The thorough theoretical study can be found in \cite{MatculevichRepin2014, 
MatculevichNeitaanmakiRepin2015, MatculevichRepinDiffUr2016,MatculevichRepin2016}.

From the authors point of view, it is important to provide fast, automated, and efficient 
algorithms of reconstructing error estimates and indicators, for 
both time-stepping and space-time approaches that deal with evolutionary equations.
Moreover, it is important to compare the performance of the error estimates for 
both methods and to analyse possible scenarios of the majorant behaviour for  
certain classes of the examples considered. 

In what follows, we present the structure of the paper.
Section \ref{sec:problem-statement} contains the statement of a model problem 
as well as the results on its solvability, which provides the necessary framework for 
subsequent chapters.
We also provide the definition of functional a posteriori error estimates for parabolic I-BVPs. 
In Section \ref{sec:time-stepping}, we discuss the application of the introduced 
error estimates in combination with the time-marching schemes. First, we present 
an algorithm for global minimisation of the majorant on each discretization time-cylinder. 
Numerical tests of this algorithm follow in the same section. 
We apply the same estimates in combination with the space-time approach to 
a serious of examples
and analyse the obtained numerical results in Section \ref{sec:space-time}. We conclude 
that in the cases when time adaptivity is needed, the space-time approach of handling 
the I-BVP combined with suggested error estimates is preferable. It can be observed 
from the performed test-examples, that the majorant provides sharper error estimation
when space-time algorithms are used. However, in engineering applications, 
where the time-incremental analysis is preferable, the functional error estimates also 
provide rather adequate error estimation and indication of its distribution over the 
computional domain.

%% file: sections/problem-statement.tex
\section{Model problem and error estimates}
\label{sec:problem-statement}

In this section, we present a model problem as well as well-posedness 
results for linear parabolic PDEs, which have been thoroughly studied in 
\cite{Ladyzhenskaya1985,Friedman1964,Zeidler1990A,Wloka1987}. We also 
introduce a functional a posteriori error estimate for the stated model and discuss 
its crucial properties. 

Let $Q$ be a space-time cylinder with a boundary surface $\Sigma$ as defined in 
the introduction (see also Figure \ref{eq:space-time-cylinder}). The general parabolic 
I-BVP \eqref{eq:equation}--\eqref{eq:initial-condition} can be re-written as the 
system
\begin{alignat}{3}
	\sigma \, u_t - \dvrg_x \, \vectorp + \vectorb \cdot \nabla_x u + c \, u & =\, f,	      
	& \quad (x, t) \in Q,\label{eq:parabolic-equation}\\
  \vectorp & =\, A \nabla_x u, & \quad (x, t) \in Q,\label{eq:dual-part}\\
  u(x, 0) & =\, u_0,		
	& \;\quad x \in \Sigma_0,
	\\
  u & =\, 0,	& \quad (x, t) \in \Sigma,\label{eq:parabolic-dirichlet-bc}
\end{alignat}
where $u_t := \partial_t u$ is a partial derivative with respect to the time variable, and 
\begin{equation}
f \in \L{2}(Q) \quad  \mbox{and} \quad 
u_0 \in \Ho{1}(\Sigma_0).
\label{eq:problem-condition}
\end{equation}
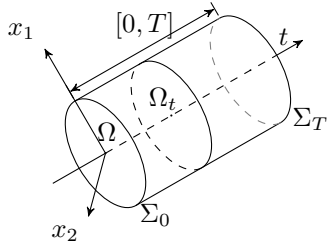
\begin{figure}
\centering
\begin{tikzpicture}[scale=0.8,
comment/.style={rectangle, inner sep= 2pt, text width=6cm, node distance=0.25cm, color=black, font=10pt},
]
\begin{scope}[rotate=30]
\node[color=black, 
transform shape,ellipse,minimum height=2cm,minimum width=1cm,draw,outer sep=0]  (a) {};
\draw (0,1) -- (2.75,1);
\draw (0,-1) -- (2.75,-1);
\draw[dashed,color=gray] (2.75,-1) arc (-90:90:-0.5 and 1);
\draw (2.75,-1) arc (-90:90:0.5 and 1);
\node at (0.2,0.3) {$\Omega$};
\node at (0.2,-1.3) {${\Sigma_0}$};
\node at (3.2,-1.2) {$\Sigma_T$};

\draw[dashed,color=black] (1.25,-1) arc (-90:90:-0.5 and 1);
\draw[color=black] (1.25,-1) arc (-90:90:0.5 and 1);
\node at (1.27,0.3) {$\Omega_t$};
\fill (1.25,0) circle (0.5pt);

\draw (-1,0) -- (-0.0,0.0);
\draw[dashed] (0,0) -- (3.25,0.0);
\draw[->] (3.25,0) -- (3.75,0.0);
\fill (2.75,0) circle (0.4pt);
\fill (0,0) circle (0.5pt);
\draw[<->] (0.01,1.175) -- (2.74,1.175);
\node[above] at (1.03+0.4,1.175+0.02) {$[0, T]$};
\draw (0.0,  1.05) -- (0.0, 1.3);
\draw (2.75, 1.05) -- (2.75,1.3);
\node at (3.5,0.2) {$t$};

\draw[->] (0,0,0) -- (0,2,0) node[anchor=south east]{$x_1$};
\draw[->] (0,0,0) -- (0,0,2) node[anchor=north east]{$x_2$};

\end{scope}
\end{tikzpicture}
\caption{Space-time cylinder $Q$.}
\label{eq:space-time-cylinder}
\end{figure}
%
We assume that $\sigma$ is a positive constant and that the operator $A$ is 
symmetric and satisfies the condition of uniform ellipticity for almost all (a.a.)
$x \in \Omega$, which reads
\begin{equation}
\underline{\nu}_A |\xi|^2 \leq A(x) \: \xi \cdot \xi \leq \overline{\nu}_{A} |\xi|^2, 
\quad \mbox{for} \quad 
\xi \in \Rd,
\quad \mbox{and} \quad 
\quad 0 < \underline{\nu}_A \leq \overline{\nu}_A < \infty.
\label{eq:operator-a}
\end{equation}
We use the notation
\begin{equation*}
  \| \, \vectortau \, \|^2_{A, \Omega} := (A \vectortau, \vectortau)_{\Omega}, \quad 
  \| \, \vectortau \, \|^2_{A^{-1}, \Omega} := (A^{-1} \vectortau, \vectortau)_{\Omega}, 
  \quad 
  \mbox{for all}
  \quad 
  \vectortau \in [\L{2}(\Omega)]^{d}, 
\end{equation*}
where $(\vectoru, \vectorv)_{\Omega} := \int_\Omega \vectoru \cdot \vectorv \dx$ 
(and $(A \vectoru, \vectorv)_{\Omega}$), stands for a (weighted) $\L{2}$ scalar-product
for all $\vectoru, \vectorv \in [\L{2}(\Omega)]^{\rm d}$.
The functions $\vectorb$ and $c$, representing the convection and 
reaction, satisfy the following conditions for a.a.
$t \in (0, T)$
\begin{equation}
\begin{array}{rl}
	\vectorb \in [\L{\infty} (\Omega)]^{d}, \quad 
	\dvrg_x \, \vectorb \in\L{\infty} (\Omega), & \quad |\vectorb| \leq \overline{\vectorb}, \\[4pt]
	c \in \L{\infty} (\Omega), & \quad c \leq \overline{c}, \\[4pt]
	0 < \delta_0 \leq \delta^2 := c - \tfrac{1}{2} \dvrg_x \, \vectorb, &
\end{array}
%
%
%
\label{eq:coefficients-condition}
\end{equation}
where ${\overline{{\vectorb}}}$ and $\overline{c}$ are positive constants.
%
After multiplying (\ref{eq:parabolic-equation}) by a test function 
$$\eta \in \HD{1}{0}(Q) := \big\{ u \in \L{2}(Q) \; \mid \;
\nabla_x u \in \L{2}(Q), \; u_t \in \L{2}(Q) , \; u |_{\Sigma} = 0 \, \big\},$$
we arrive at the generalised formulation of 
(\ref{eq:parabolic-equation})--(\ref{eq:parabolic-dirichlet-bc}): 
find $u \in \HD{1}{0}(Q)$
satisfying the integral identity
\begin{equation}
	(A \nabla_x{u}, \nabla_x{\eta})_Q 
	+ (\vectorb \cdot \nabla_x u, \eta)_Q 
	+ (c u, \eta)_Q 
	-  \sigma (u, \eta_t)_Q
	+ \sigma \big( (u, \eta)_{\Sigma_T} - (u, \eta)_{\Sigma_0} \big) 
	= (f, \eta)_Q, \quad \forall \eta \in \HD{1}{0}(Q).
	\label{eq:generalized-statement}
\end{equation}
According to \cite{Ladyzhenskaya1985}, the generalised problem 
\eqref{eq:generalized-statement} has a unique solution in $\HD{1}{0}(Q)$, 
provided that conditions \eqref{eq:problem-condition}, 
\eqref{eq:operator-a}, and (\ref{eq:coefficients-condition}) hold. 

We present a functional error estimate, which provides a guaranteed upper 
bound of $e := u - v$ for the generalised solution $u$ of the I-BVP 
(\ref{eq:generalized-statement}) and any function $v \in \HD{1}{0}(Q)$. We 
emphasise on the universality of the suggested error estimates, which makes them unique 
in comparison to other existing approaches. The estimates we present in the following are strongly 
independent of the method of the approximation reconstruction. Later on, the considered $v$ 
is generated numerically, and the distance to $u$ is measured in terms of the norm
\begin{equation}
	\! \int_0^T \! \! \big( \nu \left \| \, \nabla_x e \, \right\|^2_{A, \Omega} \,
	+ \, \theta \left\| \, \delta \, e \right\|^2_\Omega \, \big) \dt
	+ \, \zeta \! \, \left \| \, e (\cdot, T) \, \right \|^2_{\Omega} 
	= \nu \left \| \, \nabla_x e \, \right\|^2_{A, Q} \,
	+ \, \theta \left\| \, \delta \, e \right\|^2_Q \, 
	+ \, \zeta \! \, \left \| \, e \, \right \|^2_{\Sigma_T},	
	\label{eq:energy-norm}
\end{equation}
where $\nu$, $\theta$, and $\zeta$ are positive weights and 
the function $\delta$ satisfies (\ref{eq:coefficients-condition}). 
By selecting weights to balance the components of (\ref{eq:energy-norm}) with 
a desired proportion, we generate a collection of error measures, which can be used 
to control $e$.

To derive the upper bounds, we first need to transform \eqref{eq:generalized-statement} 
by subtracting the terms with approximation $v$ 
from left- (LHS) and right-hand side (RHS) and setting $\eta = e$ (see 
details in \cite{MatculevichRepinDiffUr2016}), which implies the following 
{error identity}
\begin{equation}
	\left \| \, \nabla_x e \, \right\|^2_{A, Q} \, + \,
	\left\| \, \delta \, e \right\|^2_Q \, + \,
	\tfrac{\sigma}{2} \! \, \left \| \, e \, \right \|^2_{\Sigma_T} \,
	= \tfrac{\sigma}{2} \! \, \left \| \, e \, \right \|^2_{\Sigma_0}
	+ \big(f - \sigma \, v_t -  c \, v - \vectorb \cdot \nabla_x v, e \big)_Q 
	- \big(A \nabla_x{v}, \nabla_x e \big)_Q.
	\label{eq:energy-balance-equation}
\end{equation}
%
Next, we rearrange the RHS of (\ref{eq:energy-balance-equation}) by introducing a 
`free' vector-valued function 
\begin{alignat*}{2}
	\flux \in H^{\dvrg_x}(Q) := 
	\Big\{  \flux \in \L{2} \big(0, T; [\L{2}(\Omega)\big]^{\rm d} \big)\;\big|\; 
	& \dvrg_x \flux \in \L{2} (Q) \Big \}
\end{alignat*} 	
satisfying $(\dvrg_x \flux, v)_Q + (\nabla_x v, y)_Q = 0$, and, as a result, we arrive at 
\begin{equation}
	\left \| \, \nabla_x e \, \right\|^2_{A, Q} \, + \,
	\left\| \, \delta \, e \right\|^2_Q \, + \,
	\tfrac{\sigma}{2} \! \, \left \| \, e \, \right \|^2_{\Sigma_T} \,
	= \tfrac{\sigma}{2} \! \, \left \| \, e \, \right \|^2_{\Sigma_0}
	+ \big(f + \dvrg_x \flux - \sigma \, v_t -  c \, v - \vectorb \cdot \nabla_x v, e \big)_Q 
	+ \big(\flux - A \nabla_x{v}, \nabla_x e \big)_Q.
	\label{eq:energy-balance-equation-with-flux}
\end{equation}
%
Here, the {\em residuals} correspond to equations \eqref{eq:parabolic-equation} 
and \eqref{eq:dual-part} and are denoted by
\begin{alignat}{2}
	\R_{\rm eq}  (v, \flux) & 
	:= f + \dvrg_x \flux - \sigma \, v_t -  c \, v - \vectorb \cdot \nabla_x v, \label{eq:r-f}\\
	\R_{\rm d}  (v, \flux) & := \flux - A \nabla_x{v}, \label{eq:r-d}
\end{alignat}
respectively. Moreover, we define the weighted residuals
\begin{alignat}{2}
	\R_{\rm eq}^{\mu} (v, \flux) & := \mu \, \R_{\rm eq} \quad{\rm and} \quad
	\R_{\rm eq}^{1 \minus \mu}  (v, \flux) := (1 - \mu) \, \R_{\rm eq},
	\label{eq:r-mu}
\end{alignat}
%
where $\mu(x)$ is a real-valued function taking its values in $[0, 1]$, used in order to 
split the residual $\R_{\rm eq}$ with the reaction and convection parameters into two 
parts. This way, 
the resulting estimate becomes robust even for the cases, when the values of $c$ 
change form the low to high orders of magnitude in different parts of 
$\Omega$.
The detailed 
numerical analysis of the majorant with balancing parameter $\mu$ can be found 
in \cite[Sections 2, 5]{MatculevichRepin2014}.
The theorem below recalls that a certain weighted combination of norms 
measuring the residuals \eqref{eq:r-f}--\eqref{eq:r-mu} bounds the error $e$. 

\begin{theorem}
\label{th:theorem-minimum-of-majorant-I}
(a) For any functions $v \in \HD{1, 1}{0}(Q)$ and $\flux \in H^{\dvrg_x}(Q)$, 
parameters $\nu \in (0, 2]$ and $\gamma \in \big[\tfrac12, +\infty\big[$, 
we have the estimate
\begin{multline}
	(2 - \nu) \, \left \| \, \nabla_x e \, \right\|^2_{A, Q} \, + \,
	\left\| \, \sqrt{2 - \tfrac{1}{\gamma}} \, \delta \, e \right\|^2_Q \, + \,
	\sigma \, \left \| \, e \, \right \|^2_{\Sigma_T} \\
	\leq \maj{} (v, \flux; \nu, \gamma, \mu, \alpha_i) \! 
	:= \sigma \, \| e \|^2_{\Sigma_0}
	+ \int_0^T \! \! \Big ( 
	\gamma \,\left\|\,\tfrac{1}{\delta}\,\R_{\rm eq}^{\mu}\,\right\|^2_{\Omega} \,	
	+ \alpha_1\, \|\,\R_{\rm d} \,\|^2_{A^{- 1}, \Omega} 
	+ \alpha_2 \,\tfrac{\CFriedrichs^2}{\,\underline{\nu}_A}
		\big\| \, \R_{\rm eq}^{1 \minus \mu} \, \big\|^2_{\Omega} \, 
	\Big ) \dt,
	\label{eq:majorant-1}
\end{multline}
%
where $\mu(x) \in [0, 1]$ is a real-valued function, and 
$\alpha_i$, $i = 1, 2$ are positive parameters satisfying the 
relation $\sum\limits_{i = 1}^{2} \tfrac{1}{\alpha_i} = \nu$.
Here, $\CFriedrichs$ is the constant in the Friedrichs inequality \cite{Friedrichs1937} 
\begin{equation*}
	\| v \|_{\Omega} \leq \CFriedrichs \| \nabla_x v \|_{\Omega},
	\quad \forall v \in \HD{1}{0}(\Omega) 
	:= \big\{ \, v \in \L{2}(\Omega) \; \mid \;
	\nabla_x v \in \L{2}(\Omega) , \; v |_{\partial \Omega} = 0 \, \big\},
	\label{eq:friedrichs-inequality}
\end{equation*}
and the residuals $\R_{\rm d}$, $\R_{\rm eq}^{\mu}$, and 
$\R_{\rm eq}^{1 \minus \mu}$ are defined in \eqref{eq:r-d} and \eqref{eq:r-mu}, 
respectively. 

\noindent
(b) For any parameters $\nu \in (0, 2]$, $\gamma \in \big[\tfrac{1}{2}, + \infty)$, 
$\alpha_i \in (0, + \infty)$, and any real-valued function $\mu(x) \in [0, 1]$, the 
variational problem
\begin{equation*}
\inf\limits_{
\scriptsize
\begin{array}{c}
v\in \HD{1}{0}(Q)\\[1pt]
\flux \in H^{\dvrg_x}(Q)
\end{array}
}
\maj{} (v, \flux) 
\end{equation*}
has a solution (with the corresponding zero-value for the functional), and 
its minimum is attained if and only if $v = u$ and $\flux = A \nabla_x u$.
\end{theorem}
\ProofBegin The detailed proof of the theorem can be found in, e.g., 
\cite[Theorem 1]{MatculevichRepinDiffUr2016}.
\ProofEnd
\begin{remark}
\rm
Unlike analogous error estimates for elliptic problems, the majorant is not sharp 
with respect to the error measured in the energy norm. Let parameters $\nu = 1$, 
$\mu = 0$, $\gamma = 1$,  and let the approximate solution $v$ satisfy the initial condition exactly, i.e., 
$e(\cdot, 0) = 0$. For the flux chosen as $\flux = A \nabla_x u$, 
we obtain 
\begin{alignat}{2}
	\maj{} (v, A \nabla_x u; \nu, 1, 0, \alpha_i) \! 
	& = \int_0^T \! \! \Big ( 
	\alpha_1\, \|\, \nabla_x e \,\|^2_{A^{- 1}, \Omega} 
	+ \alpha_2 \,\tfrac{\CFriedrichs^2}{\,\underline{\nu}_A}
	\big\| \, f + \dvrg_x  (A \nabla_x u) - \sigma \, v_t  - c \, v - \vectorb \cdot \nabla_x v\, \big\|^2_{\Omega} \, 
	\Big ) \dt \\
	& = \int_0^T \! \! \Big ( 
	\alpha_1\, \|\, \nabla_x e \,\|^2_{A^{- 1}, \Omega} 
	+ \alpha_2 \,\tfrac{\CFriedrichs^2}{\,\underline{\nu}_A}
	\big\| \, \sigma \, e_t + c \, e + \vectorb \cdot \nabla_x e\, \big\|^2_{\Omega} \, 
	\Big ) \dt.
\label{eq:gap}
\end{alignat}
The second term of the obtained functional contains $e_t$ and the 
constant $\tfrac{\CFriedrichs^2}{\,\underline{\nu}_A}$, which are not included 
in \eqref{eq:energy-norm} and contribute to the irremovable gap between the error and 
the estimate. 
\end{remark}


\begin{remark}
\rm
The same form of the majorant was presented and numerically tested for the evolutionary 
reaction-diffusion I-BVPs of parabolic type in \cite{MatculevichRepin2014}. The current 
work is based on implementation of the I-BVP and the majorant in Python using 
The FEniCS Project library \cite{LoggMardalWells2012} and 
contains a detailed explanation of numerical results clarifying the properties of the
functional approach to fully reliable computations. We note that the 
implementation is based on the version of the library that does not support the 
coarsening of the mesh in the time-stepping approach,
therefore it is not discussed in the paper.

\end{remark}

%% file: sections/examples-parabolic-time-stepping.tex
\section{Time-stepping approach}
\label{sec:time-stepping}

In this section, we discuss the application of the majorant for the case of a time-marching 
approach and steps of its global minimisation on each of the incremental time-cylinders. 
The summary of this scheme is contained in Algorithm \ref{alg:majorant-min}, and its 
efficiency is confirmed by numerical results presented in Examples 
\ref{ex:incr-unit-2d-t}--\ref{ex:incr-pi-shape-2d-t}.

\subsection{Global minimisation of the increment of the majorant}

For the reader's convenience, we assume that $\vectorb = {\boldsymbol 0}$, $c = 0$ 
(which implies that $\mu = 0$ and $\gamma = 1$), $\sigma = 1$, 
and $v(\cdot, 0) = u_0$. Using this rather simple formulation of the model, 
one can capture the main idea of the presented algorithm and key 
numerical properties of the error majorant. Thus, the error is defined as follows
\begin{equation}
\error := (2 - \nu)\, \ed + \et,  \quad \nu \in (0, 2], \quad \mbox{where} \quad
\ed = \int_0^T \| \nabla_x e\|^2_{A, \Omega} \dt = \| \nabla_x e\|^2_{A, Q} 
\quad \mbox{and} \quad 
 \et := \| e\|^2_{\Sigma_T}.
\label{eq:error-simplified-total-time}
\end{equation}
In the case of $\nu = 2$, it is reduced to the error measured at $\Sigma_T$, which
can be controlled by the error estimate presented in Theorem 
\ref{th:theorem-minimum-of-majorant-I}. The majorant reads as
\begin{equation}
\maj{} (v, \flux; \alpha_1, \alpha_2) 
	:= \alpha_1 \int_0^T \! \| \flux - A \nabla_x v \,\|^2_{A^{-1}, \Omega} \dt + 
		\alpha_2 \tfrac{\CFriedrichs^2}{\underline{\nu}_A}
		\!\! \int_0^T  \!\! \| \, f + \dvrg_x \flux - v_t  \, \|^2_{\Omega} \dt.
\label{eq:majorant-simplified-total-time}
\end{equation}
%
Here, the second term
$$\mfI 
:= \int_0^T \| \R_{\rm eq} \, \|^2_{\Omega}\dt 
= \int_0^T \| \, f + \dvrg_x \, \flux - v_t  \, \|^2_{\Omega}\dt$$ 
assures the 
reliability of the majorant and measures the violation of equation
\eqref{eq:parabolic-equation}, whereas the first term 
$$\mdI := 
\int_0^T \| \, \R_{\rm d} \,\|^2_{A^{-1}, \Omega} \dt = 
\int_0^T \| \, \flux - A \nabla_x v \,\|^2_{A^{-1}, \Omega} \dt$$ 
mimics the 
residual in \eqref{eq:dual-part} and can be used as a robust and efficient 
indicator. The reliability and accuracy of $\maj{}$ is measured by
the so-called efficiency index 
$\Ieff := \sqrt{\tfrac{{\overline{\mathrm M}}}{[e]}}$.

In order to adapt the majorant \eqref{eq:majorant-simplified-total-time} 
to the time-stepping approach of the approximation reconstruction, we 
introduce the following discretisation of the time-interval $[0,\; T]$: 
\begin{equation}
\mathcal{T}_K = \cup_{k = 0}^{K-1} \overline{{I}^{(k)}}, \quad {\rm where} 
\quad 
{I}^{(k)} = (t^k, t^{k + 1}),
\end{equation}
where $K$ corresponds to the number of sub-intervals. 
Generally, the spatial domain $\Omega_t$ can change its shape over time, 
i.e., $Q := \{ (x, t) : x(t) \in \Omega_t, \, t \in (0, T)\}$. In this case, 
the space-time FEM approach is rather logical. However, in the scope of this paper, 
we consider only problems on the non-moving spatial domain in time (non-moving 
space-time cylinders).
Then, representation of the space-time cylinder can be defined as follows:
\begin{equation}
\overline{Q} = \cup_{k = 0}^{K-1} \overline{Q^{(k)}}, \quad \mbox{where} \quad
Q^{(k)} := {I}^{(k)} \times \Omega.
\end{equation}
We emphasise that throughout the paper $Q^{(k)}$ are referred to as 
time-cylinders, time-slabs, or time-slices (that have volume with respect to (w.r.t.) time), 
whereas $\Omega^{(k)} := \Omega \times t^k$ are understood as instant time-cuts.
Let $\mathcal{T}_{N}$ be a mesh selected on $\Omega$, where $N$ is the number 
of elements in the space discretisation. Then, $\Theta_{K \times N} = 
\mathcal{T}_K \times \mathcal{T}_{N}$ denotes the mesh on $Q$. 
%

From now on, we assume that the approximate solution is reconstructed on a 
particular time-slice $Q^{(k)}$, such 
that $\flux \in H^{\dvrg_x} (Q^{(k)})$ and
$v \in \HD{1}{0}(Q^{(k)})$. We set $\alpha_1 = \tfrac{1}{\nu} (1 + \beta)$ 
and $\alpha_2 = \tfrac{1}{\nu} (1 + \tfrac{1}{\beta})$, where $\beta$ is a 
positive parameter. On each $Q^{(k)}$, {\em the increment of majorant} 
\eqref{eq:majorant-simplified-total-time} is denoted by $\incrmaj{}$, i.e., 
\begin{equation}
	\incrmaj (v, \flux; \beta) 
	:= \tfrac{1}{\nu}\Big((1 + \beta) \, \incrmdI{k} + 
	\big ( 1 + \tfrac{1}{\beta} \big ) \tfrac{\CFriedrichs^2}{\underline{\nu}_A}  \,
	\incrmfI \Big).
	\label{eq:majorant-simplified}
\end{equation}
%
We define the optimal $\flux_{\rm min}$ by minimisation of this increment, 
i.e.,
\begin{equation*}
\flux_{\rm min} := {\rm arg}
\Minimum_{\flux \in H^{\dvrg_x} (Q^{(k)})} 
\incrmaj (v, \flux; \beta).
\end{equation*}
%
The corresponding {\em increment of the error} is denoted by $\errorincr$.
The minimum of $\incrmaj (\flux; \beta) $ w.r.t. $\beta$ is attained at 
%
$\beta_{\min} := 
\left( \tfrac{\CFriedrichs^2 \incrmfI }{ 
\underline{\nu}_A \incrmdI{k}} \right)^{\rfrac{1}{2}}$.
%
After $\beta$ is fixed, the necessary condition for the minimiser $\flux$ reads as 
\begin{equation}
\tfrac{{\rm d}\incrmaj (v, \, \flux + \zeta \vectorw; \, \beta)}{{\rm d} \zeta} 
\Big|_{\zeta = 0} = 0, 
\label{eq:majorant-derivative-equal-to-zero}
\end{equation}
where $\vectorw \in H^{\dvrg_x} (Q^{(k)})$.
%
Condition \eqref{eq:majorant-derivative-equal-to-zero} implies that
\begin{equation}
	\tfrac{\CFriedrichs^2}{\beta \, \underline{\nu}_A} \,
	(\dvrg_x \flux, \dvrg_x \vectorw)_{Q^{(k)}} 
         +  (A^{-1}\flux, \vectorw)_{Q^{(k)}}
	= - \tfrac{\CFriedrichs^2}{\beta \, \underline{\nu}_A} 
	   \big((f - v_t), \dvrg_x \vectorw)_{Q^{(k)}}
	   + (\nabla_x v, \vectorw)_{Q^{(k)}}.
\label{eq:majorant-derivative}
\end{equation}
%
%
\begin{figure}
\centering
\begin{tikzpicture}[scale=0.8,
comment/.style={rectangle, inner sep= 2pt, text width=6cm, node distance=0.25cm, font=\rmfamily, color=black},
]
\begin{scope}[rotate=30]
\node[color=black, 
transform shape,ellipse,minimum height=2cm,minimum width=1cm,draw,outer sep=0]  (a) {};
\draw (0,1) -- (4.75,1);
\draw (0,-1) -- (4.75,-1);
\draw[dashed,color=gray] (4.75,-1) arc (-90:90:-0.5 and 1);
\draw (4.75,-1) arc (-90:90:0.5 and 1);
\node at (0.2,0.3) {$\Sigma_0$};
\node at (0.2,-1.3) {$u_0$};
\node[above] at (0.03+0.04,1.175+0.02) {$t^0$};

\draw[dashed,color=black] (1.25,-1) arc (-90:90:-0.5 and 1);
\draw[color=black] (1.25,-1) arc (-90:90:0.5 and 1);
\node at (1.27,0.3) {$\Omega_{t^{k-1}}$};
\node at (1.5,-1.3) {$v^{k-1}$};
\node[above] at (1.03+0.4,1.175+0.02) {$t^{k-1}$};
\draw (1.25, 1.05) -- (1.25,1.3);
\fill (1.25,0) circle (0.5pt);

\draw[dashed,color=black] (2.25,-1) arc (-90:90:-0.5 and 1);
\draw[color=black] (2.25,-1) arc (-90:90:0.5 and 1);
\node at (2.27,0.3) {$\Omega_{t^{k}}$};
\node at (2.5,-1.3) {$v^{k}$};
\node[above] at (2.03+0.4,1.175+0.02) {$t^{k}$};
\draw (2.25, 1.05) -- (2.25,1.3);
\fill (2.25,0) circle (0.5pt);

\draw[dashed,color=black] (3.25,-1) arc (-90:90:-0.5 and 1);
\draw[color=black] (3.25,-1) arc (-90:90:0.5 and 1);
\node at (3.37,0.3) {$\Omega_{t^{k+1}}$};
\node at (3.5,-1.3) {$v^{k+1}$};
\node[above] at (3.03+0.4,1.175+0.02) {$t^{k+1}$};
\draw (3.25, 1.05) -- (3.25,1.3);
\fill (3.25,0) circle (0.5pt);

\draw (-1,0) -- (-0.0,0.0);
\draw[dashed] (0,0) -- (5.25,0.0);
\draw[->] (5.25,0) -- (5.75,0.0);
\fill (4.75,0) circle (0.4pt);
\fill (0,0) circle (0.5pt);
\draw[<->] (0.01,1.175) -- (4.74,1.175);
\node at (5.27,1.175+0.02) {$\Sigma_T$};
\draw (0.0,  1.05) -- (0.0, 1.3);
\draw (4.75, 1.05) -- (4.75,1.3);
\node at (5.5,0.2) {$t$};

\draw[->] (0,0,0) -- (0,2,0) node[anchor=south east]{$x_1$};
\draw[->] (0,0,0) -- (0,0,2) node[anchor=north east]{$x_2$};
\end{scope}
\end{tikzpicture}
\caption{Incremental approach to the approximation reconstruction.}
\end{figure}
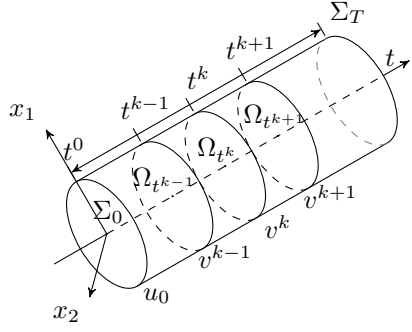
%
We reduce the integration w.r.t time by the following linear 
approximation of $v$ and $\flux$ on the increment $Q^{(k)}$:
\begin{equation}
v = v^k \tfrac{t^{k+1} - t}{\tau^k} + v^{k + 1} \tfrac{t - t^k}{\tau^k}, \quad
\flux = \flux^k \tfrac{t^{k+1} - t}{\tau^k} + \flux^{k + 1} \tfrac{t - t^k}{\tau^k}, 
\quad \tau^k = t^{k + 1} - t^k,
\label{eq:linerization}
\end{equation}
such that $v^k$, $v^{k+1}$ $\in \HD{1}{0}(\Omega)$, $\flux^k$, $\flux^{k + 1}$ 
$\in H^{\dvrg_x}(\Omega)$, and $\vectorw (x, t) := \vectoreta(x) \, T(t)$ with 
$T(t) = \tfrac{t - t^k}{\tau^k}$ and \linebreak $\vectoreta \in H^{\dvrg_x}(\Omega)$. By using 
substitution of \eqref{eq:linerization} and $\vectorw (x, t)$ into 
\eqref{eq:majorant-derivative}, we arrive at 
\begin{multline*}
\tfrac{\CFriedrichs^2}{\underline{\nu}_A}
\bigg(\Big( \dvrg_x \flux^k \tfrac{t^{k+1} - t}{\tau^k} + \dvrg_x \flux^{k+1} \tfrac{t - t^k}{\tau^k} \Big), 
\dvrg_x \vectoreta \, \tfrac{t - t^k}{\tau^k} \bigg)_{Q^{(k)}}
+ \beta \bigg(A^{-1} \Big( \flux^k \tfrac{t^{k+1} - t}{\tau^k} + \flux^{k + 1} \tfrac{t - t^k}{\tau^k}
 \Big), \vectoreta \, \tfrac{t - t^k}{\tau^k}\bigg)_{Q^{(k)}} \\
= - \tfrac{\CFriedrichs^2}{\underline{\nu}_A}  
   \bigg(\Big( f - \tfrac{ v^{k+1} -  v^k}{\tau^k} \Big), 
            \dvrg_x \vectoreta \, \tfrac{t - t^k}{\tau^k}\bigg)_{Q^{(k)}}
+ \beta 
\bigg(\Big( \nabla_x v^k \tfrac{t^{k+1} - t}{\tau^k} 
+ \nabla_x v^{k + 1} \tfrac{t - t^k}{\tau^k}\Big), 
         \vectoreta \, \tfrac{t - t^k}{\tau^k}\bigg)_{Q^{(k)}}.
\end{multline*}
%
Taking into account integration results 
%
$$\int_{t^k}^{t^{k + 1}} (t - t^k) \dt
= \int_{t^k}^{t^{k + 1}}  (t^{k + 1} - t) \dt = \tfrac{(\tau^k)^2}{2}, \quad
\int_{t^k}^{t^{k + 1}}  (t - t^k)(t^{k + 1} - t) \dt = \tfrac{(\tau^k)^3}{6}, \quad
\int_{t^k}^{t^{k + 1}}  \tfrac{1}{\tau^k}(t - t^k)^2 \dt = \tfrac{(\tau^k)^2}{3},
$$
we obtain the following identity
\begin{multline}
\tfrac{\CFriedrichs^2}{\underline{\nu}_A} \tfrac{\tau^k}{3}
\Big(\big( \dvrg_x {\flux^{k+1}} + \tfrac12 \, \dvrg_x \flux^k \big),  \dvrg_x \vectoreta\Big)_\Omega
+ \tfrac{\beta \, \tau^k}{3} 
  (A^{-1} \big( \mathbf{\flux^{k+1}} + \tfrac12 \, \flux^k \big), \vectoreta)_\Omega \\
= - \tfrac{\CFriedrichs^2}{\underline{\nu}_A} 
   \bigg( \Big(\tfrac{1}{\tau^k} F_{(t - t^k)}(x) - \tfrac{v^{k+1} - v^k}{2} \Big), \dvrg_x \vectoreta\bigg)_{\Omega} 
   + \tfrac{\beta \, \tau^k}{3} \Big(\big(\nabla_x v^{k+1} + \tfrac12 \, \nabla_x v^k \big), \vectoreta \Big)_\Omega,
\label{eq:majorant-derivative-linerized-3}
\end{multline}
where 
$$F_{(t \minus t^k)} (x) := \int_{t^k}^{t^{k+1}} f (x, t) \, (t - t^k) \dt.$$
After multiplying (\ref{eq:majorant-derivative-linerized-3}) by 
$\tfrac{3}{\beta\, \tau^k}$, we have 
\begin{multline}
	\tfrac{\CFriedrichs^2}{\beta\, \underline{\nu}_A}
	\Big( \big( \dvrg_x {\flux}^{k+1} + \tfrac12 \, \dvrg_x \flux^k \big), \dvrg_x \vectoreta \Big)_\Omega
	+ (A^{-1} \big({\flux}^{k+1} + \tfrac12 \, \flux^k \big), \vectoreta\Big)_\Omega \\
	= - \tfrac{\CFriedrichs^2}{\beta\, \underline{\nu}_A} 
            \bigg(3\, \Big( \tfrac{1}{(\tau^k)^2} F_{(t - t^k)}(x) - \tfrac{v^{k+1} - v^k}{2 \tau^k} \Big), \dvrg_x \vectoreta\bigg)_\Omega
	  + \bigg( \big(\nabla_x v^{k+1} + \tfrac12 \, \nabla_x v^k \big), \vectoreta\bigg)_\Omega,	
	\label{eq:majorant-derivative-linerized-4}
\end{multline}
where ${\flux^{k+1}}$ is the unknown function that we are aiming to optimise, 
and 
$F_{(t \minus t^k)} (x)$ is approximated by the Gauss quadratures of $4$-th order 
(see, e.g., \cite{StroudSecrest1966, Stroud1974}).

Let us assume that $\flux^{k}, {\flux^{k+1}}$, and $\vectoreta \in {\rm span} \, 
\big\{ \, \vectorvarphi_1, ..., \vectorvarphi_M \big\} \subset H^{\dvrg_x}(\Omega)$, i.e., 
$\flux^{k} := \sum_{i = 1}^M Y^{k}_i \vectorvarphi_i$, where 
$Y^{k} \in {\mathds R}^M$ is a vector of degrees of freedom (DOFs) 
for $\flux^k$,  and $\vectoreta$ can be chosen as 
$\vectoreta := \vectorvarphi_j$, $j = 1, \ldots, M$.
%
Then, condition \eqref{eq:majorant-derivative-equal-to-zero} implies the system 
of linear algebraic equations 
\begin{equation}
\left( \tfrac{\CFriedrichs^2}{\beta \, \underline{\nu}_A} S + K \right) {Y^{k+1}} = 
- \tfrac12 \left( \tfrac{\CFriedrichs^2}{\beta \,\underline{\nu}_A} S + K \right) Y^k
- \tfrac{\CFriedrichs^2}{\beta \, \underline{\nu}_A} \tfrac{3}{(\tau^k)^2} z + g,
\label{eq:system}
\end{equation}
where the components of the matrices $S$ and $K \in {\mathds R}^{M \times M}$ and 
the vectors $z$ and $g \in {\mathds R}^{M}$ are defined as follows:
\begin{alignat}{2} \!\!
\{ S_{ij} \}_{i, j=1}^M  & := 
(\dvrg_x \vectorvarphi_i, \dvrg_x \vectorvarphi_j)_{\Omega}, 
\label{eq:div-div}\\
\{ K_{ij} \}_{i, j=1}^M  & := 
(A^{-1} \vectorvarphi_i, \vectorvarphi_j)_{\Omega}, 
\label{eq:vector-mass}\\
\{ z_{j} \}_{j=1}^M & := 
\Big(\big( F_{(t - t^k)} + \tfrac{( v^k - v^{k+1}) \tau^k }{2} \big),  
\dvrg_x \vectorvarphi_j \Big)_\Omega,  
\label{eq:rhs-div-div}\\
\{ g_{j} \}_{j=1}^M & := 
\Big(\big( \nabla_x v^{k+1} + \tfrac{1}{2} \, \nabla_x v^k \big), \vectorvarphi_j\Big)_\Omega.
\label{eq:rhs-mass}
\end{alignat}
%
The sequence of the above-listed arguments is summarised in Algorithm \ref{alg:majorant-min}, 
which presents the procedure of reconstructing the optimal $\flux^{k+1}$ and the
corresponding $\incrmaj (v, \flux; \beta)$, 
such that on each time-step the increment of the majorant is approximated by means 
of the iteration procedure. 
The sequence of fluxes, produced by each iterative loop, helps to generate the sequence 
of the optimal upper bounds that approach the error as close as possible. 
Each 
$\flux^{k + 1}$ reconstructed on $Q^{(k)}$ is used as the initial data on the time-slab 
$Q^{(k+1)}$. It is important to note that the matrices $S$, $K$ and the vectors 
$z$, $g$ need to be assembled only once since they do not change in the  
minimisation cycle. Generally, Algorithm \ref{alg:majorant-min} can be extended to  
approximations that have jumps w.r.t. the time variable 
(see, e.g., \cite{RepinTomar2010}). 
Moreover, the upper bound can be used as a refinement criterion for schemes 
adaptive in time. 
\begin{algorithm}[!t]
\caption{\quad Global minimisation of $\incrmaj$ (in the case of a time-stepping scheme)}
\label{alg:majorant-min}
\begin{algorithmic} 
\STATE {\bf Input:} $Q^{(k)}$: $v^k, v^{k+1}, \flux^k$ 
\COMMENT{approximations at fixed moments of time and 
flux coefficients on $\Omega \times t^k$}
\STATE $\quad \qquad$ $\vectorvarphi_i$, $i = 1, \ldots, M$ 
\COMMENT{$H^{\dvrg_x}(\Omega)$-conforming basis functions}
\STATE $\quad \qquad$ $L^{\rm iter}_{\rm max}$ \COMMENT{number of inner optimisation iterations}
\STATE
\STATE Assemble the matrices $S$, $K$ $\in {\mathds R}^{M \times M}$and the vectors $z$, $g$ $\in {\mathds R}^{M}$ by using 
\begin{alignat*}{2}
\{ S_{ij} \}_{i, j=1}^M  & = 
(\dvrg_x \vectorvarphi_i, \dvrg_x \vectorvarphi_j)_\Omega, \; \quad
\{ z_{j} \}_{j=1}^M      = 
\Big(\big( F_{(t - t^k)} + \tfrac{( v^k - v^{k+1}) \tau^k }{2} \big), \, 
\dvrg_x \vectorvarphi_j \Big)_\Omega, \\
\{ K_{ij} \}_{i, j=1}^M  & = 
\Big(A^{-1} \vectorvarphi_i, \vectorvarphi_j \Big)_\Omega, \qquad \quad
\{ g_{j} \}_{j=1}^M      = 
\Big(\big( \nabla_x v^{k+1} + \tfrac{1}{2} \, \nabla_x v^k \big), 
\vectorvarphi_j \Big)_\Omega.
\end{alignat*}
\STATE Approximate the flux $\flux^{k} = \sum\limits_{i=1}^M Y^{k}_i \vectorvarphi_i$.
\STATE Initialise $\beta$, e.g., $ \beta = 1$.
\vspace{4pt}
\FOR{$l = 1$ {\bf to} $L^{\rm iter}_{\rm max}$}
\vspace{4pt}
\STATE Solve the system 
$
\left( \tfrac{\CFriedrichs^2}{\beta \, \underline{\nu}_A} S + K \right){Y^{k+1}} = 
	- \tfrac12 \left( \tfrac{\CFriedrichs^2}{\beta \,\underline{\nu}_A} S + K \right) Y^k
	- \tfrac{\CFriedrichs^2}{\beta \, \underline{\nu}_A} \tfrac{3}{(\tau^k)^2} z + g.
$
\STATE Approximate the flux $\flux^{k+1} = \sum\limits_{i=1}^M Y^{k+1}_i \vectorvarphi_i$.
\STATE Reconstruct $v$ and $\flux$ on $Q^{(k)}$ by using 
\begin{equation*}
v = v^k \tfrac{t^{k+1} - t}{\tau^k} + v^{k + 1} \tfrac{t - t^k}{\tau^k}, \quad
\flux = \flux^k \tfrac{t^{k+1} - t}{\tau^k} + \flux^{k + 1} \tfrac{t - t^k}{\tau^k}, \quad 
\tau^k = t^{k + 1} - t^k.
\end{equation*}

\vspace{3pt}
\STATE Compute the components of the majorant 
by using 
$$
\incrmdI{k} = \| \, \flux - A \nabla_x v \, \|^2_{A^{-1}, Q^{(k)}}
\quad {\rm and} \quad 
\incrmfI = \| \, f + \dvrg_x \flux - v_t  \, \|^2_{Q^{(k)}}.
$$
\vspace{-5pt}
\STATE Compute the optimal $\beta$ by using
$
\beta = 
\left( \tfrac{\CFriedrichs^2 \incrmfI }{ \underline{\nu}_A \incrmdI{k}} \right)^{\rfrac{1}{2}}.
$
\ENDFOR
\vspace{5pt}
\STATE Compute the increment of the majorant by using 
%
$\incrmaj (v, \flux; \beta) = (1 + \beta) \, \incrmdI{k} + 
    (1 + \tfrac{1}{\beta}) \tfrac{\CFriedrichs^2}{\underline{\nu}_A} \, \incrmfI$.
%
\STATE {\bf Output:} $\incrmaj (v, \flux; \beta)$ 
\COMMENT{increment of the majorant on $Q^{(k)}$}
\STATE $\qquad \qquad$ $\flux^{k+1}$ 
\COMMENT{reconstruction of the flux on $t^{k + 1} \times \Omega$}
\end{algorithmic}
\end{algorithm}

\newpage
A detailed study of the numerical application of $\maj{} (v, \flux)$ for 
$\Omega \in \Rd$, $d = \{1, 2\}$ is performed in 
\cite{MatculevichRepin2014}. In particular, the robustness of the majorant is tested w.r.t. 
different reaction functions $c$. Besides that, the numerical properties of 
the originally introduced minorant are tested and compared to the majorant. In 
Examples 3--5 of the same work
\cite{MatculevichRepin2014}, the numerical behaviour of the indicators $\incrmdI{k}$ 
is studied in details, i.e., their efficiency is verified by several criteria based on 
different marking procedures (denoted by $\Marker$), quantitative histograms and 
other means. The minimisation of the majorant in \cite{MatculevichRepin2014} was 
based on localised minimisation of the majorant, which performs considerably slower 
than the optimisation technique presented in Algorithm \ref{alg:majorant-min}.

\subsection{Numerical examples}

The current section is dedicated to the numerical examples, in which the I-BVP is discretised by 
the incremental method and the majorant is reconstructed and optimised using the 
global optimisation strategy presented by Algorithm \ref{alg:majorant-min}. From 
here on, the parameter $\nu$ in \eqref{eq:majorant-simplified} is set to $1$. We 
start from relatively basic test-problems, in order to clarify numerical behaviour of the 
majorant to the reader, and add several complications in examples towards the end 
of this section. Moreover, due to restrictions related to the graphical representation, 
only meshes of maximum $17 \cdot 10^3$ elements (EL) are presented in the examples, 
even those all problems have been tested on the meshes with up to $10^{6}$ elements. 
	
\begin{example}
\label{ex:incr-unit-2d-t}
\rm
First, we consider a benchmark problem on a unit square domain
$\Omega = (0, 1)^2 \subset \Rtwo$ with $T = 1$. We choose homogeneous 
Dirichlet BC, $A$ as a unit matrix, $\vectorb = {\boldsymbol 0}$, $c = 0$,
initial state $u_0 = x\,(1 - x)\,y\,(1 - y)$, 
and $u = x\,(1 - x)\,y\,(1 - y)\,(t^2 + t + 1)$ as the exact solution (the source 
function $f$ is calculated respectively). On this relatively simple problem, we aim
to highlight the most important numerical properties of the majorant and error indicator. 
Besides that, the reader gets the chance to become familiar with 
the systematic structure of the computational results analysis, which we follow in the other 
examples throughout the paper.

The function $v$ is reconstructed by the 
Lagrangian finite element space of order one denoted $\Pone$, and $\flux$ is 
approximated by linear Raviart-Thomas finite elements (FEs), which we reffer further as 
$\RTone$ .  
The optimal convergence test for fixed number of time-steps $K=100$ 
and a decreasing mesh size $h$ is illustrated in Figure 
\ref{fig:example-incr-unit-2d-uniform-convergence}. The time discretisation step is 
chosen small enough, in order to minimise its effect on the order of error convergence 
(o.e.c.) w.r.t. the refinement in space. Here, Figure 
\ref{fig:example-incr-unit-2d-uniform-convergence-a} depicts the total error 
$\error$ as defined in \eqref{eq:error-simplified-total-time} and the majorant $\maj{}$, 
whereas Figure \ref{fig:example-incr-unit-2d-uniform-convergence-b} illustrates the 
dominating term of the true error $\ed = \| \nabla_x e\|^2_{Q}$ and the indicator 
$\mdI$, which have different magnitude but decrease w.r.t $h$ with expected 
convergence order $O(h^2)$. 

\begin{figure}[!ht]
	\centering
	\subfloat[$\error$ and $\maj{}$]{
	\includegraphics[width=6cm]{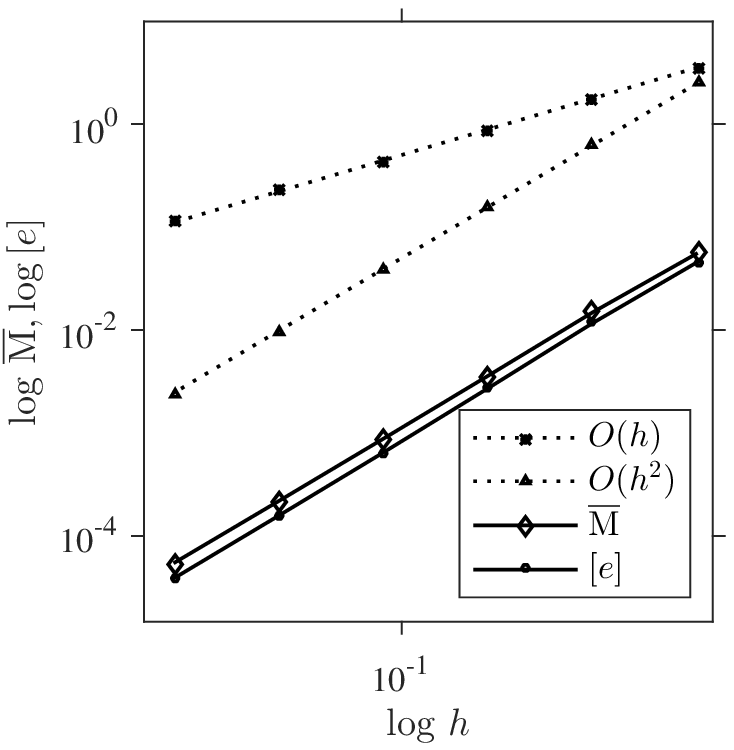}
	\label{fig:example-incr-unit-2d-uniform-convergence-a}}\quad
	\subfloat[$\ed$ and $\mdI$]{
	\includegraphics[width=6cm]{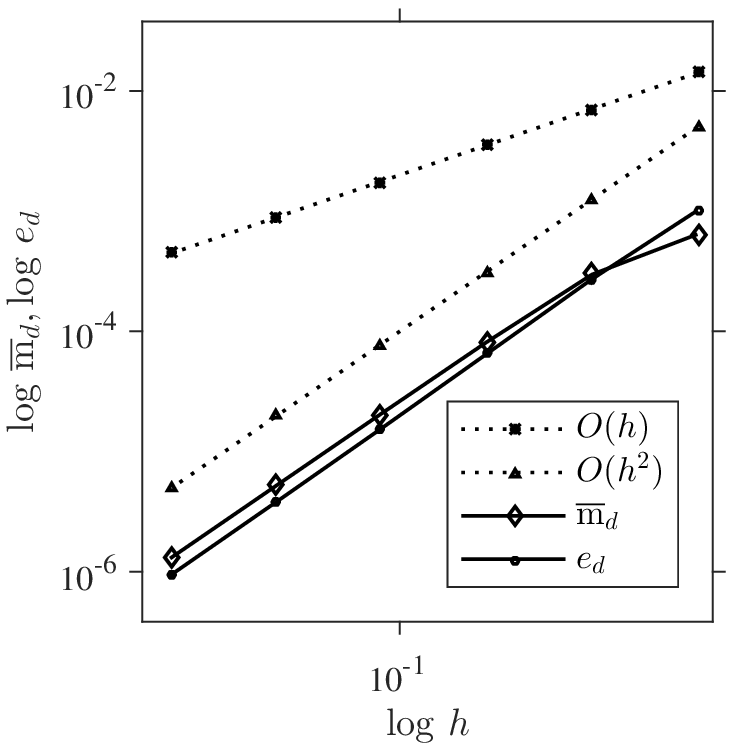}
	\label{fig:example-incr-unit-2d-uniform-convergence-b}}
	\caption{Ex. \ref{ex:incr-unit-2d-t}. The optimal convergence test of $\mdI$ and $\maj{}$.}
	\label{fig:example-incr-unit-2d-uniform-convergence}
\end{figure}
%
\begin{figure}[!ht]
	\centering
	\subfloat[$Q^{(10)}$: 128 EL, $\mathcal{T}_{9 \times 9}$ \newline
	$\incred{10} = 5.60 \cdot 10^{-4}$]{
	\includegraphics[width=5.4cm]{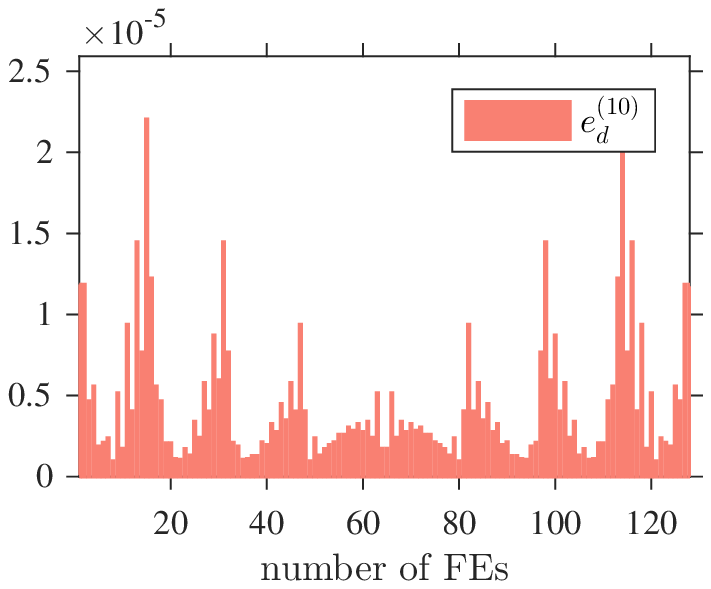}
	\label{fig:example-incr-unit-2d-e-maj-distr-128-a}}
	\subfloat[$Q^{(10)}$: 128 EL, $\mathcal{T}_{9 \times 9}$ \newline
	$\incrmdI{10} = 7.46 \cdot 10^{-4}$]{
	\includegraphics[width=5.4cm]{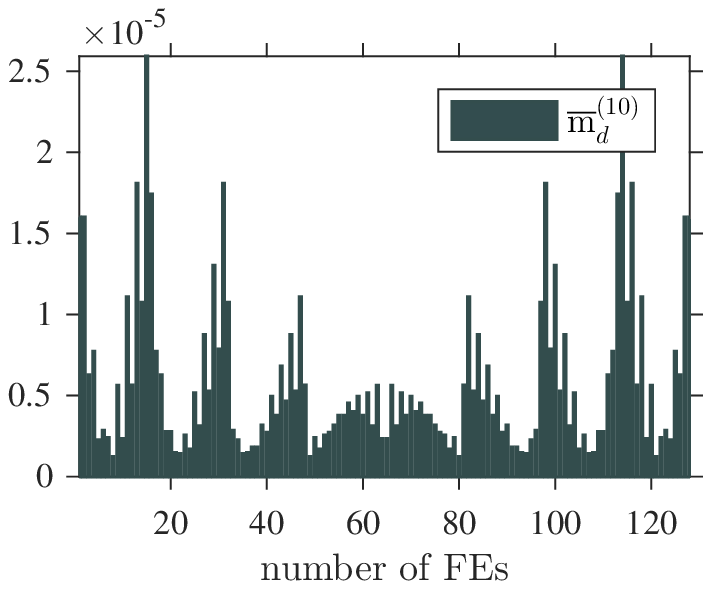}
	\label{fig:example-incr-unit-2d-e-maj-distr-128-b}}
	\subfloat[$Q^{(10)}$: 128 EL \newline
	$\incred{10}$ and $\incrmdI{10}$ overlapped]{
	\includegraphics[width=5.4cm]{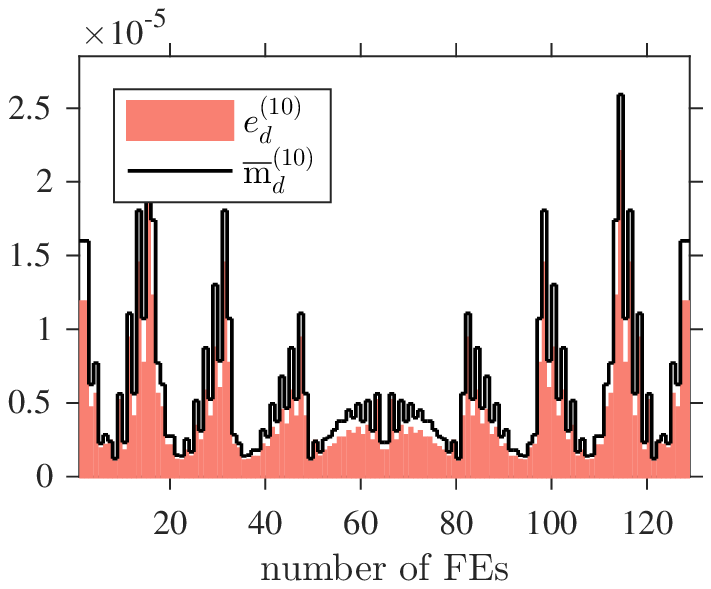}
	\label{fig:example-incr-unit-2d-e-maj-distr-128-c}}\\
	\subfloat[$Q^{(10)}$: 512 EL, $\mathcal{T}_{17 \times 17}$ \newline
	$\incred{10} = 1.36 \cdot 10^{-4}$]{
	\includegraphics[width=5.4cm]{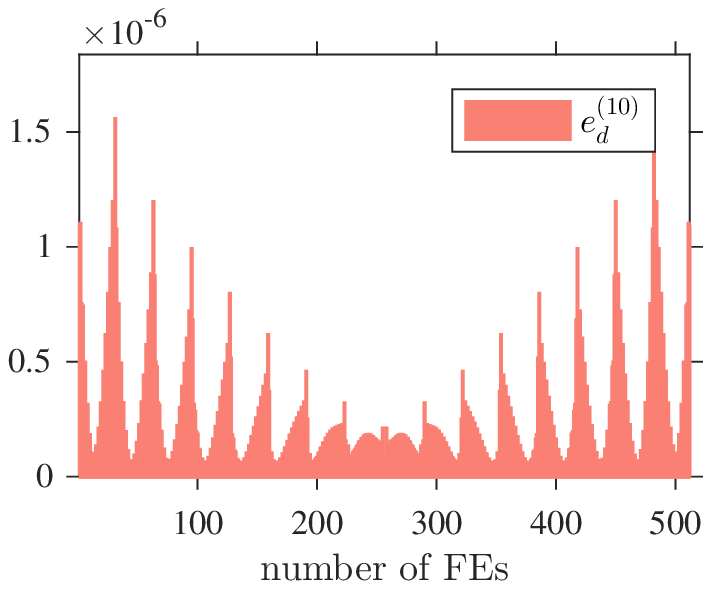}
	\label{fig:example-incr-unit-2d-e-maj-distr-512-a}}
	\subfloat[$Q^{(10)}$: 512 EL, $\mathcal{T}_{17 \times 17}$ \newline
	$\incrmdI{10} = 1.88 \cdot 10^{-4}$]{
	\includegraphics[width=5.4cm]{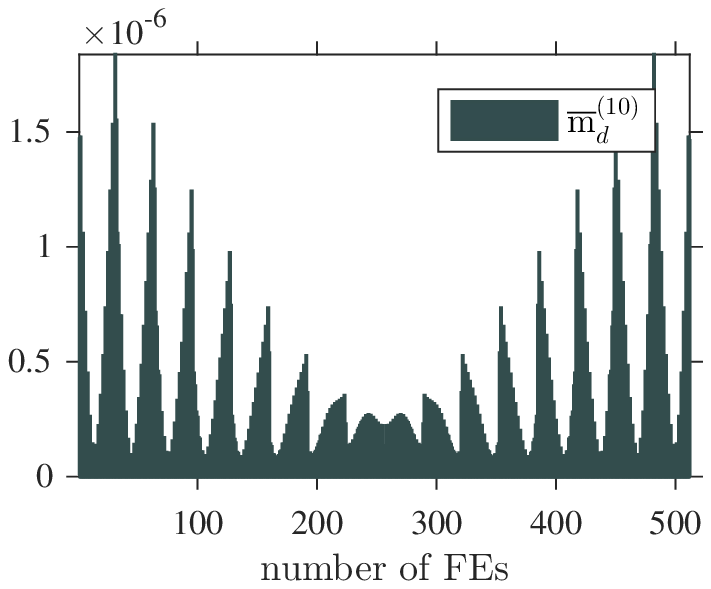}
	\label{fig:example-incr-unit-2d-e-maj-distr-512-b}}
	\subfloat[$Q^{(10)}$: 512 EL\newline
	$\incred{10}$ and $\incrmdI{10}$ overlapped]{
	\includegraphics[width=6.0cm]{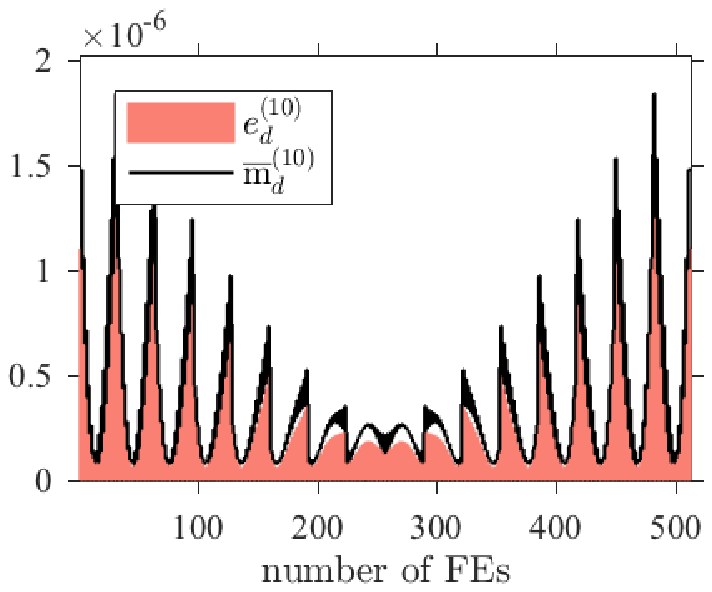}
	\label{fig:example-incr-unit-2d-e-maj-distr-512-c}}
	\caption{Ex. \ref{ex:incr-unit-2d-t}. Distribution of 
	the energy part of the error and the indicator on the slice $Q^{(10)}$.}
\end{figure}
%
Next, we solve the above formulated problem numerically with a fixed mesh on each time-step and 
confirm that $\incrmdI{k}$ represents the local distribution of $\incred{k}$ 
efficiently on each $Q^{(k)}$. We fix two different meshes $\Theta_{10 \times 18}$ 
(Figures \ref{fig:example-incr-unit-2d-e-maj-distr-128-a}--\ref{fig:example-incr-unit-2d-e-maj-distr-128-b}) 
and $\Theta_{10 \times 289}$ 
(Figures \ref{fig:example-incr-unit-2d-e-maj-distr-512-a}--\ref{fig:example-incr-unit-2d-e-maj-distr-512-b}),
where number of subintervals for discretisation in time is chosen  as $K = 10$. We 
compare the distributions of $\incred{10}$ and $\incrmdI{10}$ on the final slice $Q^{(10)}$ for both meshes. 
Here, Figures 
\ref{fig:example-incr-unit-2d-e-maj-distr-128-a} and 
\ref{fig:example-incr-unit-2d-e-maj-distr-512-a} (first column)
present the local distribution 
of $\incred{10}$ and Figures 
\ref{fig:example-incr-unit-2d-e-maj-distr-128-b} and 
\ref{fig:example-incr-unit-2d-e-maj-distr-512-b}(second column) illustrate the same element-wise 
characteristics for $\incrmdI{10}$ (vertical axis), 
whereas the cell-elements (EL) are enumerated according to the FE implementation and depicted on the horizontal axis.  
Figures \ref{fig:example-incr-unit-2d-e-maj-distr-128-c} and 
\ref{fig:example-incr-unit-2d-e-maj-distr-512-c} illustrate the same distributions 
compared (overlapped) in a single graphic. For both cases, the histograms should 
convince the reader of the quantitative efficiency of the tested error indicator. 

\begin{figure}[!ht]
	\centering
	\subfloat[$Q^{(1)}$: 128 EL \newline 
	$\incred{1} = 7.54 \cdot 10^{-5}$, $\incrmdI{1} = 1.02 \cdot 10^{-4}$]{
	\includegraphics[width=5.2cm, trim={0cm 0cm 0.1cm 0cm}, clip]{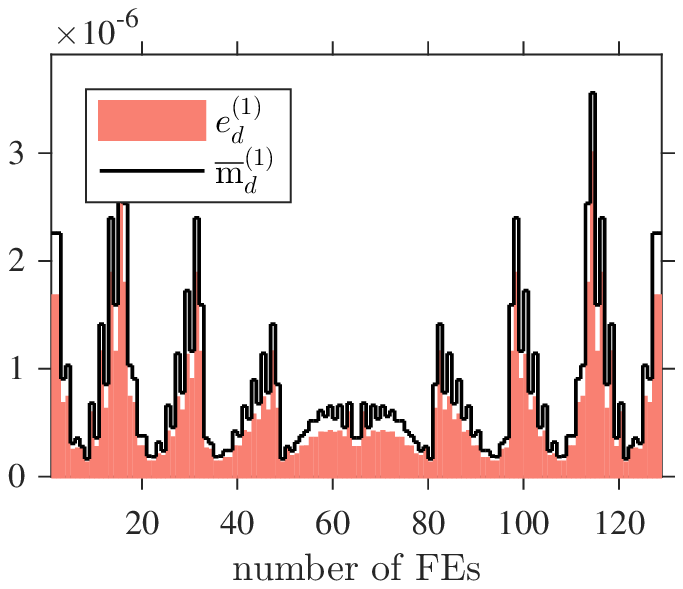}
	}
	\subfloat[$Q^{(2)}$: 280 EL \newline 
	$\incred{2} =  6.09 \cdot 10^{-5}$, $\incrmdI{2} = 7.64 \cdot 10^{-5}$]{
	\includegraphics[width=5.2cm, trim={0cm 0cm 0.1cm 0cm}, clip]{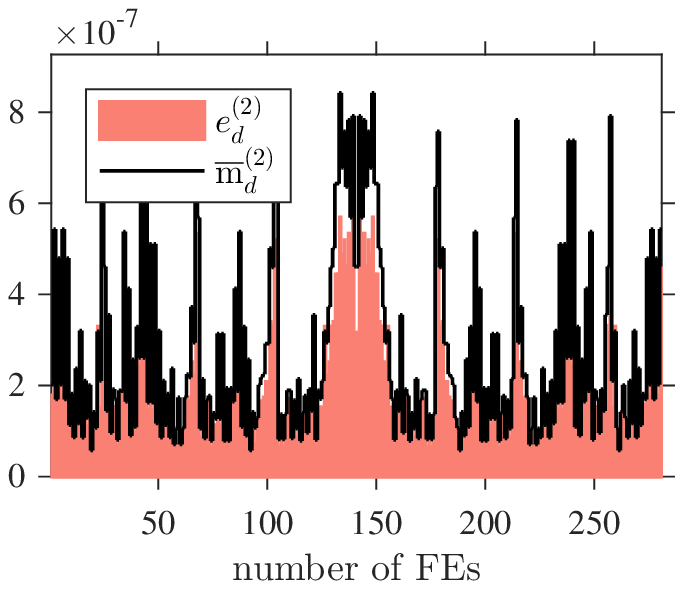}
	}
	\subfloat[$Q^{(3)}$: 688 EL \newline 
	$\incred{3} = 3.16 \cdot 10^{-5}$, $\incrmdI{3} = 3.80 \cdot 10^{-5}$]{
	\includegraphics[width=6.0cm, trim={0cm 0cm 0.1cm 0cm}, clip]{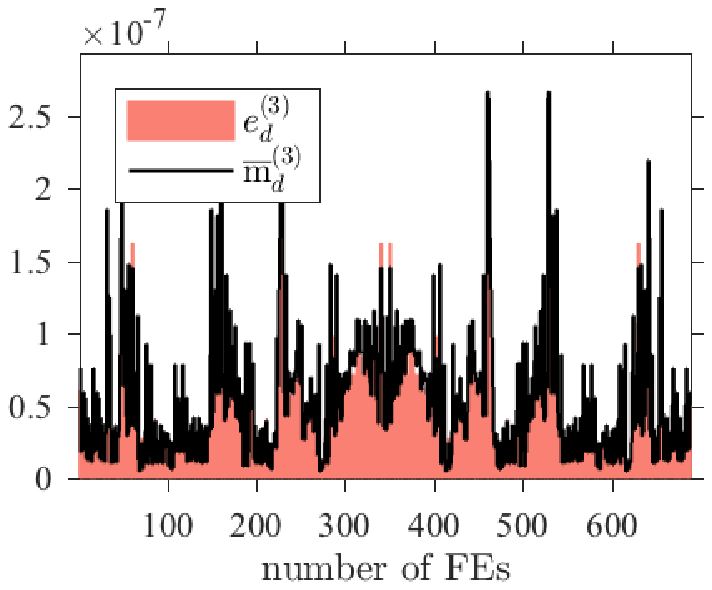}
	}
	\caption{Ex. \ref{ex:incr-unit-2d-t}. 
	True error and indicator distribution 
	on $Q^{(k)}$, $k = 1, 2, 3$, 
	(bulk marking $\Marker_{0.3}$).}
	\label{fig:example-incr-unit-2d-e-maj-distr-bulk}
\end{figure}
\begin{figure}[!ht]
	\centering
	\subfloat[$Q^{(1)}$: \; 276 EL, 165 ND 
	\newline $\incred{1} = 7.97 \cdot 10^{-5}, \incrmdI{1} = 1.02 \cdot 10^{-4}$
	]{
	\includegraphics[width=6.4cm, trim={1cm 0cm 1cm 1cm}, clip]{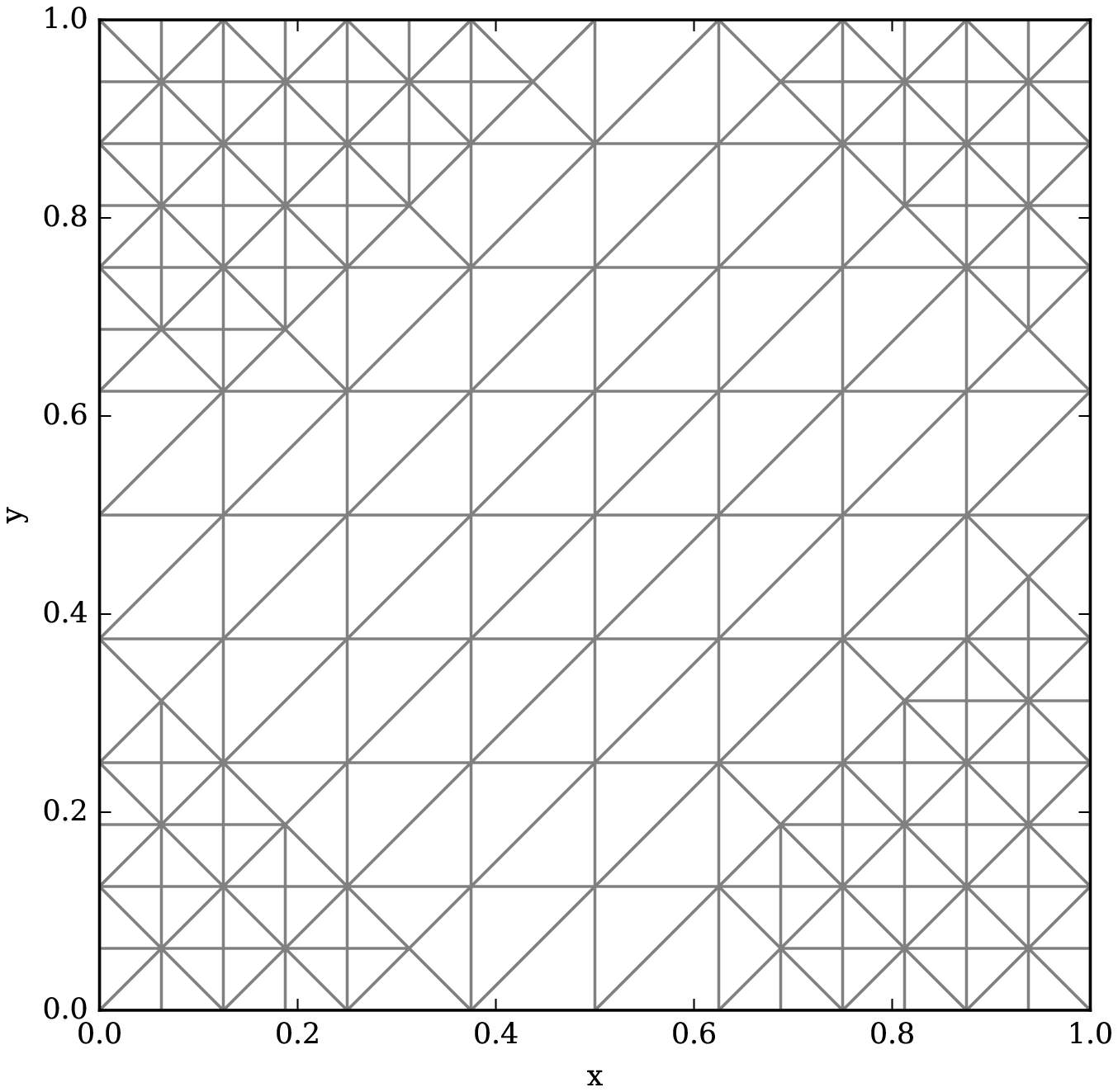}
	\label{fig:example-incr-unit-2d-error-bulk-meshes-time-1}}\!
	\subfloat[$Q^{(1)}$: \; 280 EL, 167 ND 
	\newline $\incred{1} = 6.52 \cdot 10^{-5}, \incrmdI{1} = 6.56 \cdot 10^{-5}$
	]{ 
	\includegraphics[width=6.4cm,  trim={1cm 0cm 1cm 1cm}, clip]{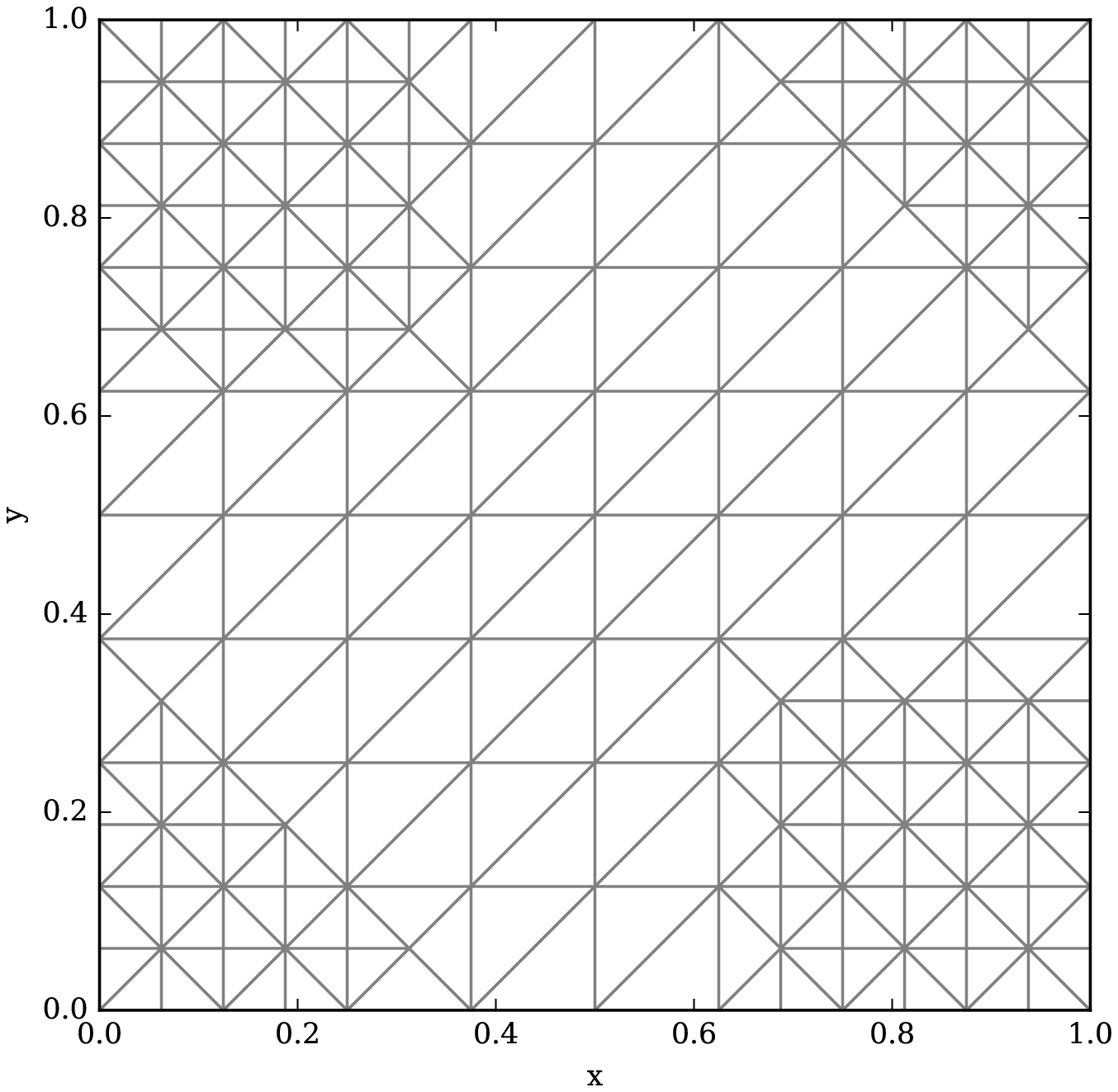}
	\label{fig:example-incr-unit-2d-majorant-bulk-meshes-time-1}}\\[-5pt]
	\subfloat[$Q^{(3)}$: \; 1544 EL, 825 ND 
	\newline $\incred{3} = 2.47 \cdot 10^{-5}, \incrmdI{3} = 2.44 \cdot 10^{-5}$]{
	\includegraphics[width=6.4cm,  trim={1cm 0cm 1cm 1cm}, clip]{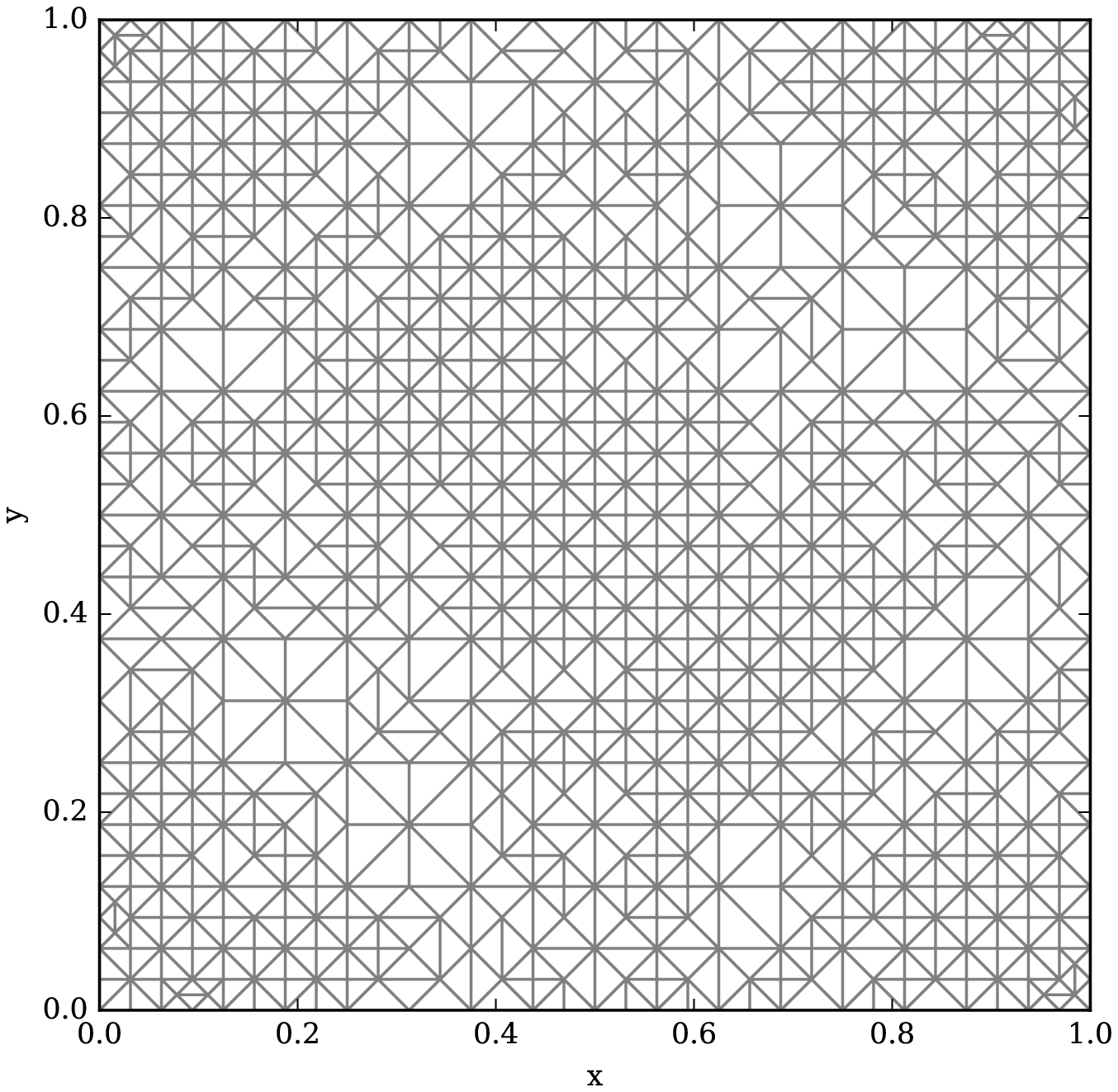}
	\label{fig:example-incr-unit-2d-error-bulk-meshes-time-3}}
	\subfloat[$Q^{(3)}$: \; 1580 EL, 843 ND 
	\newline $\incred{3} = 2.50 \cdot 10^{-5}, \incrmdI{3} = 2.47 \cdot 10^{-5}$
	]{
	\includegraphics[width=6.4cm,  trim={1cm 0cm 1cm 1cm}, clip]{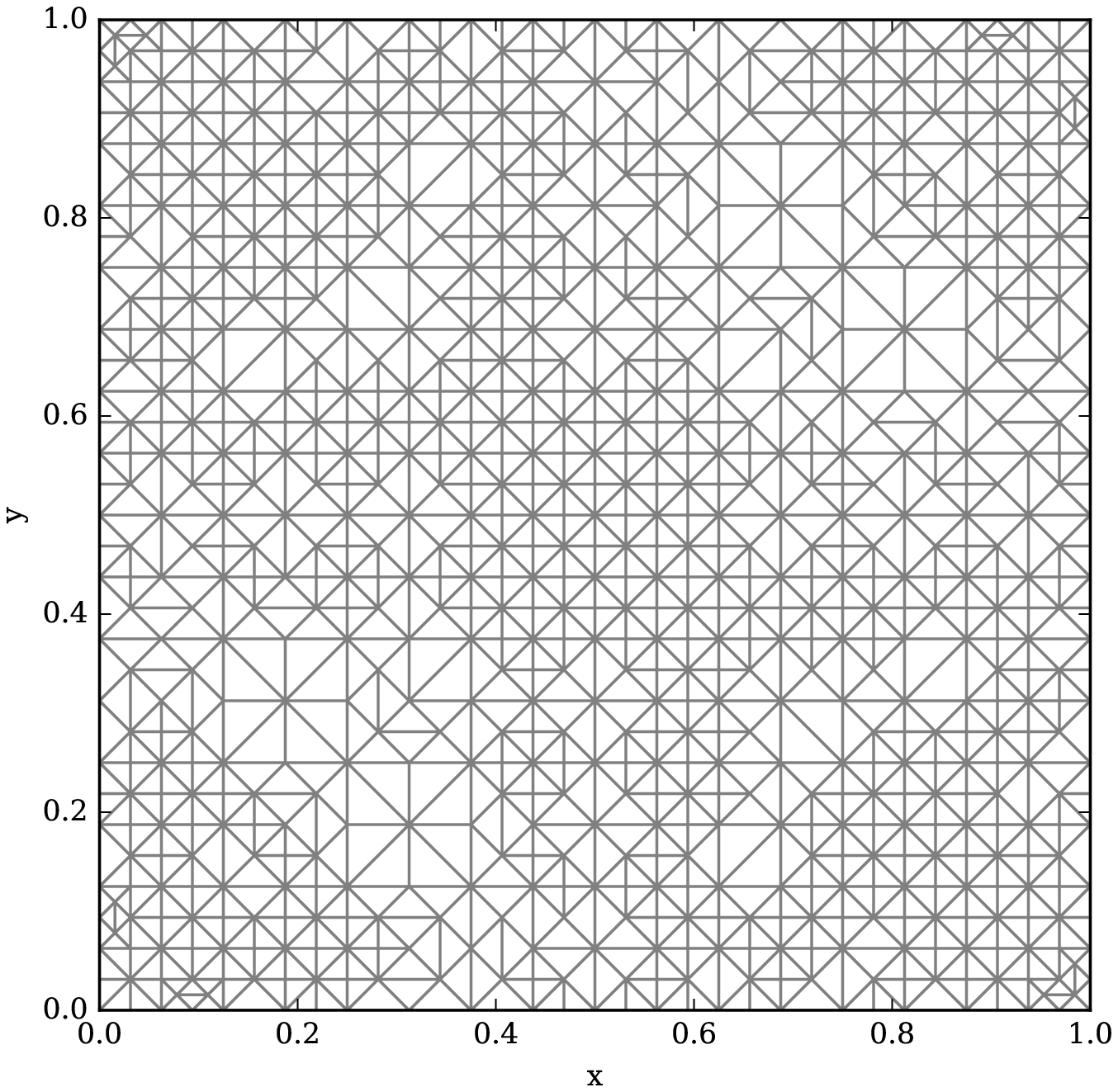}
	\label{fig:example-incr-unit-2d-majorant-bulk-meshes-time-3}}\\[-5pt]
	\subfloat[$Q^{(5)}$: \; 8819 EL, 4541 ND 
	\newline $\incred{5} = 7.48 \cdot 10^{-6}, \incrmdI{5} = 7.48 \cdot 10^{-6}$
	]{
	\includegraphics[width=6.4cm,  trim={1cm 0cm 1cm 1cm}, clip]{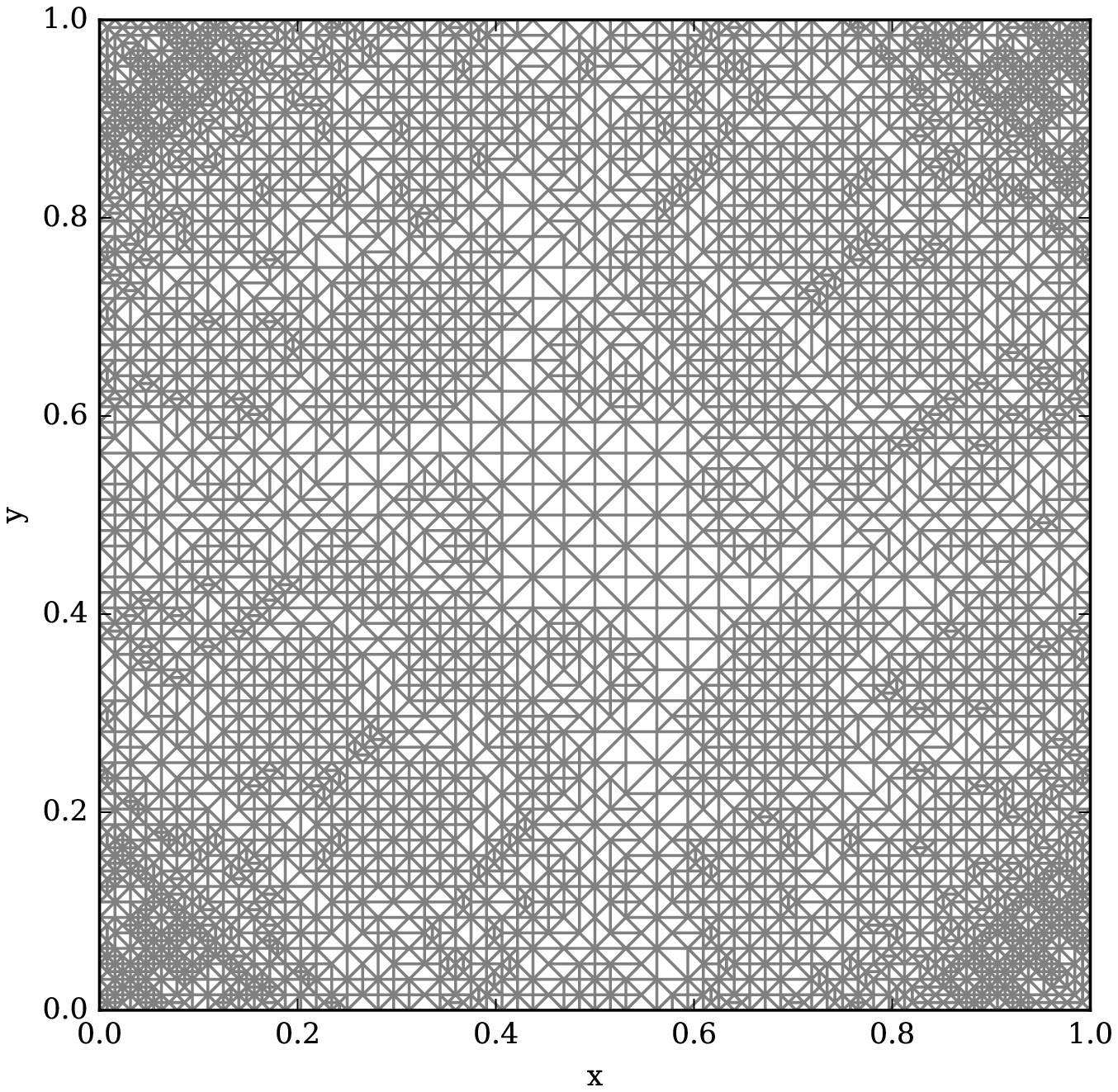}
	\label{fig:example-incr-unit-2d-error-bulk-meshes-time-5}}
	\subfloat[$Q^{(5)}$: \; 9044 EL, 4649 ND 
	\newline $\incred{5} = 7.49 \cdot 10^{-6}$, $\incrmdI{5} = 7.49 \cdot 10^{-6}$
	]{
	\includegraphics[width=6.4cm, trim={1cm 0cm 1cm 1cm}, clip]{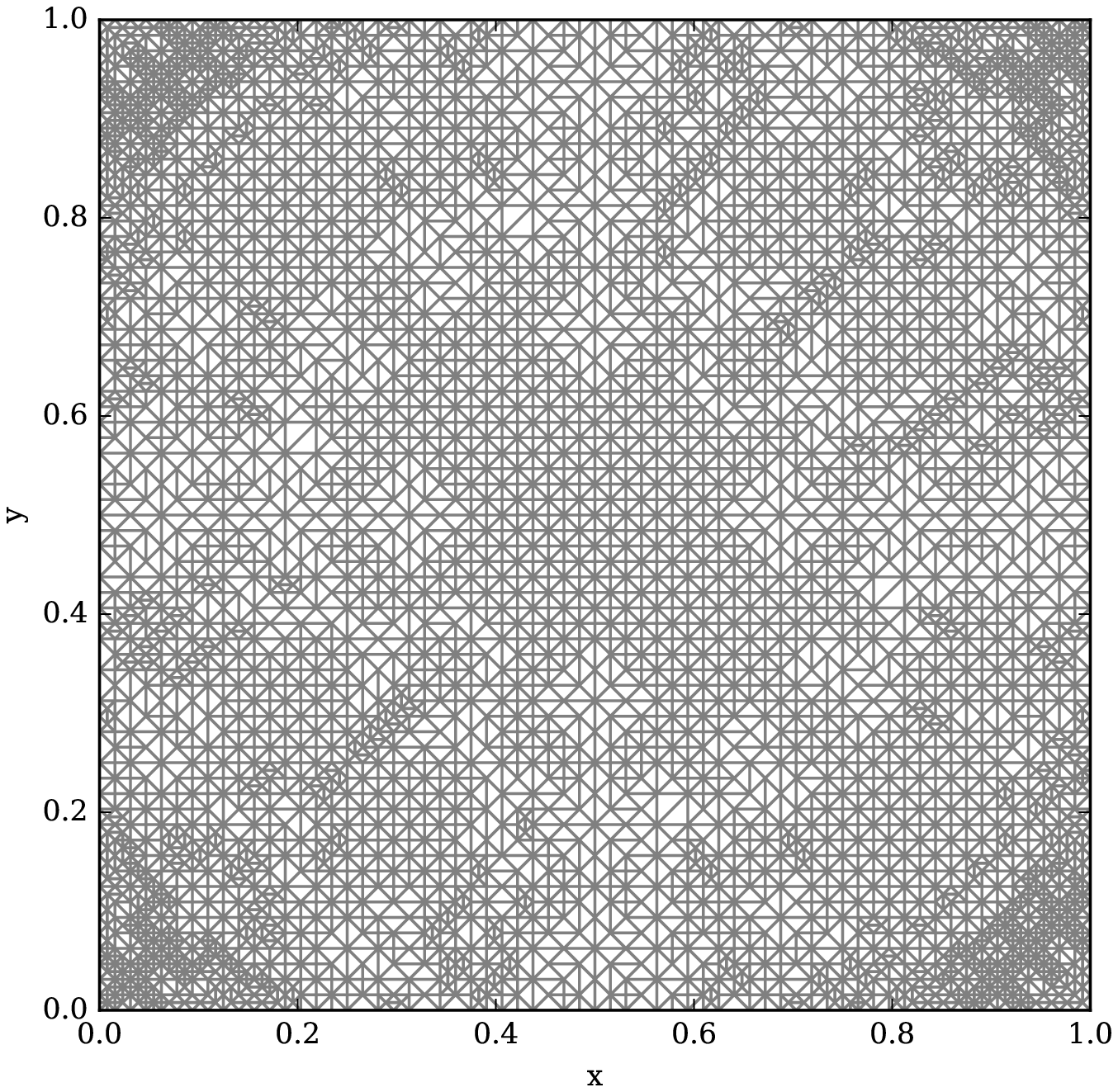}
	\label{fig:example-incr-unit-2d-majorant-bulk-meshes-time-5}}
	\caption{Ex. \ref{ex:incr-unit-2d-t}. 
	Evolution of meshes on the time-slices $Q^{(k)}$, $k = 1, 3, 5$. The refinement 
	is based on the true error (a), (c), (e) 
	and the indicator (b), (d), (f) (bulk marking $\Marker_{0.3}$).}
	\label{fig:example-incr-unit-2d-bulk-meshes}
\end{figure}

\begin{figure}[!ht]
	\centering
	\subfloat[$Q^{(1)}$: \; 284 EL,169 ND 
	\newline $\incred{1} = 7.97 \cdot 10^{-5}, \incrmdI{1} = 1.02 \cdot 10^{-4}$]{
	\includegraphics[width=6.4cm,  trim={1cm 0cm 1cm 1cm}, clip]{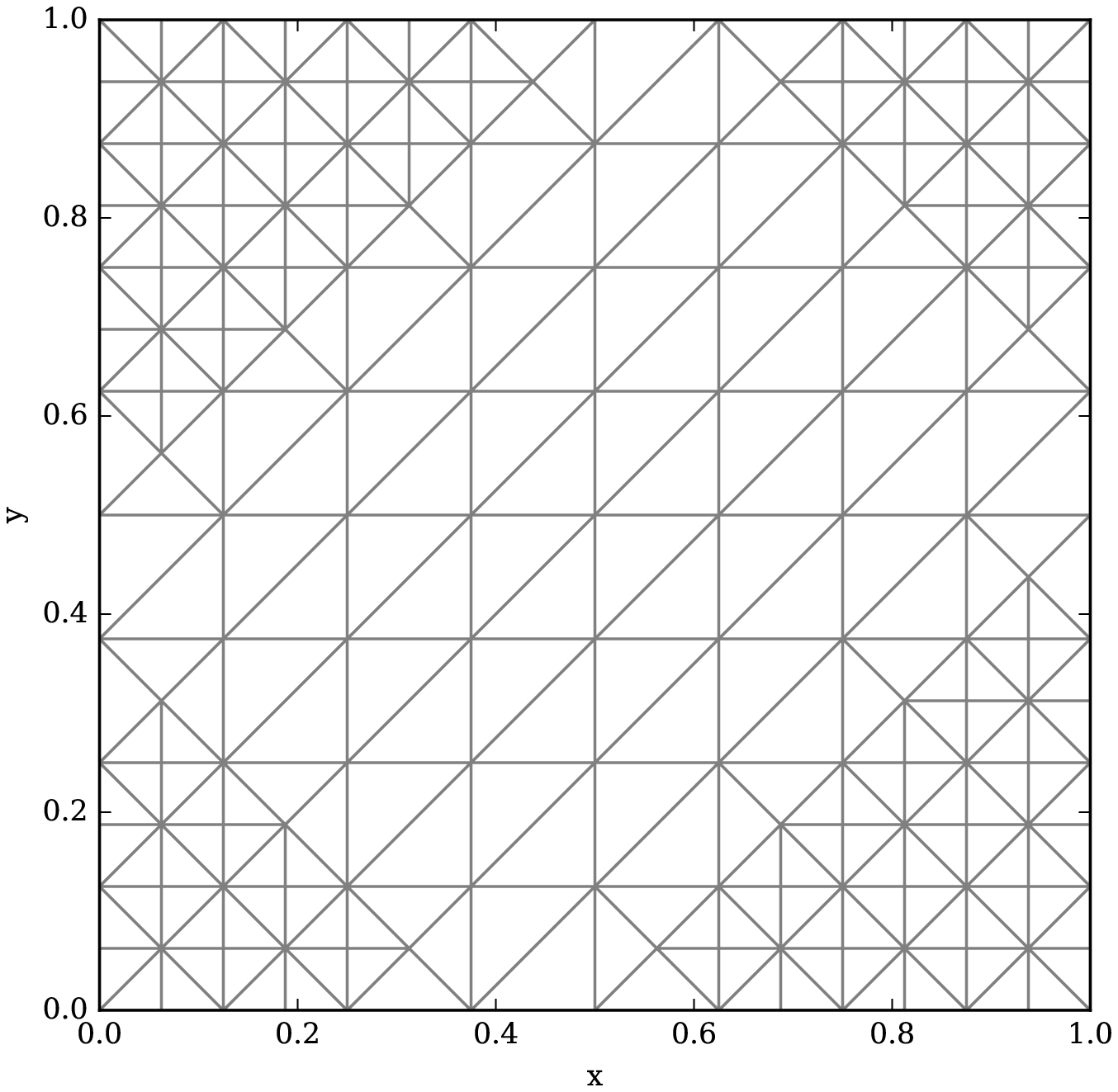}
	}
	\subfloat[$Q^{(1)}$: \; 296 EL, 175 ND 
	\newline $\incred{1} = 7.97 \cdot 10^{-5}, \incrmdI{1} = 1.02 \cdot 10^{-5}$]{
	\includegraphics[width=6.4cm,  trim={1cm 0cm 1cm 1cm}, clip]{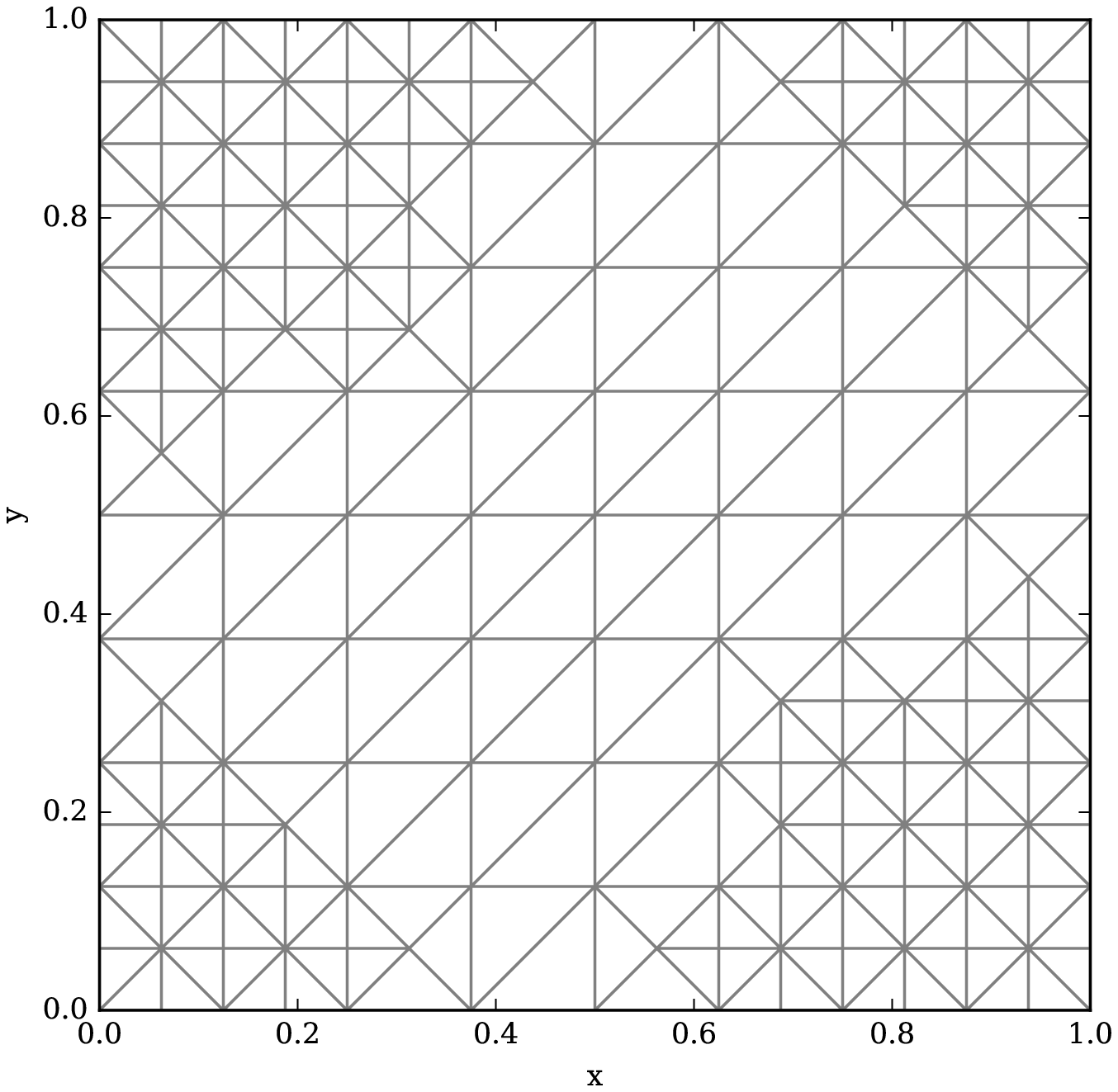}
	}\\[-5pt]
	\subfloat[$Q^{(3)}$: \; 2012 EL, 1065 ND 
	\newline $\incred{3} = 3.18 \cdot 10^{-5}, \incrmdI{3} = 3.60 \cdot 10^{-5}$]{
	\includegraphics[width=6.4cm,  trim={1cm 0cm 1cm 1cm}, clip]{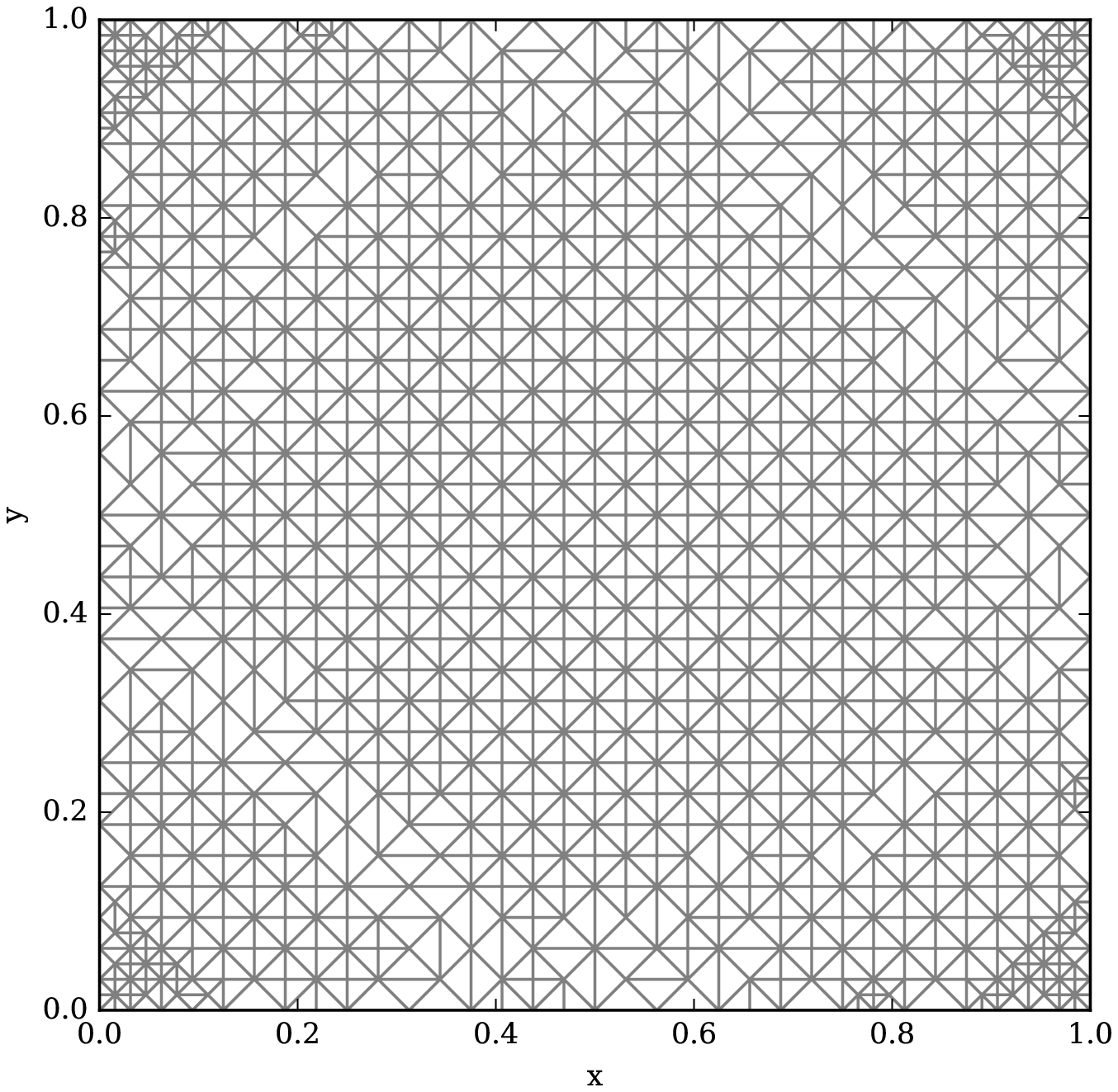}
	}
	\subfloat[$Q^{(3)}$: \; 2000 EL, 1057 ND 
	\newline $\incred{3} = 3.08 \cdot 10^{-5}, \incrmdI{3} = 3.49 \cdot 10^{-5}$]{
	\includegraphics[width=6.4cm,  trim={1cm 0cm 1cm 1cm}, clip]{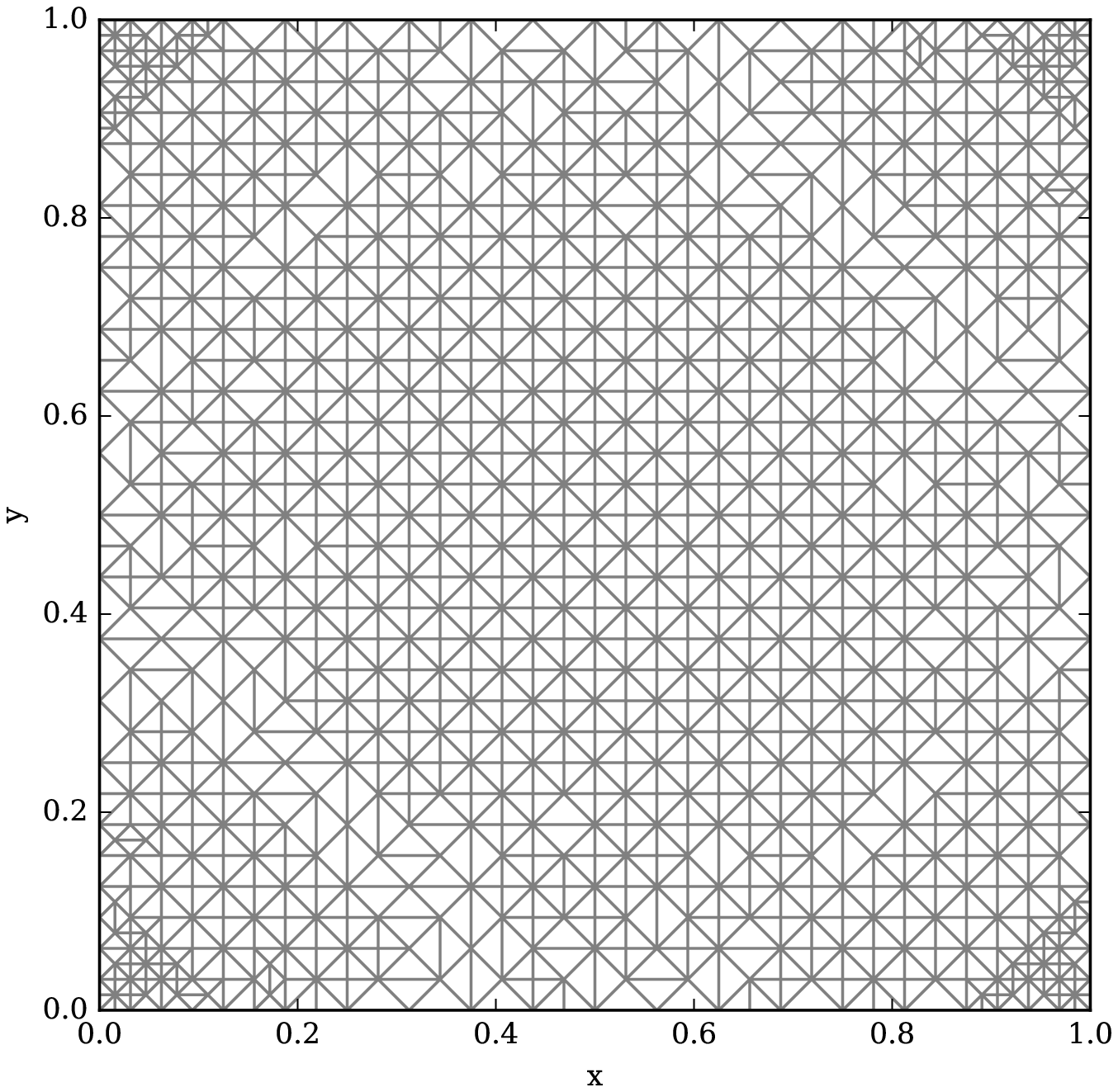}
	}\\[-5pt]
	\subfloat[$Q^{(5)}$: \; 13832 EL, 7059 ND 
	\newline $\incred{5} = 6.91 \cdot 10^{-6}, \incrmdI{5} = 8.15 \cdot 10^{-6}$]{
	\includegraphics[width=6.4cm,  trim={1cm 0cm 1cm 1cm}, clip]{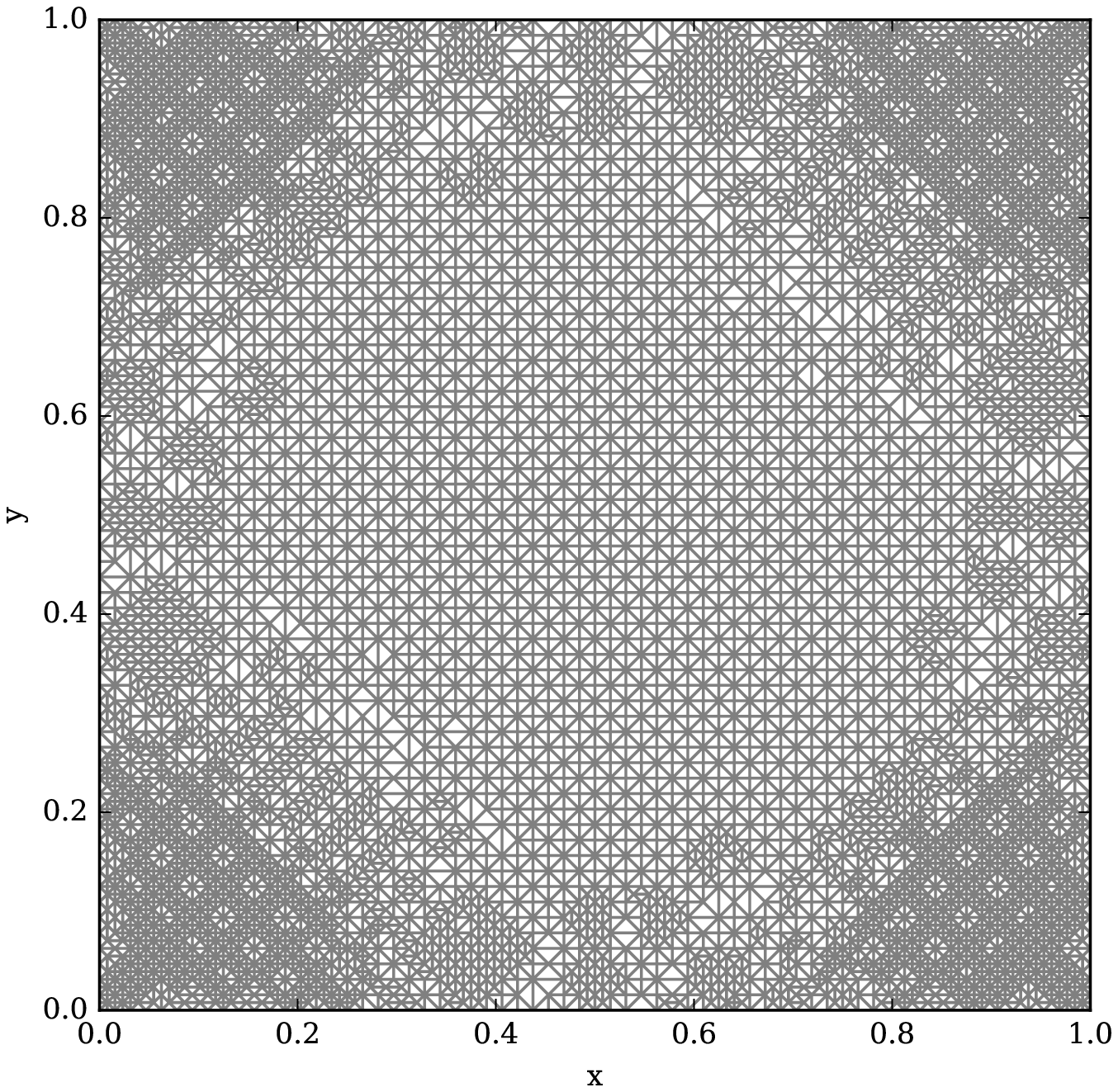}
	}
	\subfloat[$Q^{(5)}$: \; 13752 EL, 7027 ND 
	\newline $\incred{5} = 6.94 \cdot 10^{-6}, \incrmdI{5} = 5.17 \cdot 10^{-6}$]{
	\includegraphics[width=6.4cm,  trim={1cm 0cm 1cm 1cm}, clip]{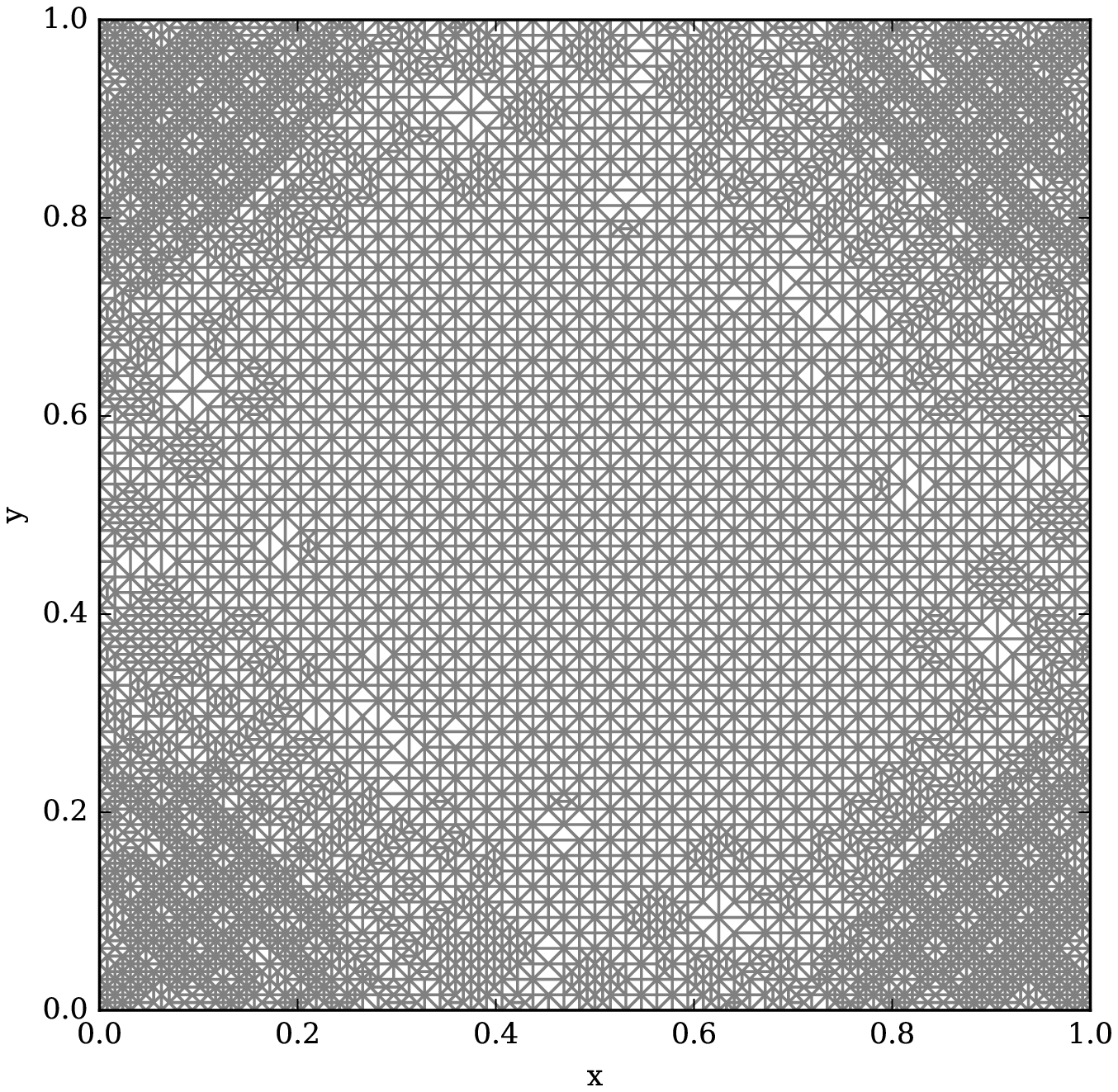}
	}\
	\caption{Ex. \ref{ex:incr-unit-2d-t}. 
	Evolution of meshes on the time-slices $Q^{(k)}$, $k = 1, 3, 5$. 
	The refinement is based on the true error (a), (c), (e) 
	and the indicator (b), (d), (f) using the marker $\Marker_{\rm AVR}$.}
	\label{fig:example-incr-unit-2d-average-meshes}
\end{figure}

As a next step, we consider an adaptive refinement strategy with the bulk marking criterion 
$\Marker_{\theta}$ introduced in \cite{Dorfler1996}. In this strategy, we form a subset of 
elements, which contain the 
highest indicated errors. The subset is formed in such a way that we keep adding elements
until the accumulated error is reached.
The selection process stops, when the error accumulated on 
the previous step exceeds the defined `bulk' level (threshold).
In this example, we set $\theta = 0.3$ (for more details on alternative marking strategies, 
we refer the reader to \cite{Malietall2014}). Let the initial mesh be $\mathcal{T}_{121}$ 
(200 EL, 121 ND).
Figure \ref{fig:example-incr-unit-2d-e-maj-distr-bulk} illustrates the distributions of 
$\incred{k}$ and $\incrmdI{k}$ on different cylinder slices $Q^{(k)}$, $k = 1, 2, 3$, 
and demonstrates the quantitative efficiency of the indicator provided by the majorant. 
Under every sub-plot of Figure \ref{fig:example-incr-unit-2d-e-maj-distr-bulk}, we also 
provide the information on total values of $\incred{k}$ and $\incrmdI{k}$.
%

We analyse the meshes obtained during the refinement based on either 
$\incred{k}$ or $\incrmdI{k}$ (Figure \ref{fig:example-incr-unit-2d-bulk-meshes}). 
In Figures \ref{fig:example-incr-unit-2d-error-bulk-meshes-time-1}, 
\ref{fig:example-incr-unit-2d-error-bulk-meshes-time-3}, and 
\ref{fig:example-incr-unit-2d-error-bulk-meshes-time-5} (left column), we present 
the meshes obtained after the refinement process based on the local true error distribution, 
and Figures \ref{fig:example-incr-unit-2d-majorant-bulk-meshes-time-1}, 
\ref{fig:example-incr-unit-2d-majorant-bulk-meshes-time-3}, and 
\ref{fig:example-incr-unit-2d-majorant-bulk-meshes-time-5} (right column) expose 
the meshes constructed when the refinement is based on the local indicator $\incrmdI{k}$. It is easy 
to observe that the topology of the meshes on the RHS of Figure 
\ref{fig:example-incr-unit-2d-bulk-meshes} resembles the topology of the meshes on 
the LHS. 
For 
this case, the efficiency of the total majorant is $\Ieff =  1.23$.
%

\newpage
The bulk marking strategy can be compared to marking determined by the level 
of the average error $\Marker_{\rm AVR}$ \cite[Algorithm 2.1]{Malietall2014}.  
The latter one marks those elements, on which the error exceeds the average level. 
Figure \ref{fig:example-incr-unit-2d-average-meshes} demonstrates the sequence of 
the meshes obtained as a result of the refinement based on $\incred{k}$ 
(LHS) and $\incrmdI{k}$ (RHS) on the respective slices $Q^{(k)}$, $k = 1, 3, 5$.
In this case, we obtain the efficiency $\Ieff = 1.4$. 
The majorant is not as accurate as expected, due to the fact that, unlike in the error 
majorants for the elliptic BVP, there is always a gap between $\maj{}$ and $\error$ 
caused by the term $v_t$ in $\mfI$ (see also \eqref{eq:gap}).

\begin{table}[!ht]
\centering
\footnotesize
\begin{tabular}{c|cccc|cccc}
\multicolumn{1}{c|}{$ $}
& \multicolumn{4}{c|}{ Implicit scheme }
& \multicolumn{4}{c}{ Explicit scheme } \\
\midrule
$k$ & 
DOFs$({v})$ & $\error$ & $\maj{}$ & $\Ieff$ & 
DOFs$({v})$ &  $\error$ & $\maj{}$ & $\Ieff$\\
\midrule
1 & 14641 & $2.29 \cdot 10^{-9}$ & $3.78 \cdot 10^{-9}$ & 1.29 & 14641 & $1.26 \cdot 10^{-6}$ & $7.89 \cdot 10^{-5}$ & 7.93 \\
2 & 23627 & $4.01 \cdot 10^{-9}$ & $6.84 \cdot 10^{-9}$ & 1.31 & 27175 & $2.06 \cdot 10^{-3}$ & $4.56 \cdot 10^{-3}$ & 1.49 \\
3 & 39795 & $5.05 \cdot 10^{-9}$ & $9.14 \cdot 10^{-9}$ & 1.35 & 45489 & $9.19 \cdot 10^{3}\;\;$  & $1.46 \cdot 10^{4}\;\;$ & 1.26 \\
4 & 67719 & $5.66 \cdot 10^{-9}$ & $1.06 \cdot 10^{-8}$ & 1.37 & 82344 & $1.15 \cdot 10^{12}$ & $1.63 \cdot 10^{12}$ & 1.19 \\
\end{tabular}
\caption{Ex. \ref{ex:incr-unit-2d-t}. 
Total error, majorant, and efficiency index for $v$ 
generated by implicit and explicit schemes.}
\label{tab:explicit-implicit-schemes}
\end{table}

\newpage
It is important to note that the majorant can be used as a tool to predict the 
`blow-ups' in time-dependent explicit schemes, which are less time-consuming in 
comparison to the implicit ones, but are unstable. Furthermore, if for one-dimensional (in space)
schemes, the stability condition is written explicitly, for two- and three-dimensional 
problems as well as for non-linear problems there are no such criteria. As an example, 
we consider the mesh $\Theta_{1280 \times 14641}$ (28800 EL, 14641 ND) 
and illustrate the majorant response on the instability of the explicit scheme (see Table 
\ref{tab:explicit-implicit-schemes}). Here, the column DOFs$({v})$ reflects 
the number of DOFs for $v$. The LHS of the table contains total values of the error 
and the majorant, obtained by using a stable implicit scheme, whereas the RHS 
illustrates a drastic increase of the majorant even when the `blow-up' is not 
yet obvious from the error values. 
\end{example}

\begin{example}
\label{ex:incr-unit-3d-t}
\rm
{
In order to demonstrate that the same behaviour of the error estimates can be 
observed for the problems in higher dimension, we consider 
the unit cube} $\Omega = (0, 1)^3 \subset \Rthree$ with $T = 1$, $A = I$, 
$\vectorb = {\boldsymbol 0}$, $c = 0$,  
initial condition $u_0 = x\,(1 - x)\,y\,(1 - y)\,z\,(1 - z)$, homogeneous Dirichlet BC, 
and $$u = x\,(1 - x)\,y\,(1 - y)\,z\,(1 - z)(t^2 + t + 1).$$ Analogously, we take 
$v \in \Pone$. However, in the current 
example, we compare performance of the majorant reconstructed with two different 
approximations of the flux, i.e., approximated by Raviart-Thomas FEs of the lowest order 
$\flux \in \RTzero$ and a linear one $\flux \in \RTone$. Figure 
\ref{fig:example-incr-unit-3d-uniform-convergence-a} demonstrates the uniform 
convergence of $\error$ and $\maj{} (\flux)$ with $\flux \in \RTzero$, and Figure 
\ref{fig:example-incr-unit-3d-uniform-convergence-b} depicts the same 
characteristics for $\flux \in \RTone$. They both confirm the optimal 
convergence order of the majorant constructed with $\flux \in \RTzero$ and 
$\flux \in \RTone$. 

\begin{figure}[!ht]
	\centering
	\subfloat[$\error$, $\maj{} (\flux)$, $\flux \in {\rm RT}_0$]
	{\includegraphics[width=6cm]{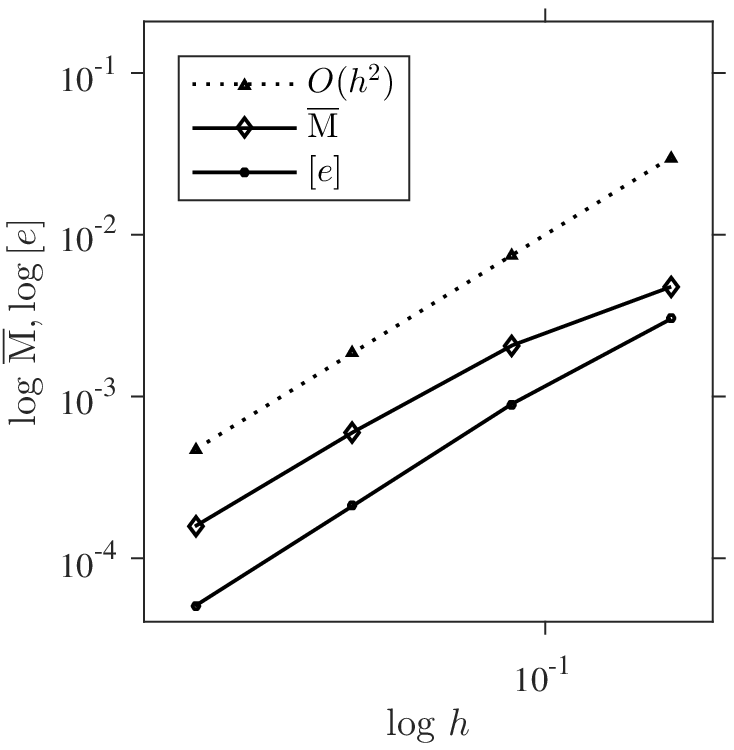}
	\label{fig:example-incr-unit-3d-uniform-convergence-a}}\quad
	\subfloat[$\error$, $\maj{} (\flux)$, $\flux \in {\rm RT}_1$]
	{\includegraphics[width=5.9cm]{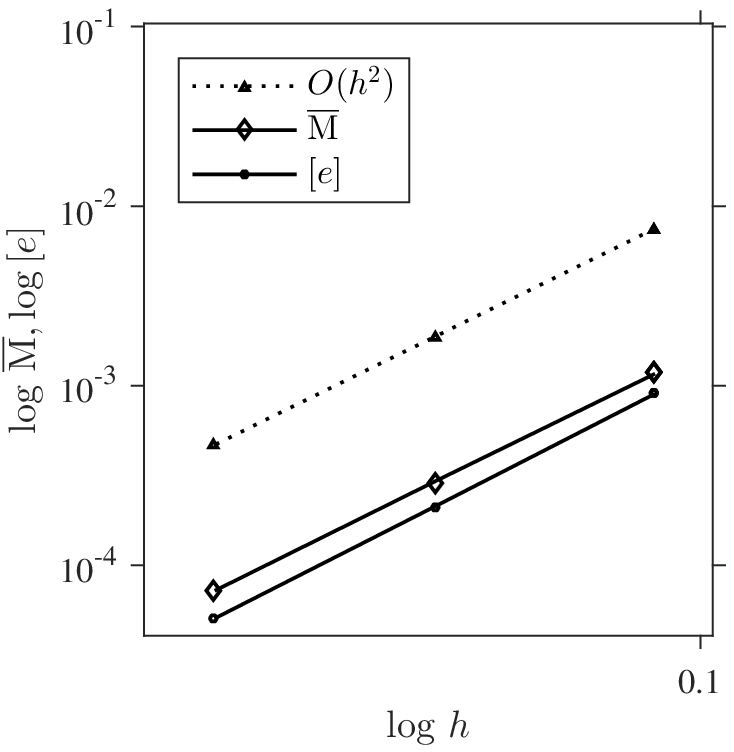}
	\label{fig:example-incr-unit-3d-uniform-convergence-b}}
	\caption{Ex. \ref{ex:incr-unit-3d-t}. 
	The optimal convergence of $[e]$ and $\maj{}$, 
	(a) $\flux \in {\rm RT}_0$ and (b) $\flux \in {\rm RT}_1$.}
	\label{fig:example-incr-unit-3d-uniform-convergence}
\end{figure}

\begin{figure}[!ht]
	\centering
	\qquad 
	\subfloat[$Q^{(10)}$: 24576 EL (sorted) \newline 
	$\incred{10} = 9.13 \cdot 10^{-5}$, $\incrmdI{10} = 1.16 \cdot 10^{-4}$]
	{\includegraphics[width=6cm]{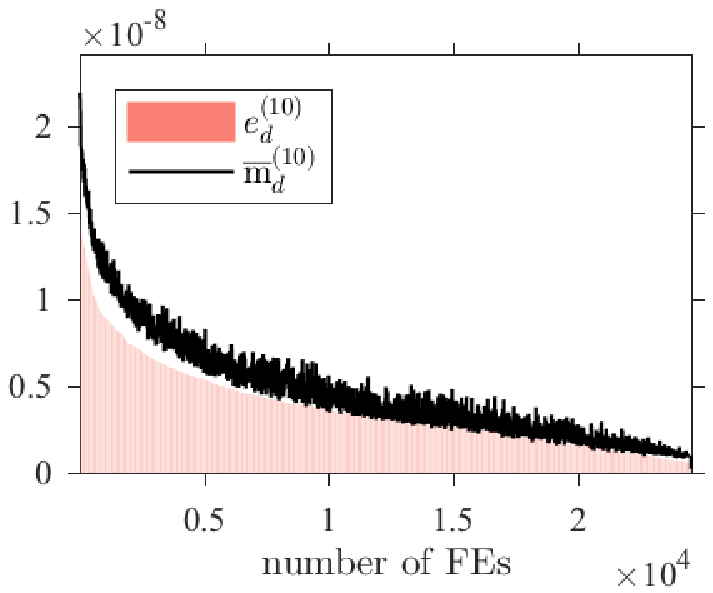}
	\label{fig:example-incr-unit-3d-e-maj-distr-a}}\quad 
	\subfloat[$Q^{(10)}$: 24576 EL (sorted) \newline 
	$\incred{10} = 9.13 \cdot 10^{-5}$, $\incrmdI{10} = 9.15 \cdot 10^{-5}$]
	{\includegraphics[width=6cm]{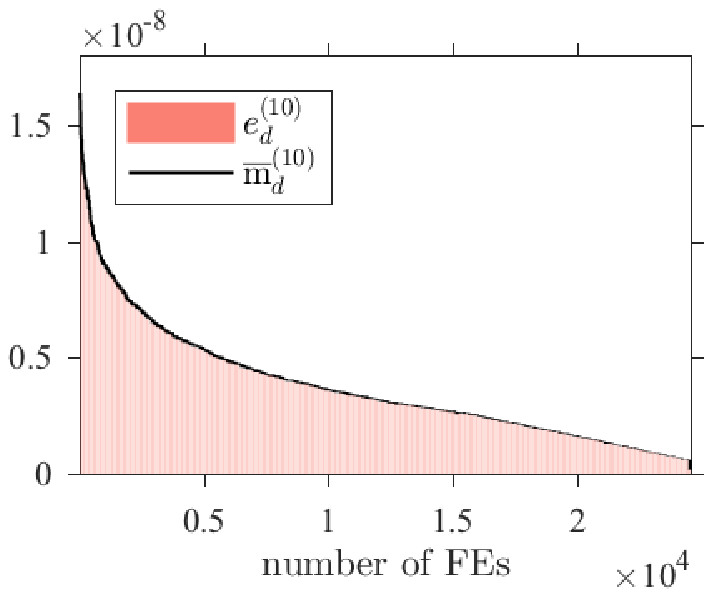}
	\label{fig:example-incr-unit-3d-e-maj-distr-b}}
	\caption{Ex. \ref{ex:incr-unit-3d-t}. Energy parts of the true error and the indicator 	
	distributions based on (a) $\flux \in {\rm RT}_0$ and (b) $\flux \in {\rm RT}_1$.}
	\label{fig:example-incr-unit-3d-e-maj-distr}
\end{figure}

\begin{figure}[!ht]
	\centering
	\vskip-5pt
	\subfloat[$Q^{(0)}$: 48 EL
	\newline $\incred{0} = 2.10 \cdot 10^{-2}$, $\incrmaj{0} = 8.43 \cdot 10^{-3}$]{
	\includegraphics[width=5.4cm]{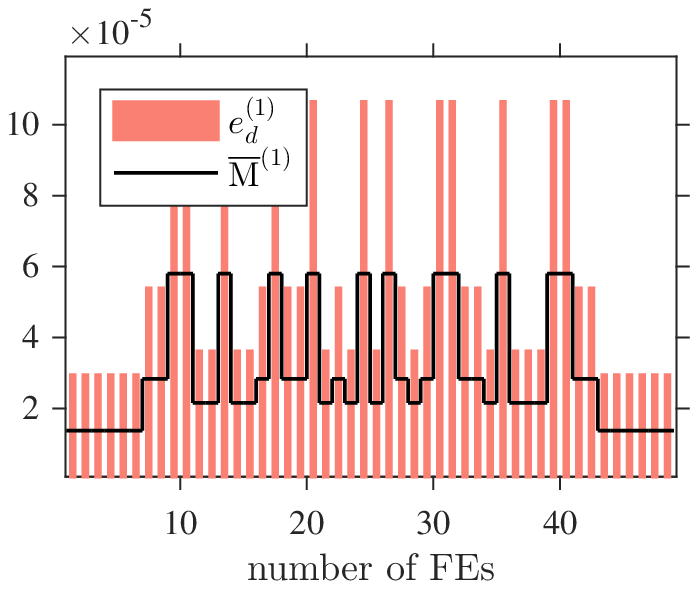}
	\label{fig:example-incr-unit-3d-e-maj-distr-bulk-marking-theta-20-a}}
	\subfloat[$Q^{(1)}$: 274 EL 
	\newline  $\incred{1} = 1.14 \cdot 10^{-1}$, $\incrmaj{1} = 7.38 \cdot 10^{-2}$]{
	\includegraphics[width=5.4cm]{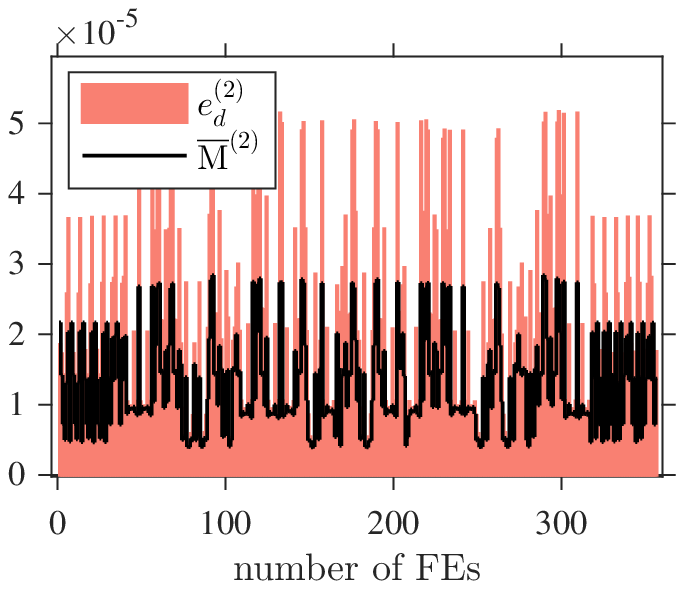}
	\label{fig:example-incr-unit-3d-e-maj-distr-bulk-marking-theta-20-b}}
	%
	\caption{Ex. \ref{ex:incr-unit-3d-t}. Error and indicator distributions for the time-slices $Q^{(k)}$, 
	$k = 0, 1$ ($\flux \in {RT}_1$).}
	\label{fig:example-incr-unit-3d-e-maj-distr-bulk-theta-20}
\end{figure}
%

%
\begin{figure}[!ht]
	\centering
	\!\!\!\!\!\!\!\!\!
	\subfloat[$Q^{(0)}$: 48 EL]{
	\includegraphics[width=5.6cm, trim={5cm 1cm 4cm 2cm}, clip]{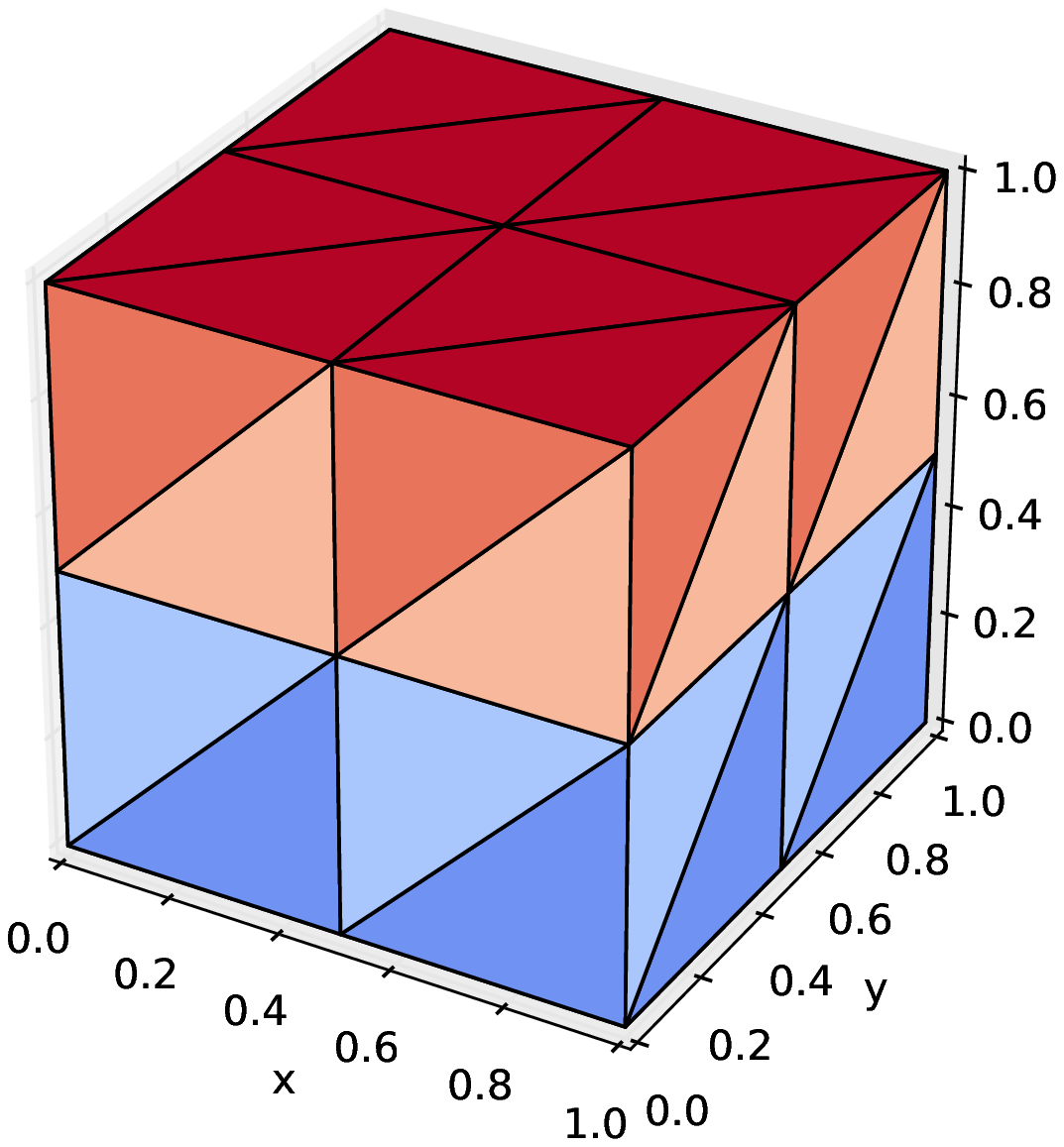}
	\label{fig:example-incr-unit-3d-e-maj-distr-bulk-marking-theta-60-a}}
	\subfloat[$Q^{(1)}$: 356 EL]{
	\includegraphics[width=5.6cm, trim={5cm 1cm 4cm 2cm}, clip]{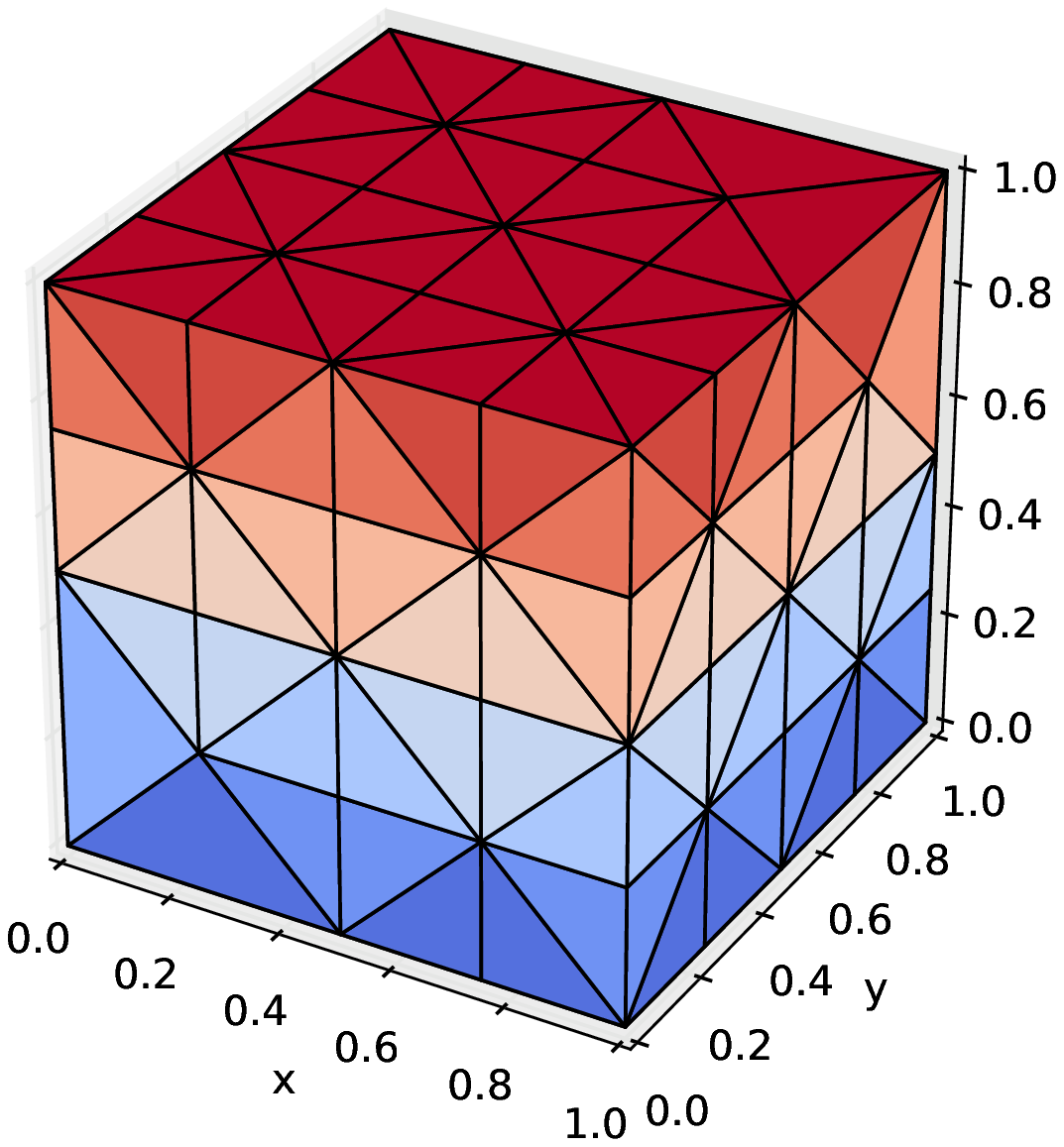}
	\label{fig:example-incr-unit-3d-e-maj-distr-bulk-marking-theta-60-b}} 
	%
	%
	\subfloat[$Q^{(3)}$: 5678 EL]{
	\includegraphics[width=5.6cm,  trim={5cm 1cm 4cm 2cm}, clip]{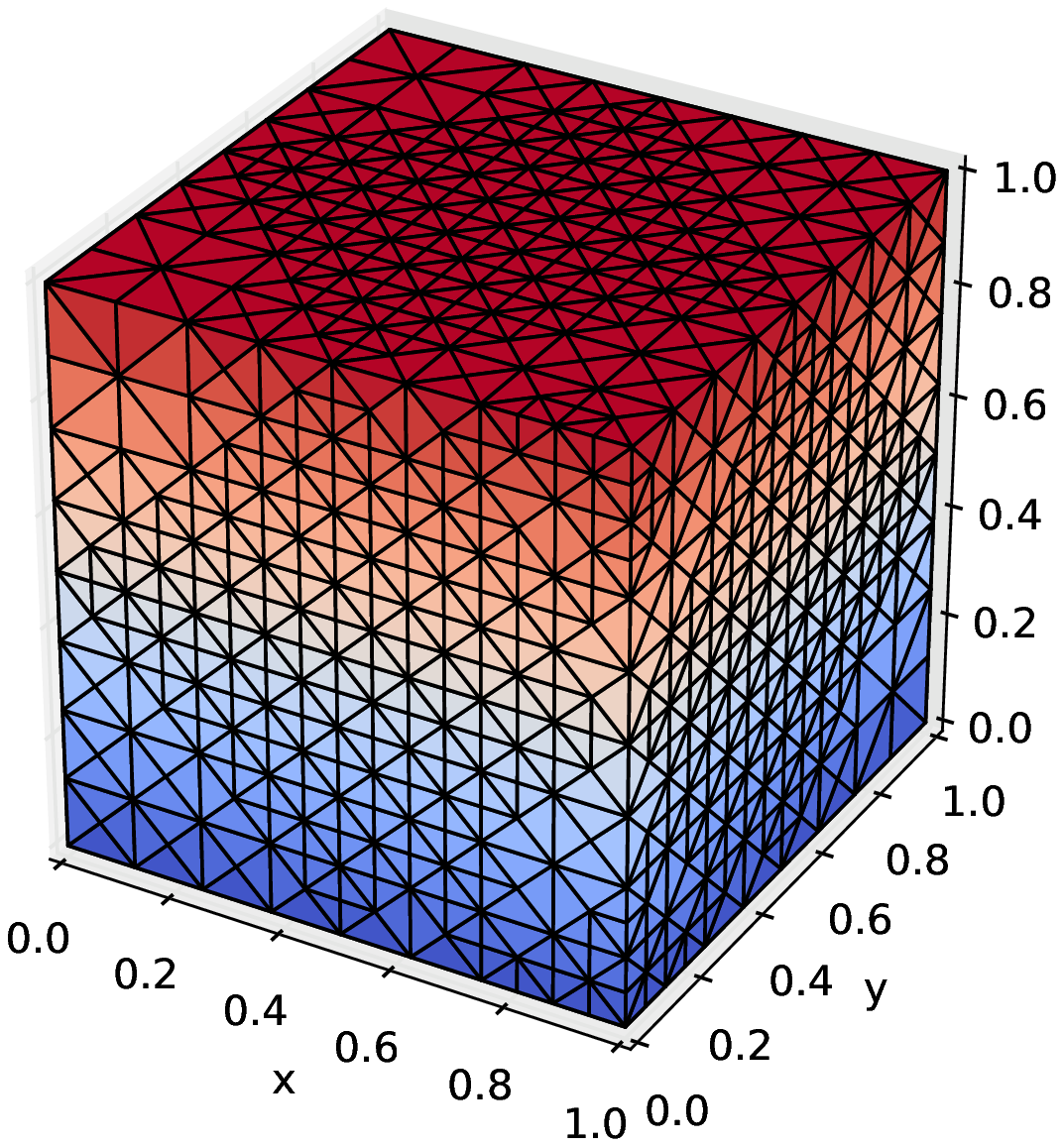}
	\label{fig:example-incr-unit-3d-e-maj-distr-bulk-marking-theta-20-d}}
	\caption{Ex. \ref{ex:incr-unit-3d-t}. Evolution of meshes on the time-slices $Q^{(k)}$, $k = 0, 1, 2, 3$, 
	the refinement is based on the majorant using the marker $\Marker_{\rm 0.6}$.}
	\label{fig:example-incr-unit-3d-mesh-bulk-marking-theta-60}
\end{figure}

Next, we compare the indicators that are reconstructed based on the fluxes, approximated 
with elements of different degrees, i.e., $\RTzero$ and $\RTone$. 
Figure \ref{fig:example-incr-unit-3d-e-maj-distr-a} presents the 
error distributions and the error indicator on the same plot, but (unlike previous 
histograms) the arrays of cells are sorted such that the values 
of the local true errors $\incred{10}$ decrease. The array with $\incrmdI{10}$ is 
depicted in the order defined by the indices obtained as a result of sorting. The 
evaluation of the error and the indicator distributions by such histograms was introduced 
in \cite[Section 3.4]{Malietall2014}. It is easy to see that $\incrmdI{10}$ 
(see Figure \ref{fig:example-incr-unit-3d-e-maj-distr-a}) is less efficient than the 
one reconstructed from $\flux \in \RTone$ in Figure 
\ref{fig:example-incr-unit-3d-e-maj-distr-b}. Latter illustrations reaffirm that 
using the fluxes of higher regularity is more advantageous for the majorant 
reconstruction, i.e., it guarantees a rather sharp and efficient prediction of the local error 
distribution. 

Finally, we consider the refinement strategy with bulk marking $\Marker_{0.6}$.
We take the initial mesh $\mathcal{T}_{3 \times 3}$, $K = 10$, and illustrate 
the obtained distributions after the refinement on the time-slices 
$k = 0, 1$ (Figure \ref{fig:example-incr-unit-3d-mesh-bulk-marking-theta-60}).
The number of obtained elements and total values of 
$\incred{k}$ and $\incrmdI{k}$ are pycharmshown below the plots. The evolution of 
meshes (obtained during the adaptation procedure from one slice to another, $k = 0, 1, 3$) is 
shown in Figure \ref{fig:example-incr-unit-3d-mesh-bulk-marking-theta-60}.
\end{example}

\begin{example}
\label{ex:incr-l-shape-2d-t}
\rm
When the question of the efficient error indication is concerned, it is important 
to demonstrate that the studied majorant efficiently catches the error-jumps related 
to various singularities.
One of the classical benchmark examples for such kind of testing 
is the problem defined on the $L$-shaped domain 
$\Omega := (\minus 1, 1) \times (\minus 1, 1) 
\backslash [0, 1) \times [0, \minus 1)$ with $T = 1$, $A = I$, 
$\vectorb = {\boldsymbol 0}$, $c = 0$,
the Dirichlet BC with the load $u_D = r^{1/3} \, \sin \theta$ with  $r = (x^2 + y^2)$ and 
$\theta = \tfrac23 \, {\rm atan2} (y, x)$ on $\Sigma$, the input source function 
$f = r^{1/3} \, \sin \theta \, \left(2 \, t + 1 \right)$, and the initial condition 
$u_0 = r^{1/3} \, \sin \theta$. 
The corresponding exact solution $u = r^{1/3} \,\sin\theta\,\left(t^2 + t + 1 \right)$ 
has a singularity at the point $(r, \theta) = (0, 0)$.

%
\begin{figure}[!ht]
	\centering
	\subfloat[$\ed$ and $\mdI$]{
	\includegraphics[width=6cm]{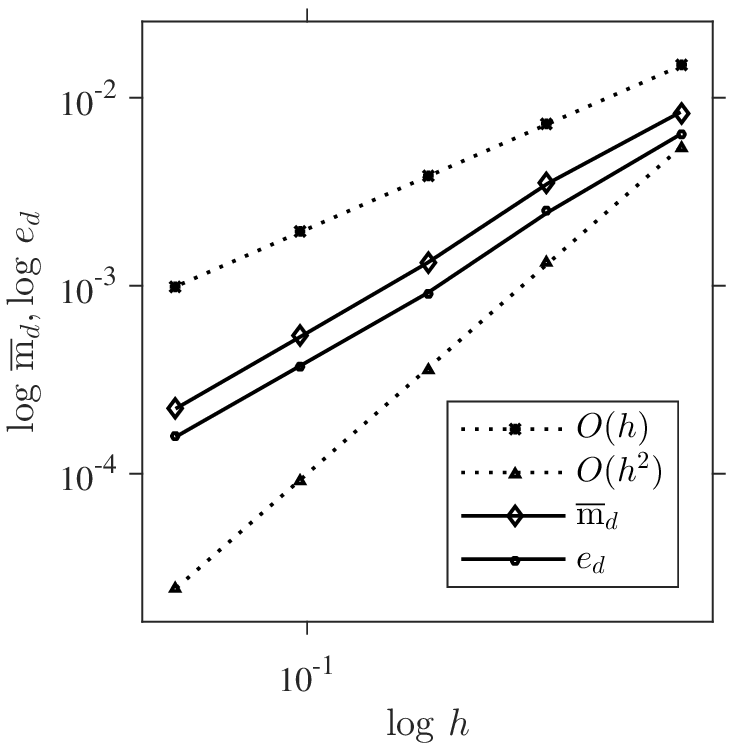}
	\label{fig:example-incr-l-shape-uniform-convergence-ed-md-a}} \quad
	\subfloat[$\error$ and $\maj{}$]{
	\includegraphics[width=6cm]{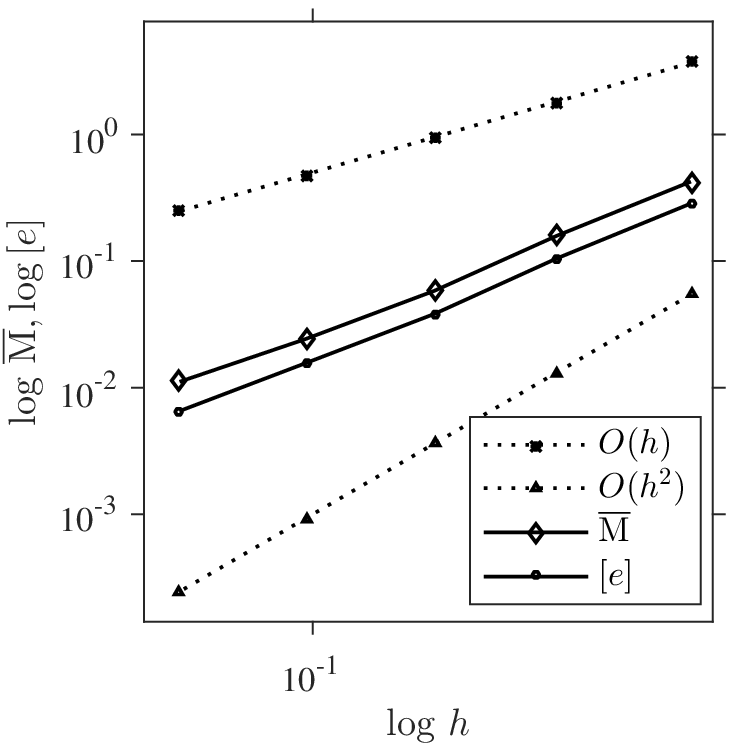}
	\label{fig:example-incr-l-shape-uniform-convergence-ed-md-b}}
	\caption{Ex. \ref{ex:incr-l-shape-2d-t}. 
	The optimal convergence rate of $\ed$ and $\mdI$ (a) and $\error$ and $\maj{}$ (b).}	
	\label{fig:example-incr-l-shape-uniform-convergence-ed-md}
\end{figure}
%
%
%
%
%
%
%
%
%
\begin{figure}[!ht]
	\quad 	
	\centering
	\subfloat[$Q^{(10)}$: 28 EL \newline 
	$\incred{10} = 2.46 \cdot 10^{-3}$,   
	$\incrmdI{10} = 3.50 \cdot 10^{-3}$]{\includegraphics[width=5.4cm]{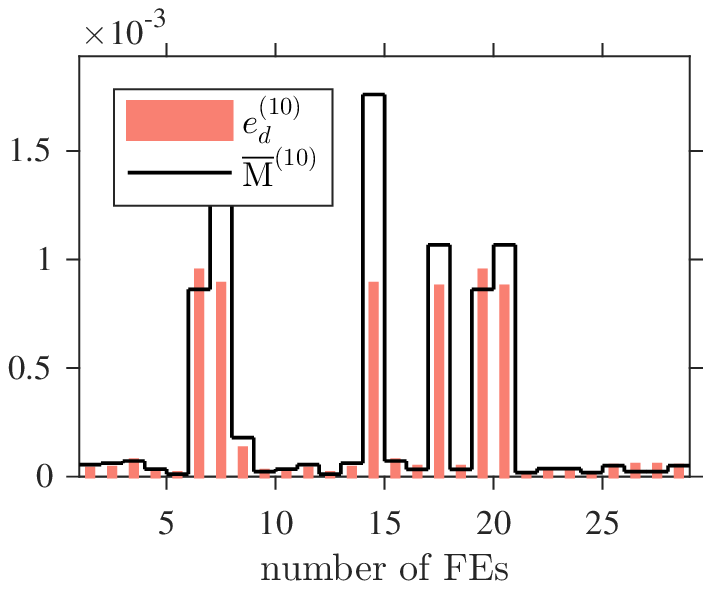}
	\label{fig:example-incr-l-shape-2d-e-maj-distr-a}} \qquad \qquad
	\subfloat[$Q^{(10)}$: 192 EL \newline 
	$\incred{10} = 2.46 \cdot 10^{-3}$,  
	$\incrmdI{10} = 3.50 \cdot 10^{-3}$]{\includegraphics[width=5.2cm]{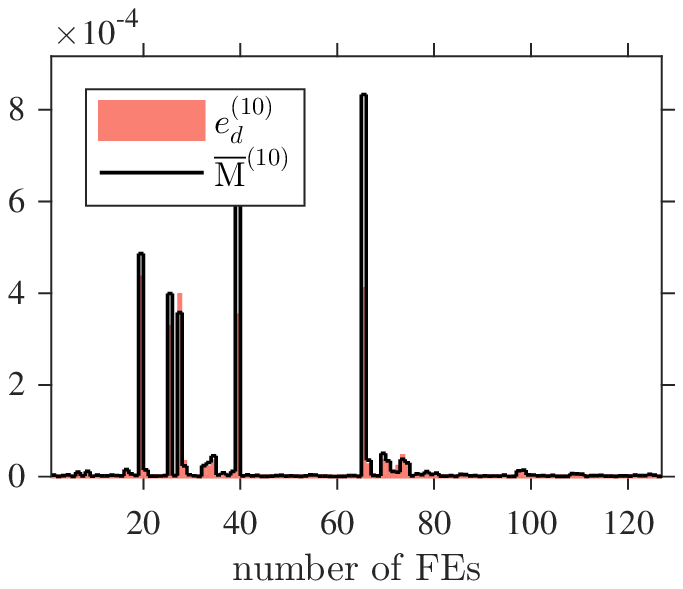}
	\label{fig:example-incr-l-shape-2d-e-maj-distr-b}}
	%
	\caption{Ex. \ref{ex:incr-l-shape-2d-t}. 
	Error and indicator distributions on $Q^{(10)}$, computed on a  
	mesh with (a) 28 EL and (b) 192 EL.}
	\label{fig:example-incr-l-shape-2d-e-maj-distr}
\end{figure}

\begin{figure}[!ht]
\centering
	%
	\subfloat[$Q^{(5)}$: \; 1188 EL, 629 DOF \newline $\incred{5} = 7.78 \cdot 10^{-4}$, $\incrmdI{5} = 1.00 \cdot 10^{-3}$]{\includegraphics[width=6.4cm,  trim={1cm 0cm 1cm 1cm}, clip]{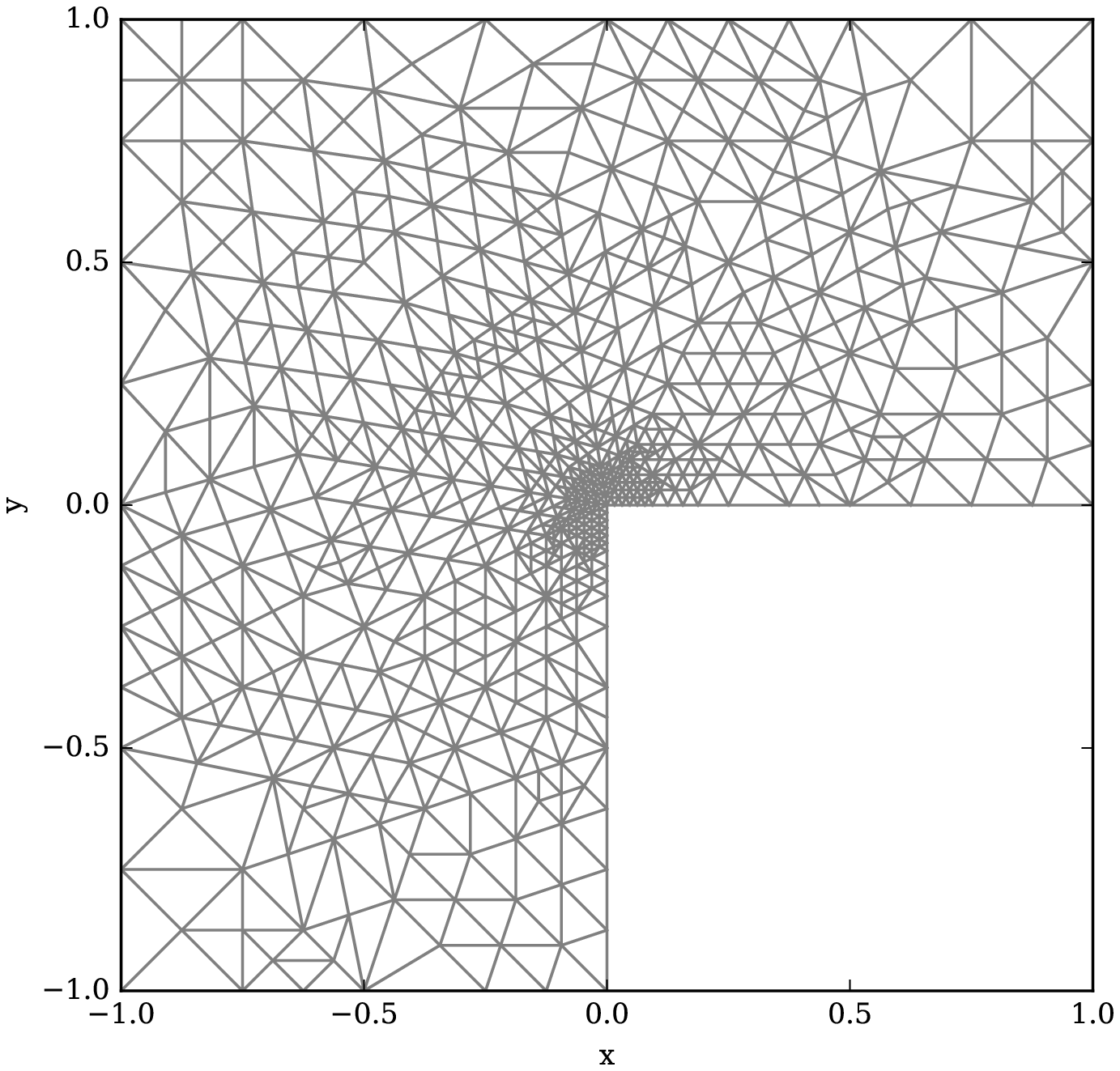}
	}
	\subfloat[$Q^{(5)}$: \; 1022 EL, 543 DOF \newline $\incred{5} =  8.35 \cdot 10^{-3}$, $\incrmdI{5} = 1.07 \cdot 10^{-3}$]{\includegraphics[width=6.4cm,  trim={1cm 0cm 1cm 1cm}, clip]{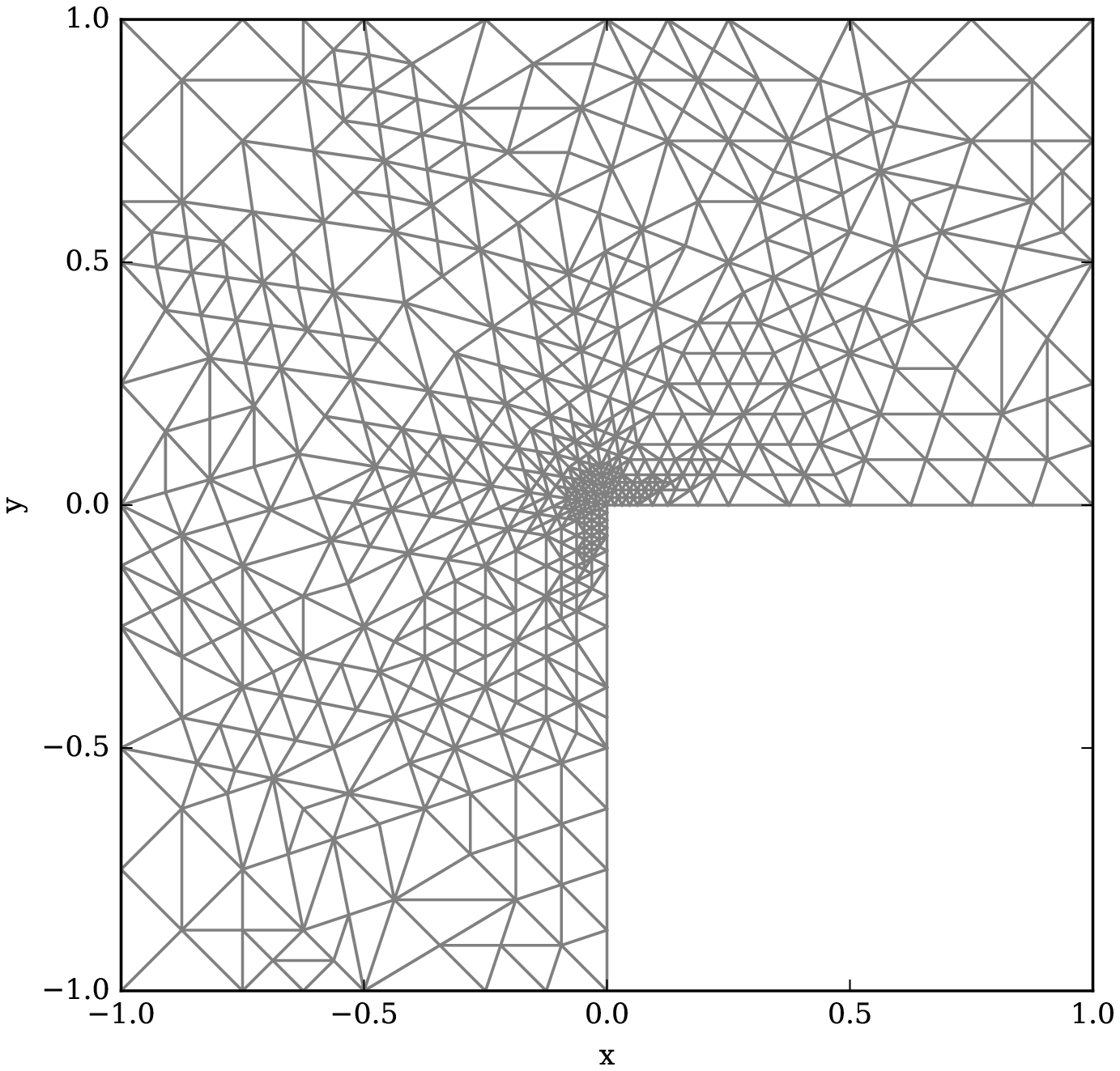}
	}\\[-5pt]
	\subfloat[$Q^{(7)}$: \; 5846 EL, 3008 DOF \newline $\incred{7} = 3.04 \cdot 10^{-4}$, $\incrmdI{7} = 3.43 \cdot 10^{-4}$]{\includegraphics[width=6.4cm,  trim={1cm 0cm 1cm 1cm}, clip]{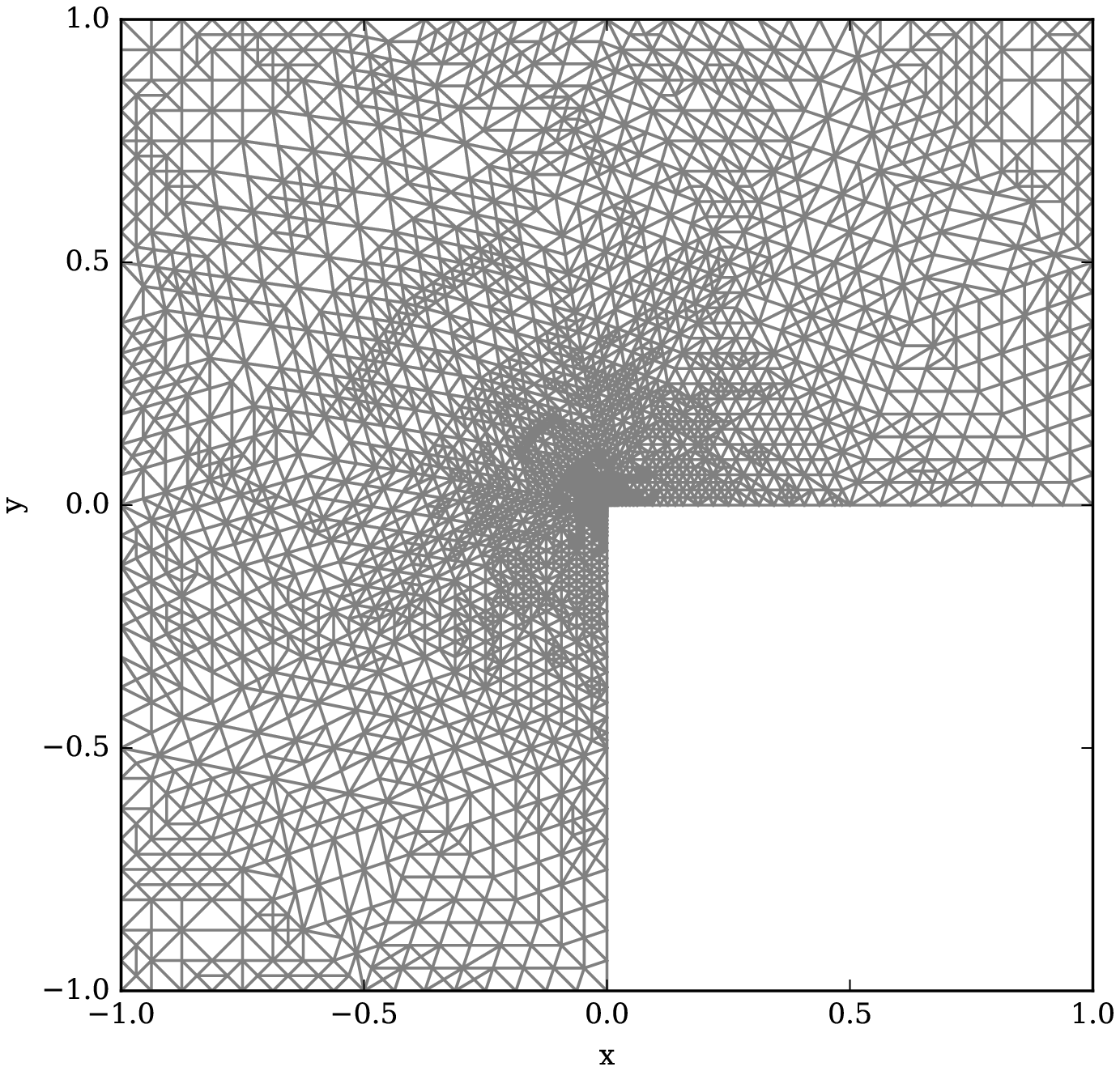}
	}
	\subfloat[$Q^{(7)}$: \; 4676 EL, 2400 DOF \newline $\incred{7} =  3.45 \cdot 10^{-4}$, $\incrmdI{7} = 3.88 \cdot 10^{-4}$]{\includegraphics[width=6.4cm,  trim={1cm 0cm 1cm 1cm}, clip]{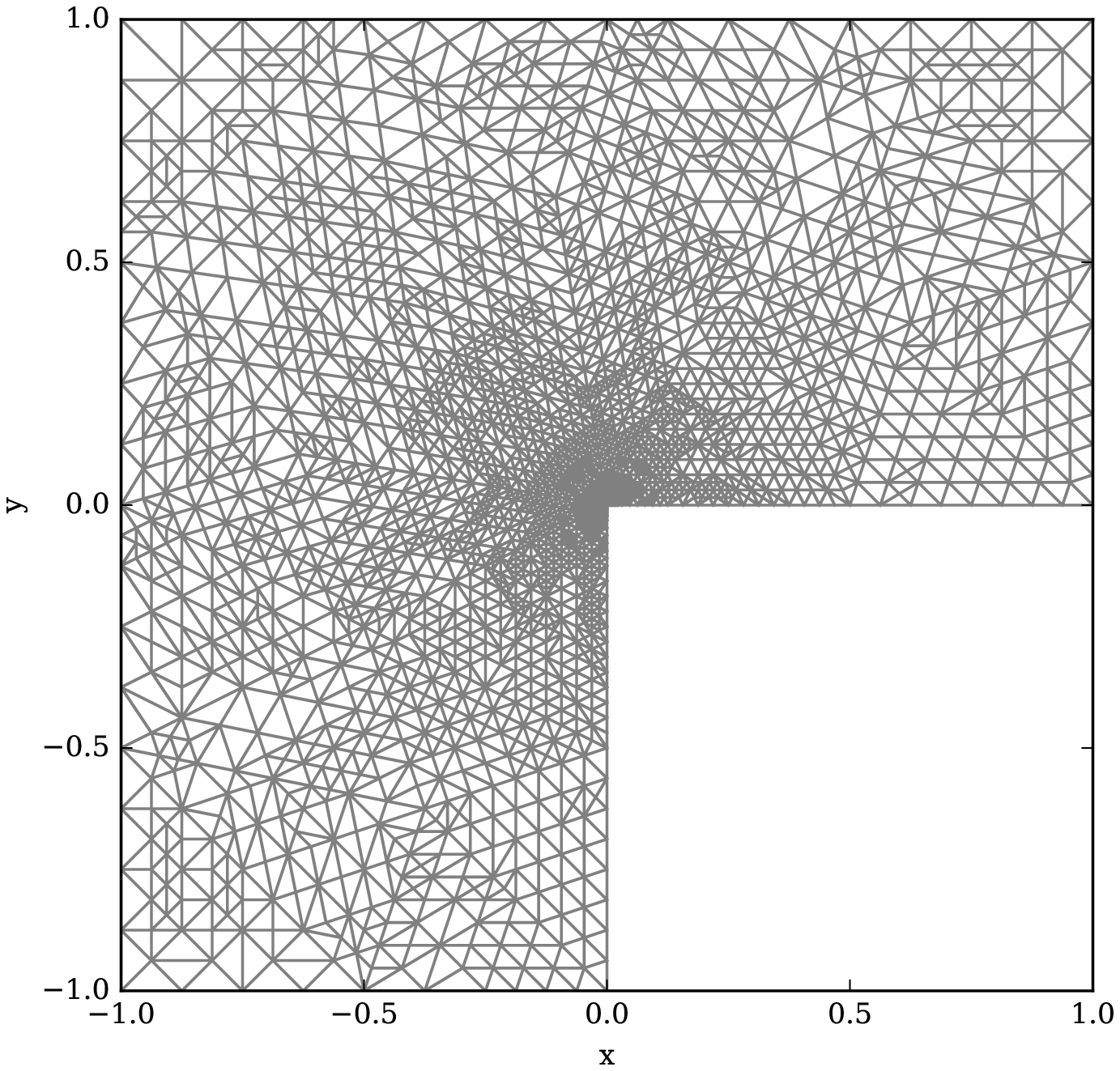}
	}
	\caption{Ex. \ref{ex:incr-l-shape-2d-t}. 
	Evolution of meshes on $Q^{(k)}$, $k = 5, 7$, after the refinements
	based on the true error (a), (c) and on the indicator (b), (d), 
	using the marking $\Marker_{\rm AVR}$.}
	\label{fig:example-incr-l-shape-average-meshes}
\end{figure}

\begin{figure}[!ht]
	\centering
\subfloat[$Q^{(3)}$: \; 470 EL, 256 ND \newline $\incred{3} = 1.88 \cdot 10^{-3}$, $\incrmdI{3} =  2.44 \cdot 10^{-3}$]{\includegraphics[width=6.4cm,  trim={1cm 0cm 1cm 1cm}, clip]{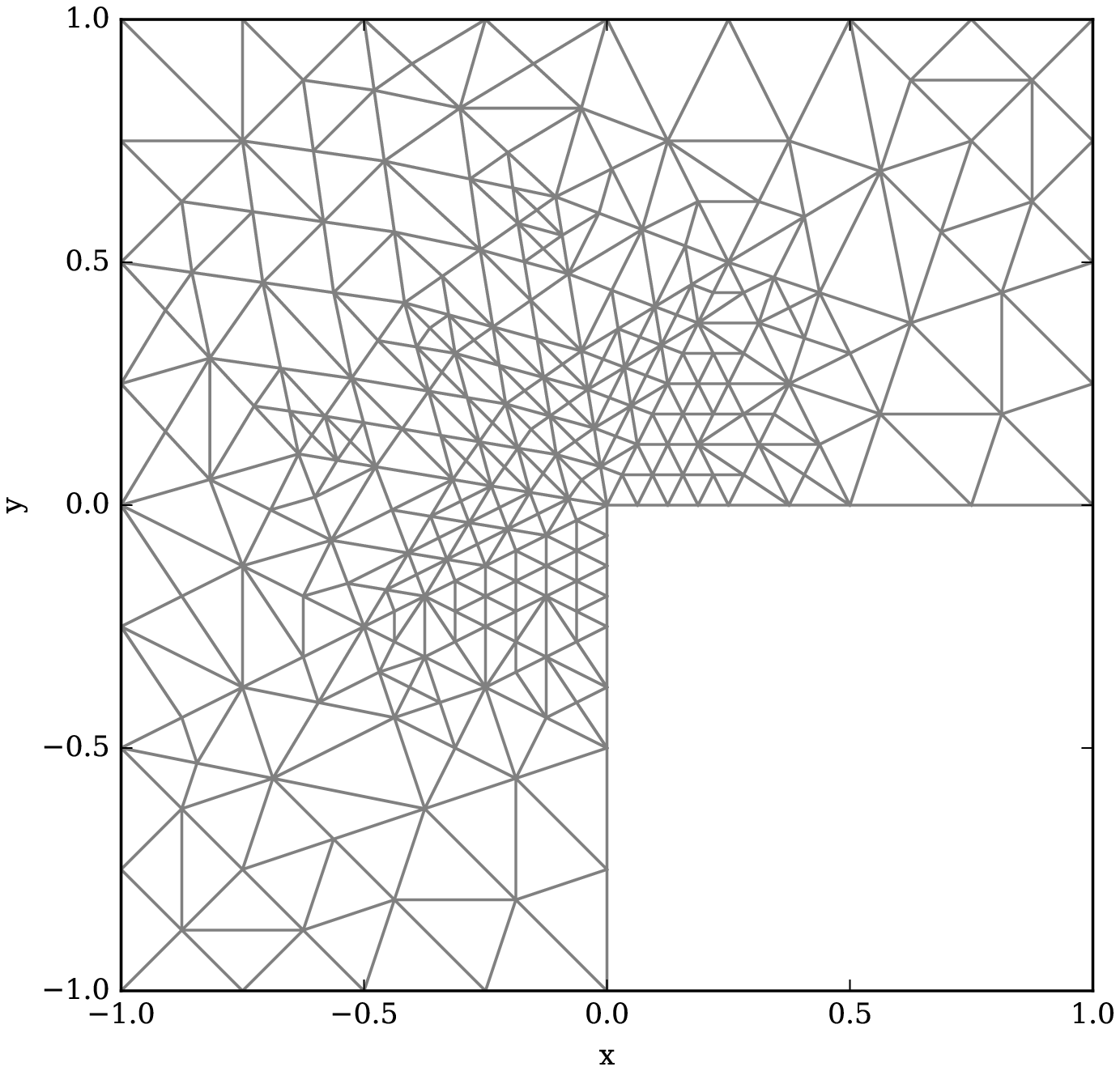}
	}\!
	\subfloat[$Q^{(3)}$: \; 489 EL, 266 ND \newline $\incred{3} = 1.80 \cdot 10^{-3}$, $\incrmdI{3} = 2.35 \cdot 10^{-3}$]{\includegraphics[width=6.4cm,  trim={1cm 0cm 1cm 1cm}, clip]{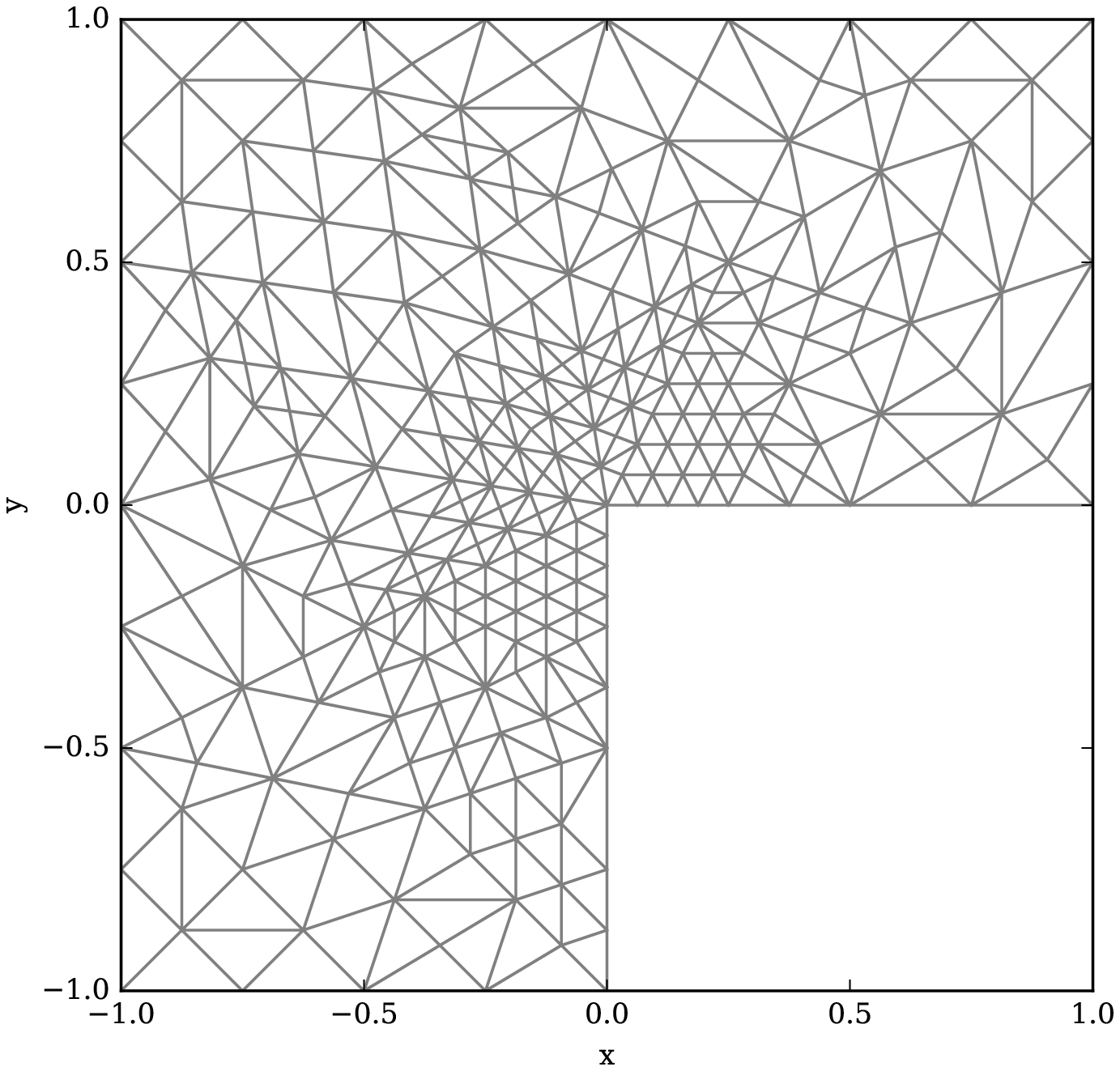}
	}\\[-5pt]
	%
	\subfloat[$Q^{(5)}$: \; 3013 EL, 1563 ND \newline $\incred{5} = 5.34 \cdot 10^{-4}$, $\incrmdI{5} = 7.12 \cdot 10^{-4}$]{\includegraphics[width=6.4cm,  trim={1cm 0cm 1cm 1cm}, clip]{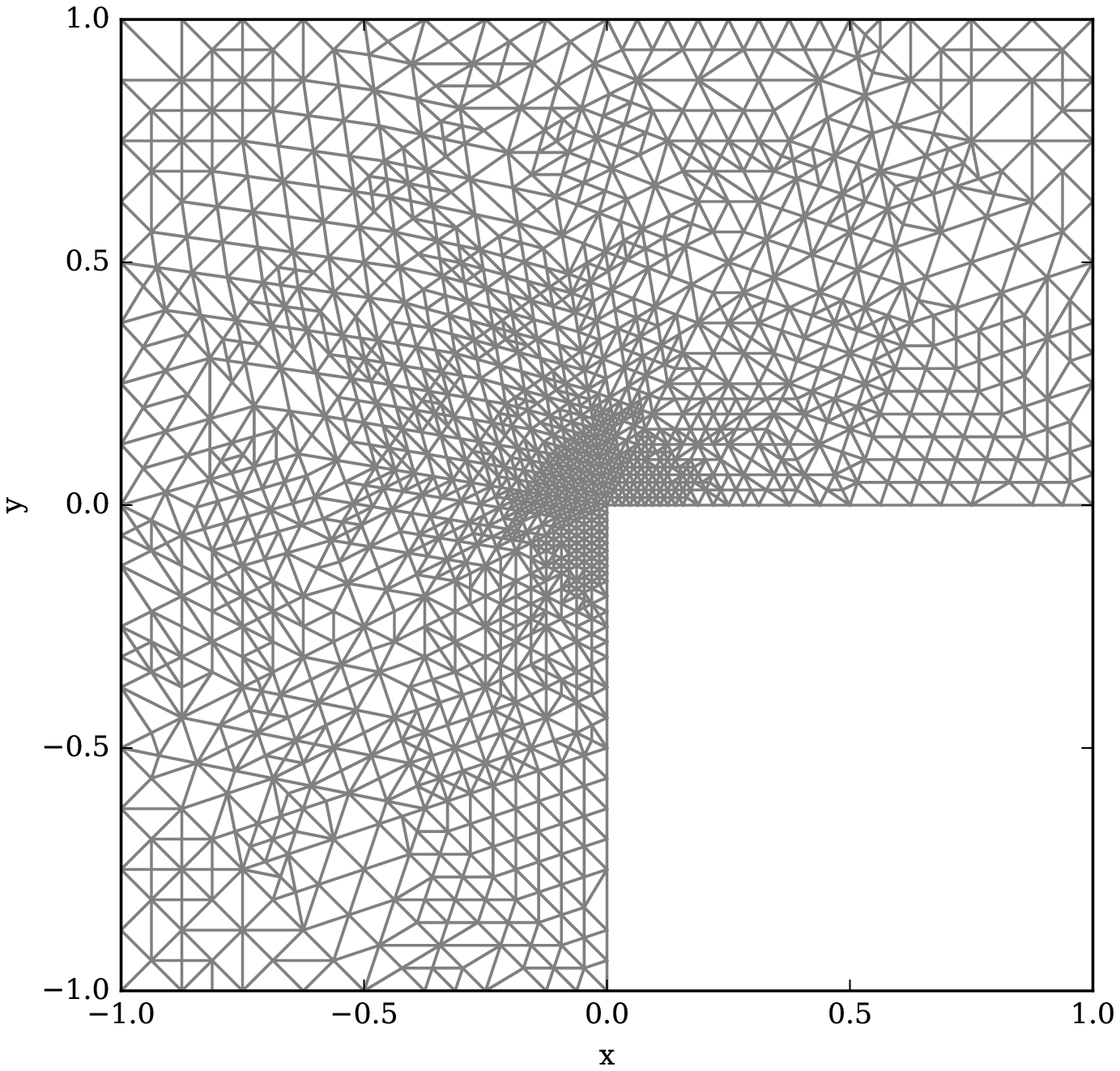}
	}\!
	\subfloat[$Q^{(5)}$: \; 3055 EL, 1579 ND \newline $\incred{5} = 5.36 \cdot 10^{-4}$, $\incrmdI{5} = 7.11 \cdot 10^{-4}$]{\includegraphics[width=6.4cm,  trim={1cm 0cm 1cm 1cm}, clip]{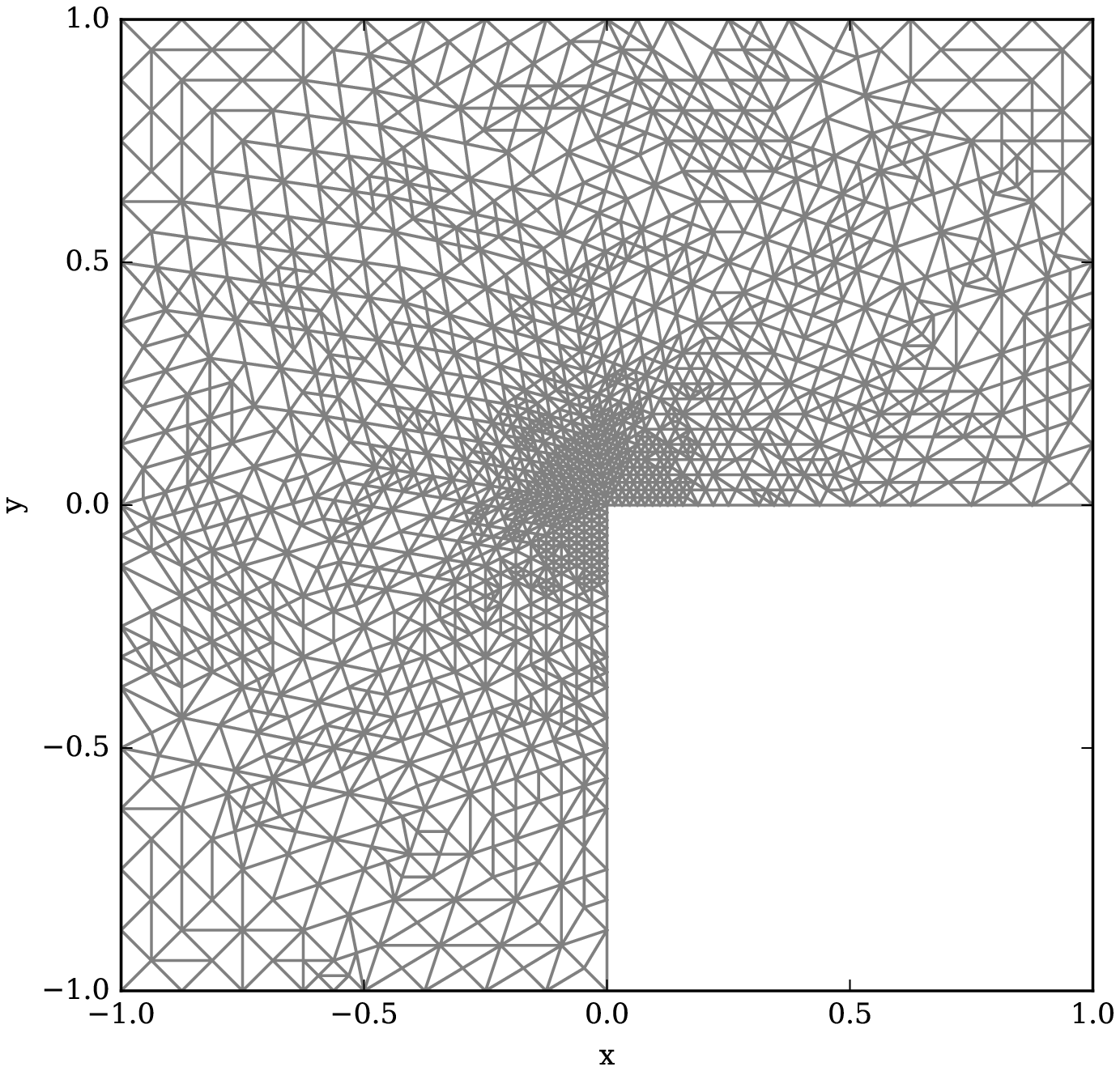}
	}\\
	\caption{Ex. \ref{ex:incr-l-shape-2d-t}. 
	Evolution of meshes on $Q^{(k)}$, $k = 3, 5$, after the refinement
	based on the true error (a), (c) and on the indicator (b), (d), 
	using bulk marking $\Marker_{0.3}$.}
	\label{fig:example-incr-l-shape-bulk-30-meshes}
\end{figure}

The results of the optimal convergence test for the indicator $\mdI$ are provided in 
Figure \ref{fig:example-incr-l-shape-uniform-convergence-ed-md-a} 
(taking into account that  $v \in \Pone$ and $\flux \in \RTone$). Analogously,
we fix the time-step ($K  = 100$) and refine the mesh discretising $\Omega$. 
As expected, the speed of convergence of both the error and the majorant is suboptimal and 
lies between $O(h)$ and $O(h^2)$.  
Figure \ref{fig:example-incr-l-shape-uniform-convergence-ed-md-b} provides the plot
illustrating the
convergence of the total error and the majorant. The difference in the decay of $\mdI$ 
and $\maj{}$ can be explained by the presence of the term $-v_t$ in $\mfI$ 
(the equilibrium part of the majorant) and possible accumulation of the error in the flux $\flux$ 
(in addition to the accumulation of the error in $v$). 

\newpage
Let $K = 10$. Then, the distribution of the local errors on $Q^{(10)}$ is indicated 
quite efficiently by $\incrmdI{10}$ (see Figure 
\ref{fig:example-incr-l-shape-2d-e-maj-distr}). 
Figure \ref{fig:example-incr-l-shape-2d-e-maj-distr-a} provides information about the 
error and majorant distribution on $Q^{(10)}$, where $\Omega$ is discretised by the 
mesh with 28 EL. Figure \ref{fig:example-incr-l-shape-2d-e-maj-distr-b} illustrates 
analogous characteristics for the refined mesh with 192 EL. Both figures confirm that 
$\incrmdI{k}$ manages to locate the errors associated with solution singularities. 

Finally, we consider the adaptive refinement with two marking procedures and 
analyse the obtained meshes. Figure \ref{fig:example-incr-l-shape-average-meshes} 
shows the meshes, obtained by using the marker $\Marker_{\rm AVR}$, and 
Figure \ref{fig:example-incr-l-shape-bulk-30-meshes} compares the meshes derived
using the bulk marking $\Marker_{0.3}$. Analogously to Ex. \ref{ex:incr-unit-2d-t}, 
one can see that the meshes generated during the refinement based on the local true 
error distribution (LHS) are similar to the corresponding meshes provided by adaptive 
algorithm based on the local indicator (RHS). 
\end{example}

In two final examples, we consider rather realistic cases, 
where the exact solution is 
not given explicitly and only the prescribed right-hand side and BC are known. 
In this setting, we compare the approximations reconstructed by 
$\Pone$ FEs to the reference solutions reconstructed with the Lagrangian polynomials 
of degree $4$. Naturally, such a comparison (without an exact solution) has 
only heuristic character, however, this is the only way to estimate the error in 
real-life problems.

\begin{example}
\label{ex:incr-pi-shape-2d-t}
\rm
In the current example, let the domain have $\Pi$-shape, i.e., 
$$\Omega := (- 1, 1) \times (- 1, 1) / \big[-\tfrac{1}{2}, -\tfrac{1}{2}\big] \times [0, -1) 
\in \Rtwo,$$
with the final time fixed to $T = 2$. We consider rather trivial initial and 
boundary conditions $u_0 = 0$ and $u_D = 0$, and the coefficients in 
\eqref{eq:parabolic-equation} are chosen the following way: $A = I$, 
$\vectorb = {\boldsymbol 0}$, $c = 0$. The source function 
$$f = t \sin t \sin \pi x + t \cos t \sin \pi y,$$ rapidly oscillating on $Q$, 
produces approximations changing in time that are
depicted in Figure \ref{fig:ncr-pi-shape-2d-t-approximation}. Here, the
approximate solution is reproduced on the fixed mesh $\mathcal{T}_{h}$ (411 ND) 
and with $K = 15$ time steps.
\begin{figure}[!ht]
	\centering
          \!\!\!\!
          \subfloat[$v^{1} \quad \mbox{on} \quad \Omega_{t^1}$]
          {\includegraphics[width=5.6cm, trim={3cm 1cm 2cm 2cm}, clip]
          {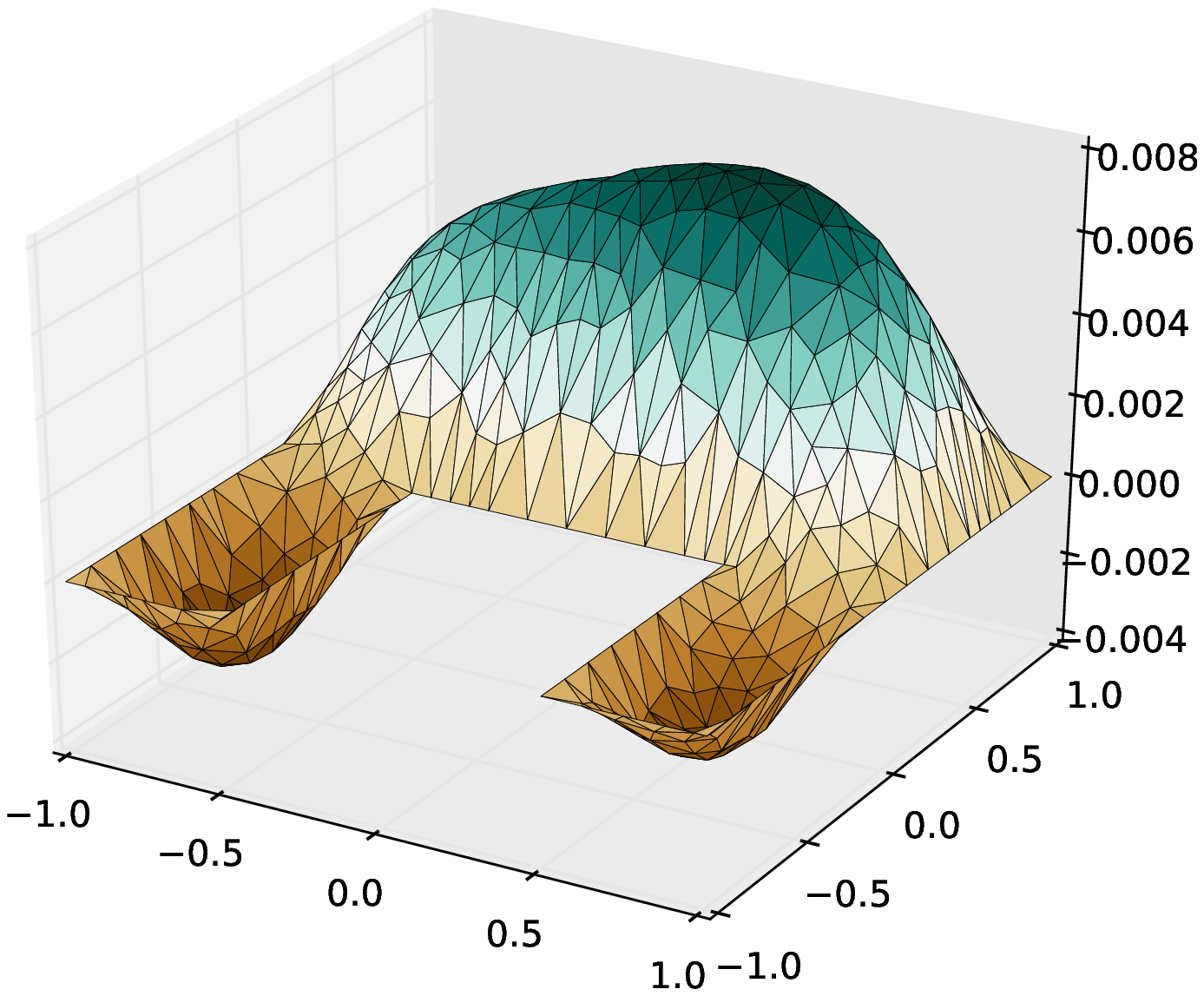}}
	\subfloat[$v^{5} \quad \mbox{on} \quad \Omega_{t^5}$]
          {\includegraphics[width=5.6cm, trim={3cm 1cm 2cm 2cm}, clip]
          {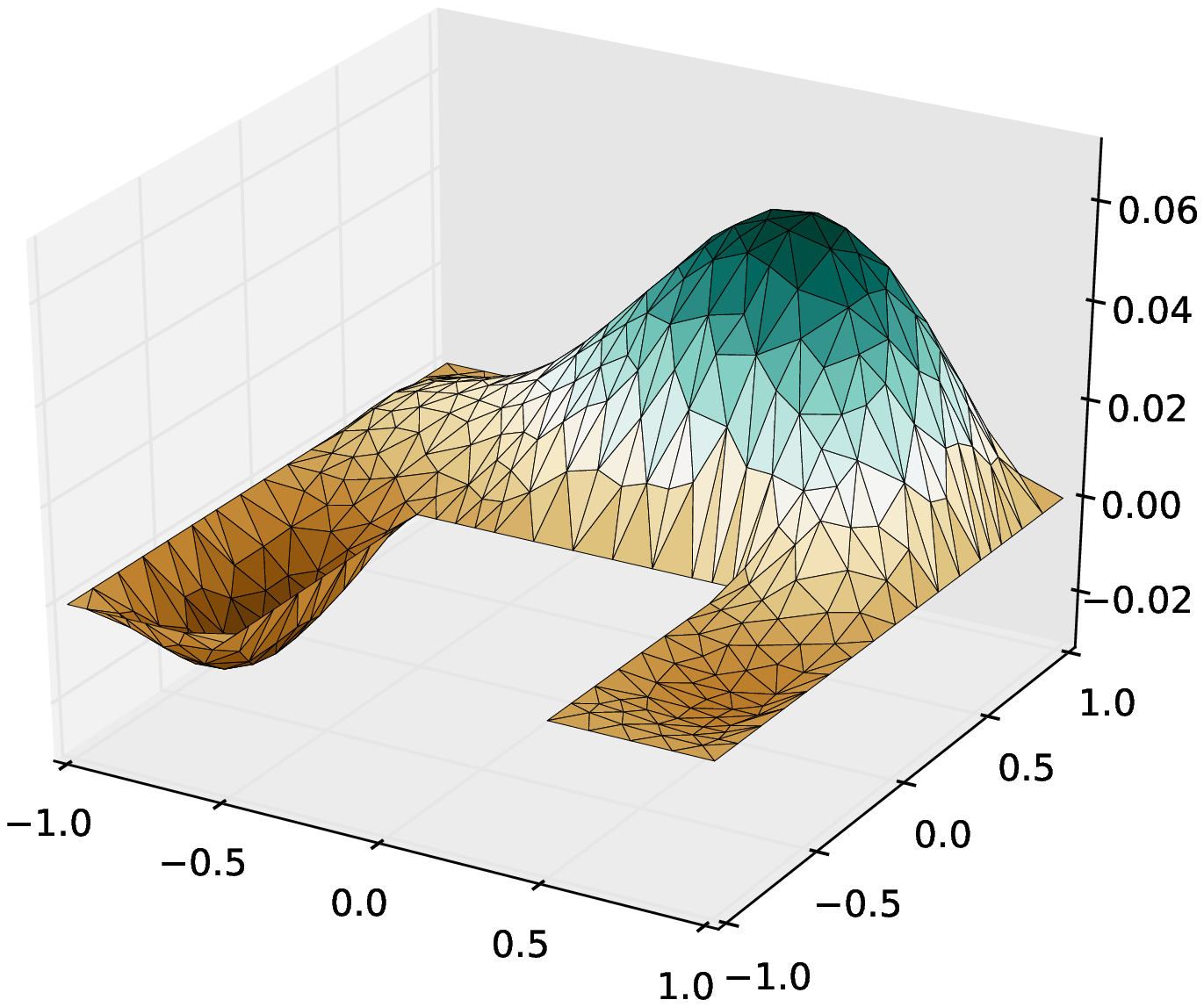}}
          \subfloat[$v^{11} \quad \mbox{on} \quad \Omega_{t^{11}}$]
          {\includegraphics[width=5.6cm, trim={3cm 1cm 2cm 2cm}, clip]
          {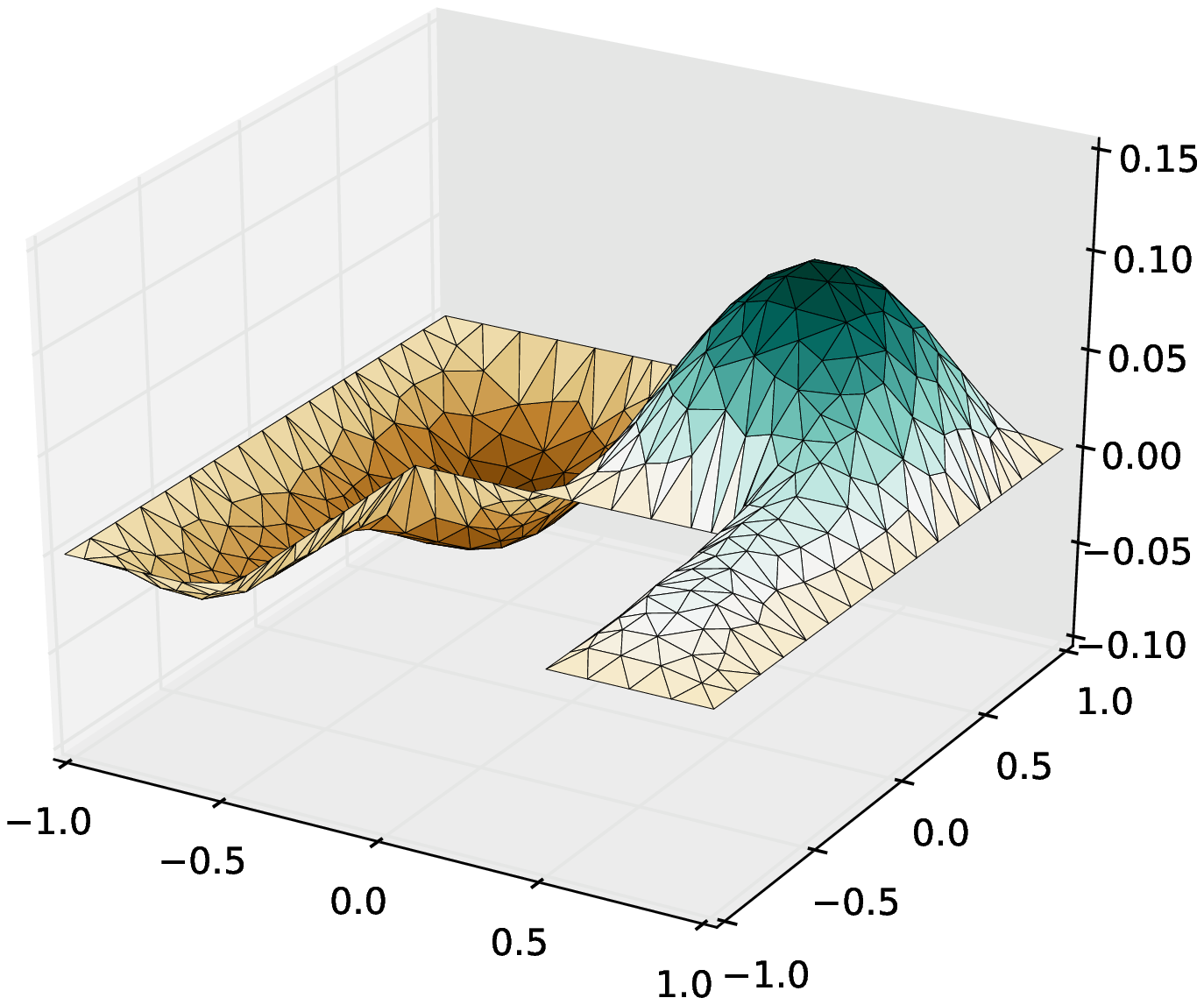}}\!
\caption{Ex. \ref{ex:incr-pi-shape-2d-t}. 
Sequence of approximate solutions reconstructed on the mesh
$\mathcal{T}_{h}$ (411 ND).}
\label{fig:ncr-pi-shape-2d-t-approximation}
\end{figure}


\begin{figure}[!t]
	\centering
	\!\!\!\!
	\subfloat[\qquad $Q^{(0)}$: initial mesh \newline 110 EL, 76 ND]{
	\includegraphics[width=5.6cm, trim={1.8cm 0cm 2.8cm 1.2cm}, clip]
	{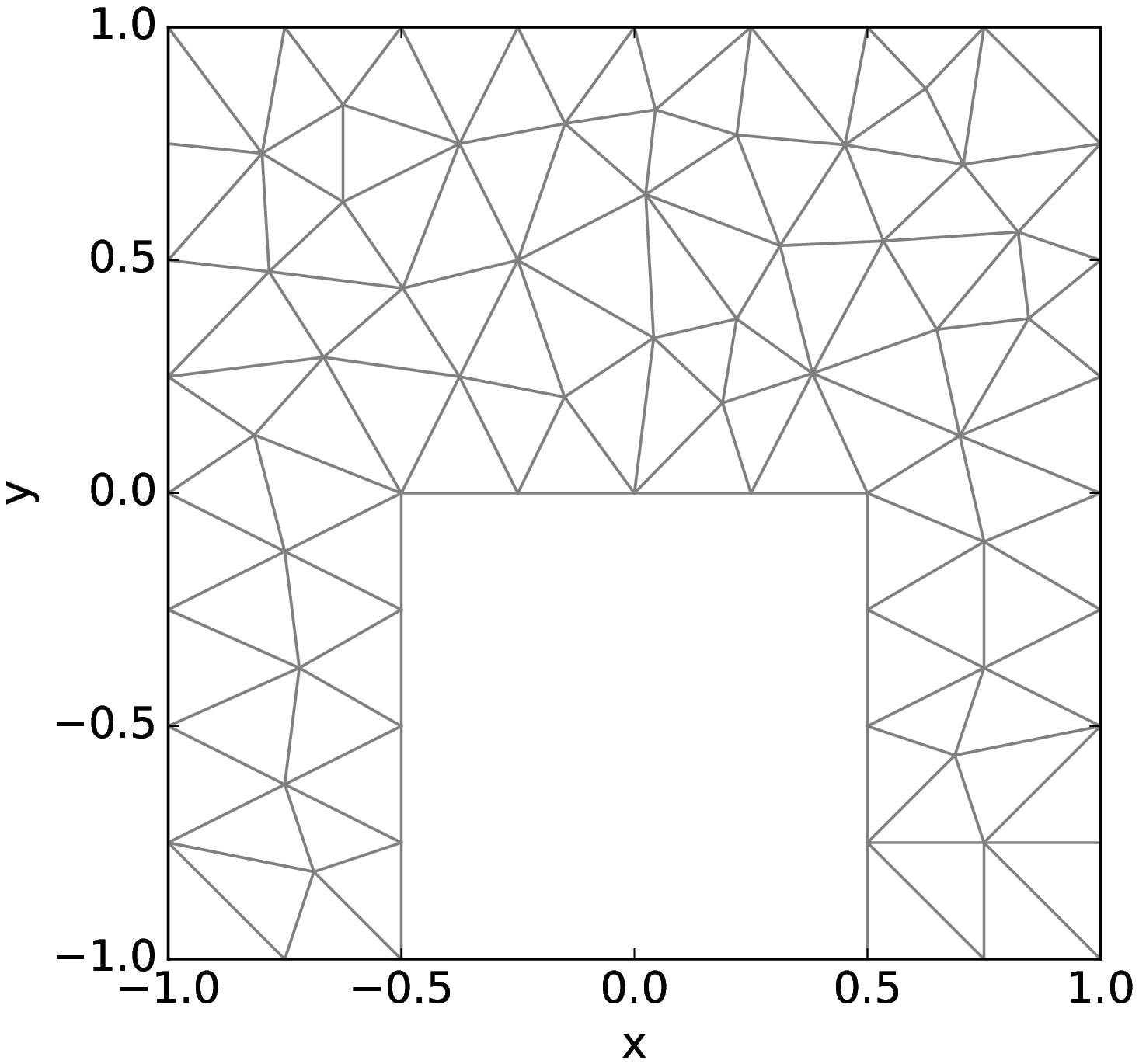}
	\label{eq:incr-pi-shape-2d-t-bulk-10-error-majorant-a}
	} 
	\subfloat[\qquad $Q^{(9)}$: refinement based on $\ed$ 
	\newline 
	$~$ \qquad 5923 EL, 3095 ND]{
	\includegraphics[width=5.6cm, trim={1.8cm 0cm 2.8cm 1.2cm}, clip]{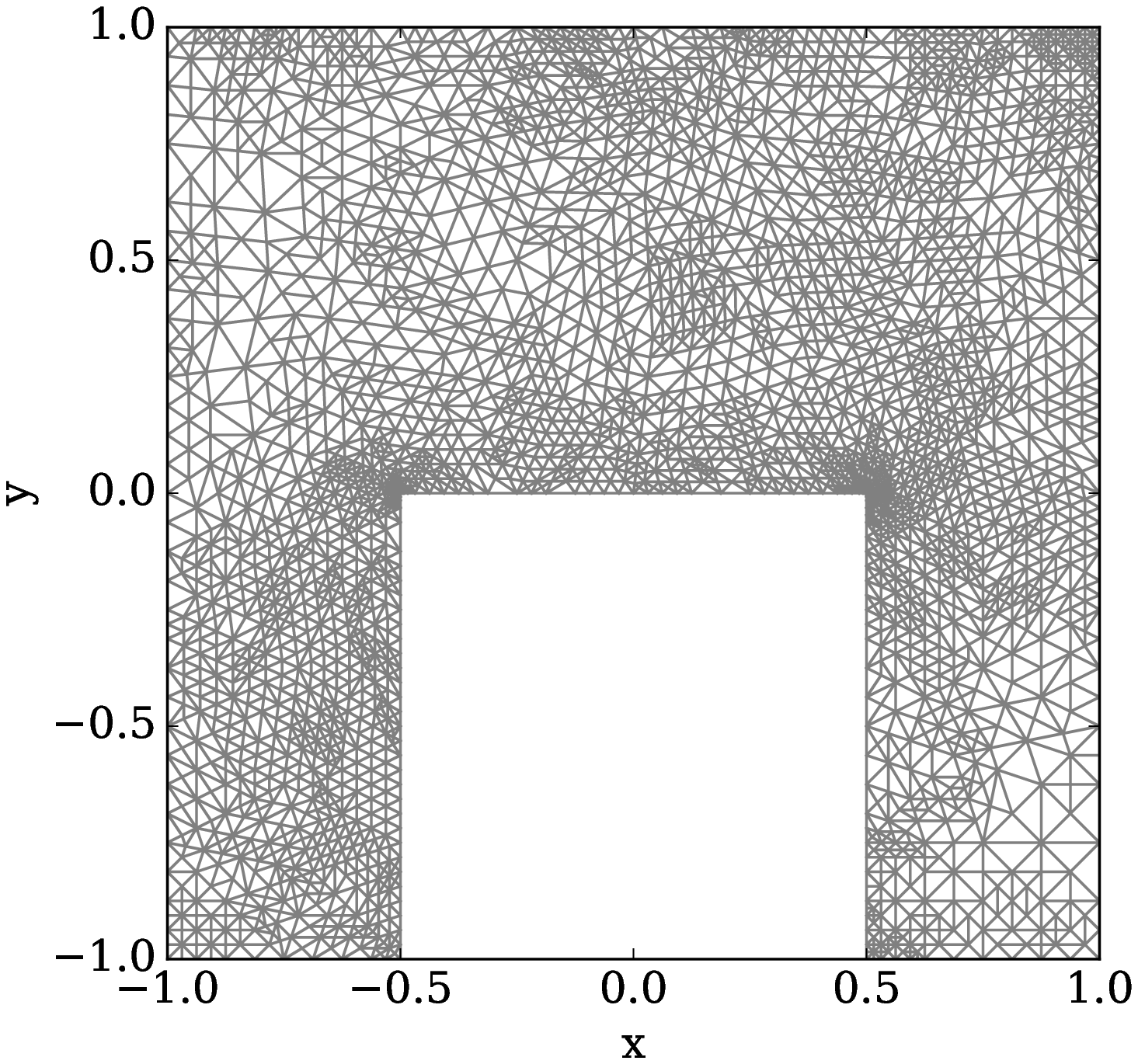}
	\label{eq:incr-pi-shape-2d-t-bulk-10-error-majorant-b}
	}
	\subfloat[\qquad $Q^{(9)}$: refinement based on $\mdI$ 
	\newline
	$~$ \qquad 5828 EL, 3048 ND]{
	\includegraphics[width=5.6cm, trim={1.8cm 0cm 2.8cm 1.2cm}, clip]{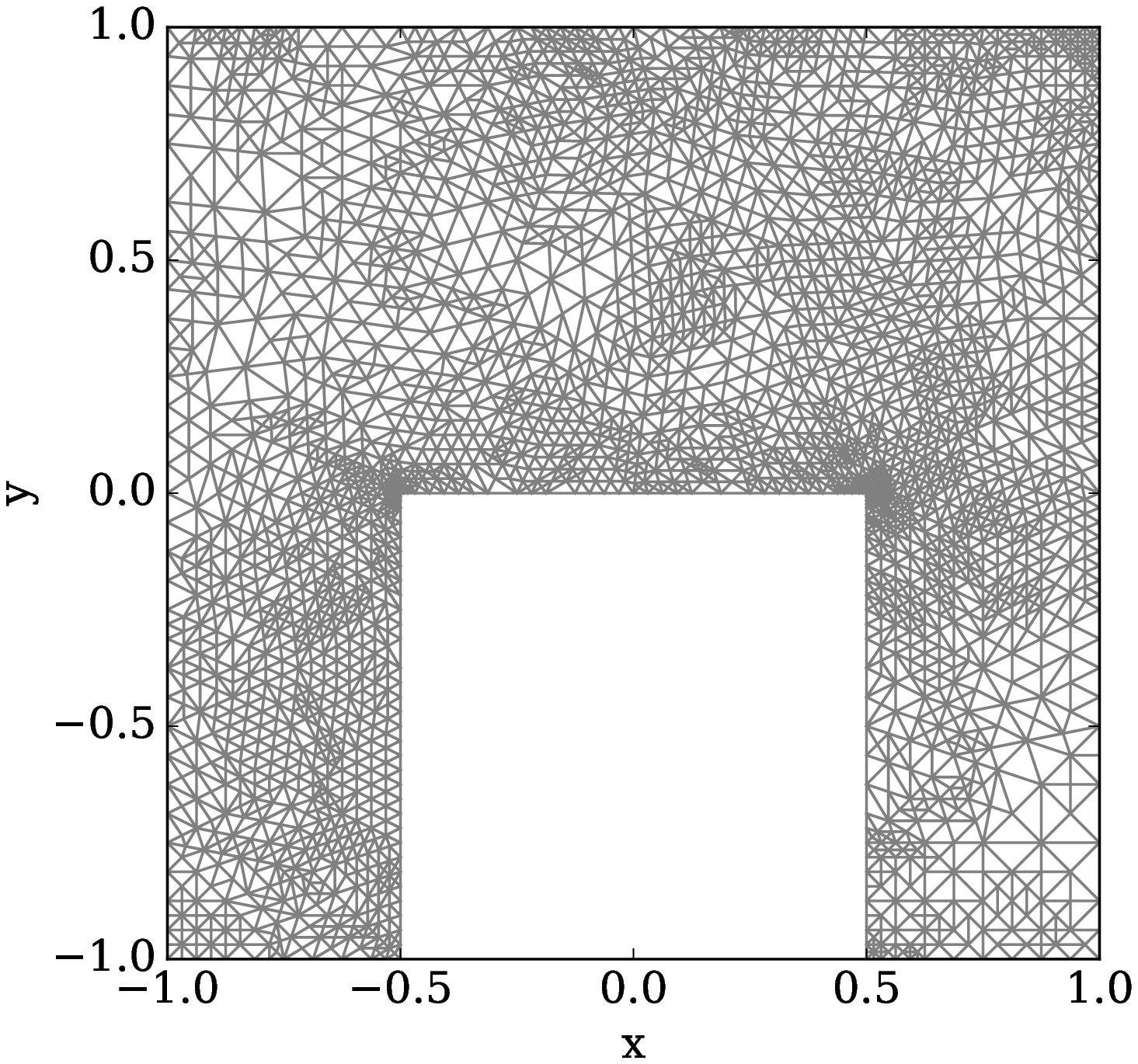}
	\label{eq:incr-pi-shape-2d-t-bulk-10-error-majorant-c}
	}\\[-5pt]
	\subfloat[\qquad $Q^{(12)}$: refinement based on $\ed$ 
	\newline 
	$~$ \qquad 21158 EL, 10827 ND]{
	\includegraphics[width=5.6cm, trim={1.8cm 0cm 2.8cm 1.2cm}, clip]{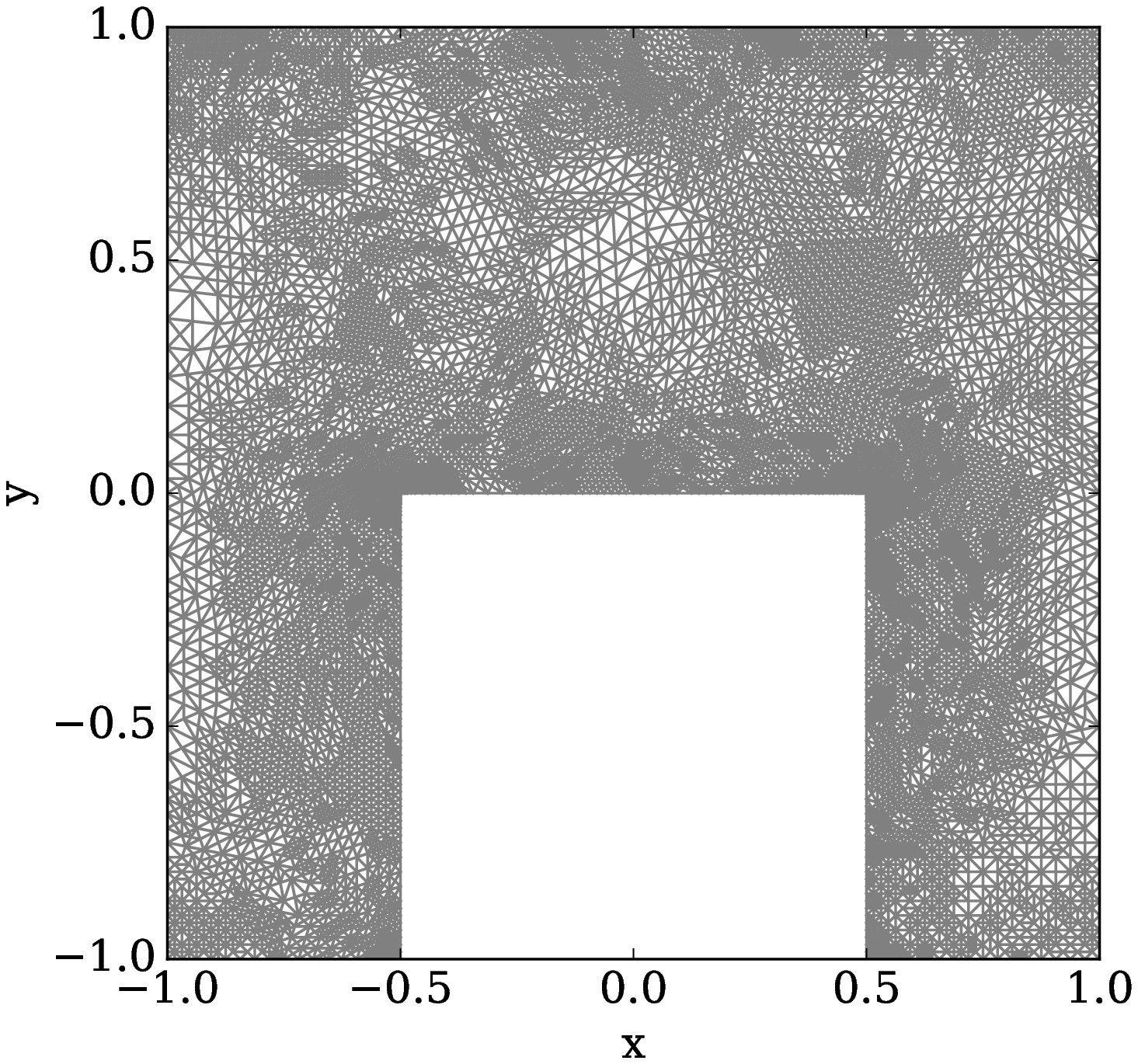}
	\label{eq:incr-pi-shape-2d-t-bulk-10-error-majorant-d}
	} \!\!\!\!
	\subfloat[\qquad $Q^{(12)}$: refinement based on $\mdI$ 
	\newline
	$~$ \qquad 21158 EL, 10828 ND]{
	\includegraphics[width=5.6cm, trim={1.8cm 0cm 2.8cm 1.2cm}, clip]{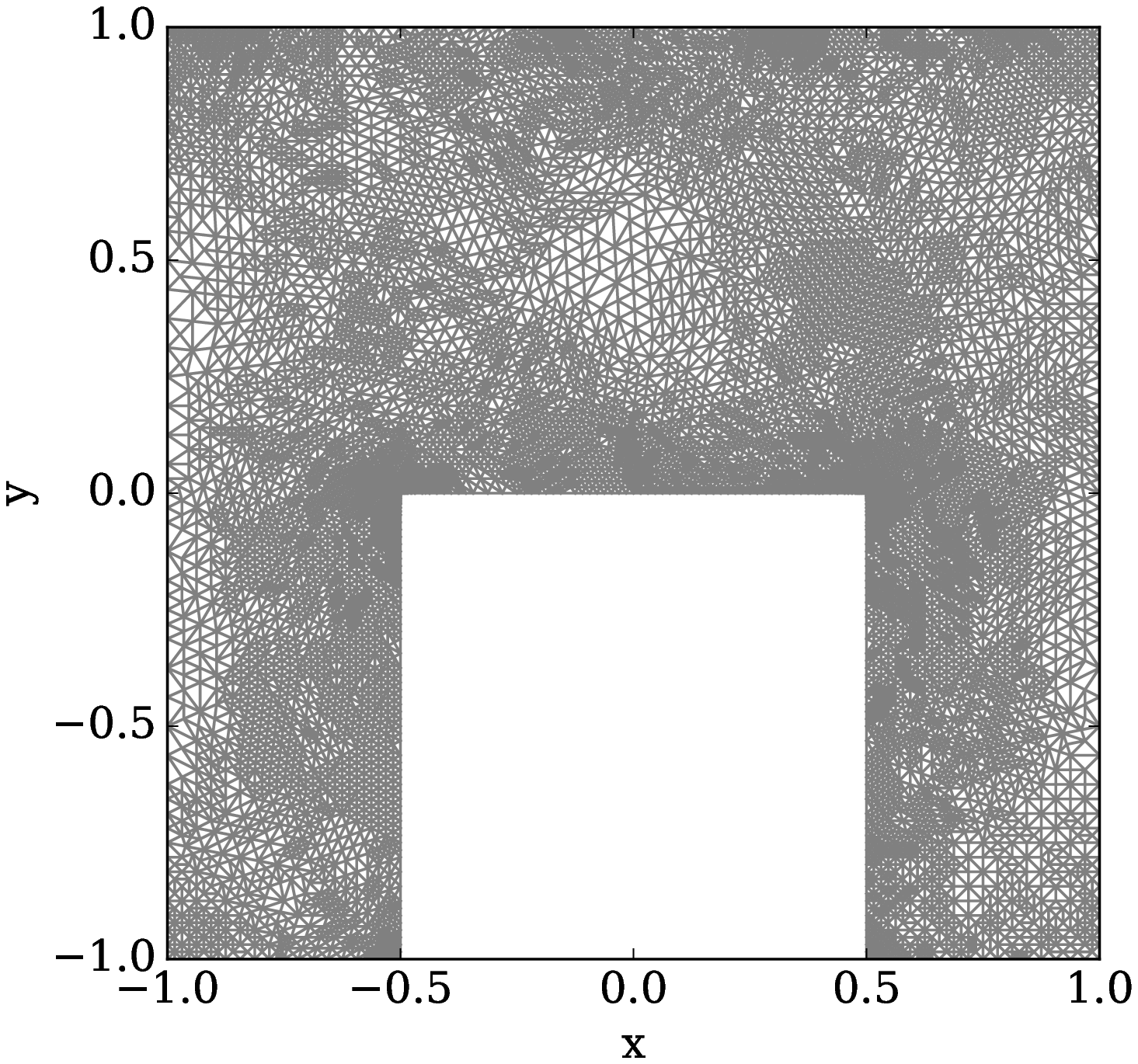}
	\label{eq:incr-pi-shape-2d-t-bulk-10-error-majorant-e}
	}
	\caption{Ex. \ref{ex:incr-pi-shape-2d-t}. 
	Initial mesh (a) and the meshes produced during the refinement on $Q^{(9)}$ and $Q^{(12)}$ 
	time-slabs based either on
	the true error (b), (d) or on the majorant (c), (e), using the bulk marking $\Marker_{0.1}$.}
	\label{eq:incr-pi-shape-2d-t-bulk-10-error-majorant}
\end{figure}

\begin{table}[!t]
\footnotesize
\centering
\begin{tabular}{c|ccc}
$Q^{(k)}$ & $[e]_{\cup_{i = 1}^k {Q^{(i)}}}$ 
& $\overline{\mathrm M}_{\cup_{i = 1}^k {Q^{(i)}}}$ 
& $\Ieff = \overline{\mathrm M}_{\cup_{i = 1}^k {Q^{(i)}}} / [e]_{\cup_{i = 1}^k {Q^{(i)}}}$\\
\midrule
$Q^{(0)}$ &  $4.23 \cdot 10^{-6}$ & $2.99 \cdot 10^{-5}$ &       2.66 \\
$Q^{(2)}$ &  $6.80 \cdot 10^{-5}$ & $1.18 \cdot 10^{-4}$ &       1.32 \\
$Q^{(4)}$ &  $1.90 \cdot 10^{-4}$ & $2.87 \cdot 10^{-4}$ &       1.23 \\
$Q^{(6)}$ &  $3.20 \cdot 10^{-4}$ & $4.68 \cdot 10^{-4}$ &       1.21 \\
 $Q^{(8)}$ &  $4.26 \cdot 10^{-4}$ & $6.21 \cdot 10^{-4}$ &       1.21 \\
 $Q^{(10)}$ &  $5.02 \cdot 10^{-4}$ & $7.37 \cdot 10^{-4}$ &       1.21 \\
$Q^{(12)}$ &  $5.49 \cdot 10^{-4}$ & $8.17 \cdot 10^{-4}$ &       1.22 \\
$Q^{(14)}$ &  $5.76 \cdot 10^{-4}$ & $8.69 \cdot 10^{-4}$ &       1.23 \\
\end{tabular}
\caption{Ex. \ref{ex:incr-pi-shape-2d-t}. Accumulation of the total error and the majorant.}
\label{tab:incr-pi-shape-2d-t-incremented-error-majorant}
\end{table}

As in previous examples, we test the efficiency of the error 
indicator (generated by the majorant) by comparing the meshes produced during the 
adaptive refinement steps performed based either on $\ed$ or on $\mdI$. We use the bulk marking 
criterion with the parameter $\theta = 0.1$ for selecting elements of the mesh 
that need to be refined. In Figure \ref{eq:incr-pi-shape-2d-t-bulk-10-error-majorant-a},
we present the initial mesh that is taken in both refinement 
procedures. Figure \ref{eq:incr-pi-shape-2d-t-bulk-10-error-majorant-b} 
illustrates the resulting mesh that is produced on the time-slice $Q^{(9)}$ after 
using to the true error distribution as the refinement criterion, and, finally, 
Figure \ref{eq:incr-pi-shape-2d-t-bulk-10-error-majorant-c} depicts the mesh 
corresponding to the refinements based on the error indicator element-wise 
distribution. Similar meshes are presented for the time-slice $Q^{(12)}$ in Figures
\ref{eq:incr-pi-shape-2d-t-bulk-10-error-majorant-d} and 
\ref{eq:incr-pi-shape-2d-t-bulk-10-error-majorant-e}.

The total energy error and the majorant (accumulated incrementally as the time 
passes) are presented in Table \ref{tab:incr-pi-shape-2d-t-incremented-error-majorant}. 
The first column contains sequential time-steps, the second one presents the error 
(including the error at previous $k\!-\!1$ steps as well as the increment 
$\incred{k}$ on the $k$-th interval), the third column corresponds to the majorant, 
and the last one shows the majorant efficiency index characterising  
its performance. It is obvious that in the current test-case, 
the efficiency index indicates the adequate behaviour of the majorant even though the 
source term $f$ strongly depends on $t$. This can be explained by looking 
closer at the error and the element-wise indicator distributions (see Figure 
\ref{fig:incr-pi-shape-2d-t-error-majorant-distr}). Here,  
Figures 
\ref{fig:incr-pi-shape-2d-t-error-majorant-distr-a}--\ref{fig:incr-pi-shape-2d-t-error-majorant-distr-d}
in the first row present the local error distributions w.r.t. numbered finite elements (the horizontal axis), 
whereas Figures 
\ref{fig:incr-pi-shape-2d-t-error-majorant-distr-e}--\ref{fig:incr-pi-shape-2d-t-error-majorant-distr-h} 
illustrate the majorant element-wise distribution. 
One can see that the estimate
mimics all local jumps of the errors not only qualitatively but also quantitatively 
sharp. That affects the efficiency index respectively, therefore it stays in the range of 
the interval $[1.21, 1.23]$.

\begin{figure}[!t]
	\centering
	\subfloat[$\incred{2}$]{\includegraphics[width=3.9cm]
	{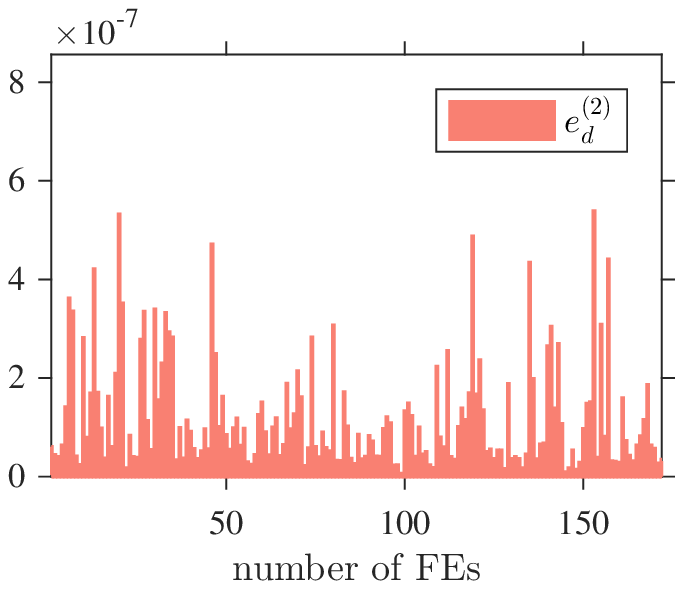}
	\label{fig:incr-pi-shape-2d-t-error-majorant-distr-a}} \,
	\subfloat[$\incred{3}$]{\includegraphics[width=4.0cm]
	{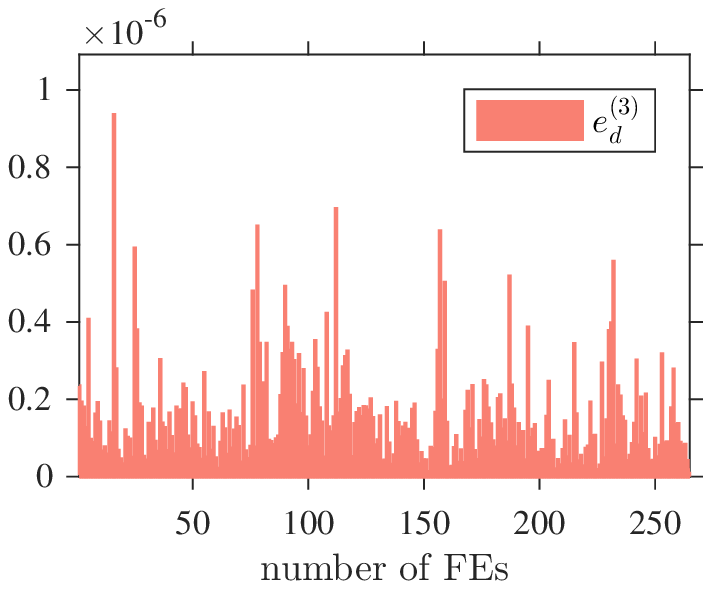}
	\label{fig:incr-pi-shape-2d-t-error-majorant-distr-b}} \,
	\subfloat[$\incred{4}$]{\includegraphics[width=3.9cm]
	{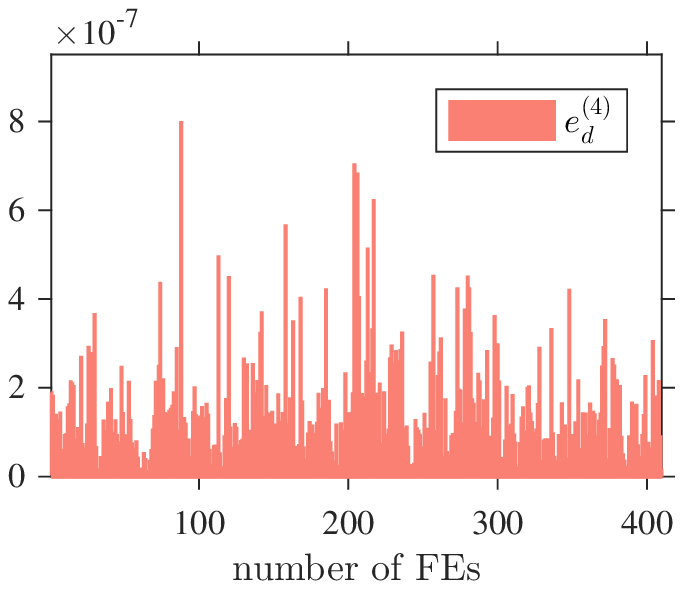}
	\label{fig:incr-pi-shape-2d-t-error-majorant-distr-c}} \,
	\subfloat[$\incred{5}$]{\includegraphics[width=4.1cm]
	{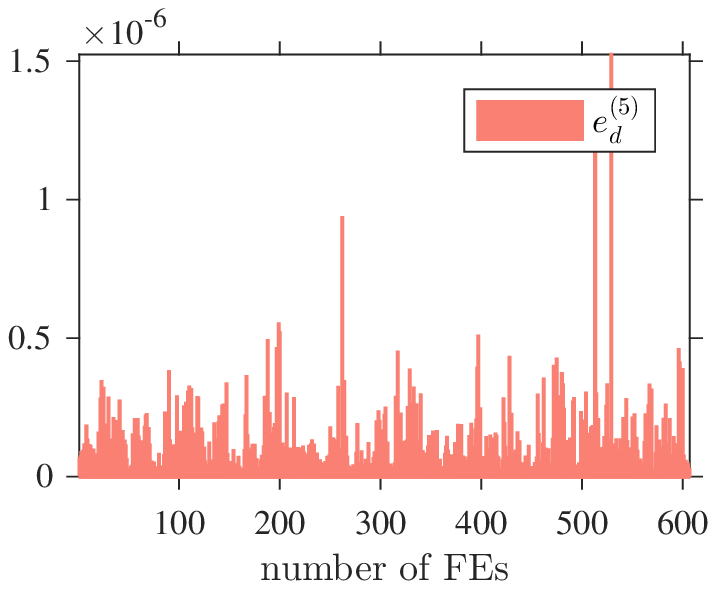}
	\label{fig:incr-pi-shape-2d-t-error-majorant-distr-d}} \,
	\\[-5pt]
	\subfloat[$\incrmdI{2}$]{\includegraphics[width=3.9cm]
	{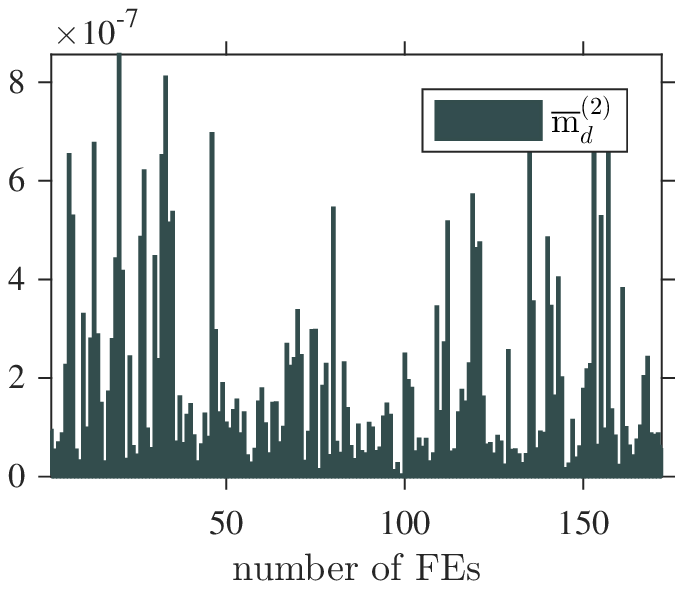}
	\label{fig:incr-pi-shape-2d-t-error-majorant-distr-e}} \,
	\subfloat[$\incrmdI{3}$]{\includegraphics[width=4.0cm]
	{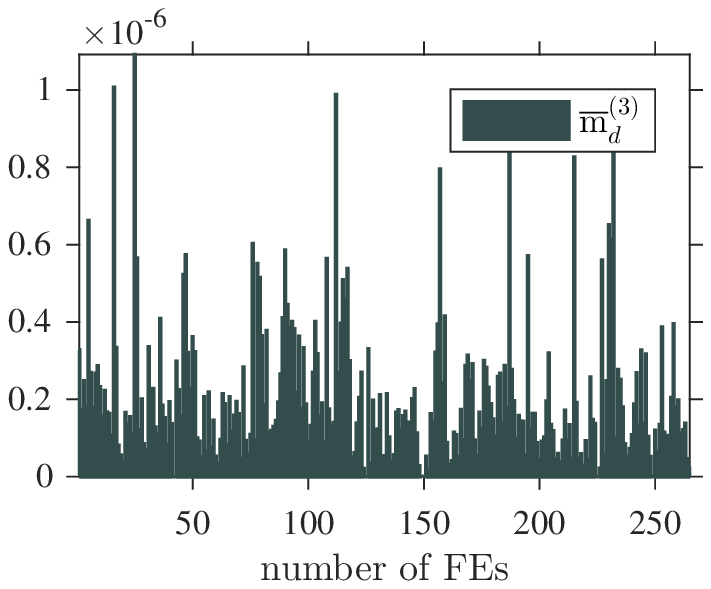}
	\label{fig:incr-pi-shape-2d-t-error-majorant-distr-f}} \,
	\subfloat[$\incrmdI{4}$]{\includegraphics[width=3.9cm]
	{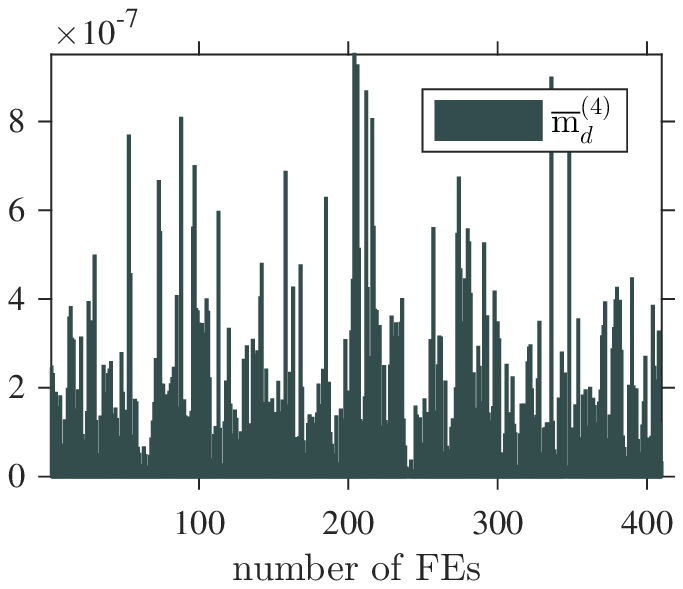}
	\label{fig:incr-pi-shape-2d-t-error-majorant-distr-g}} \,
	\subfloat[$\incrmdI{5}$]{\includegraphics[width=4.1cm]
	{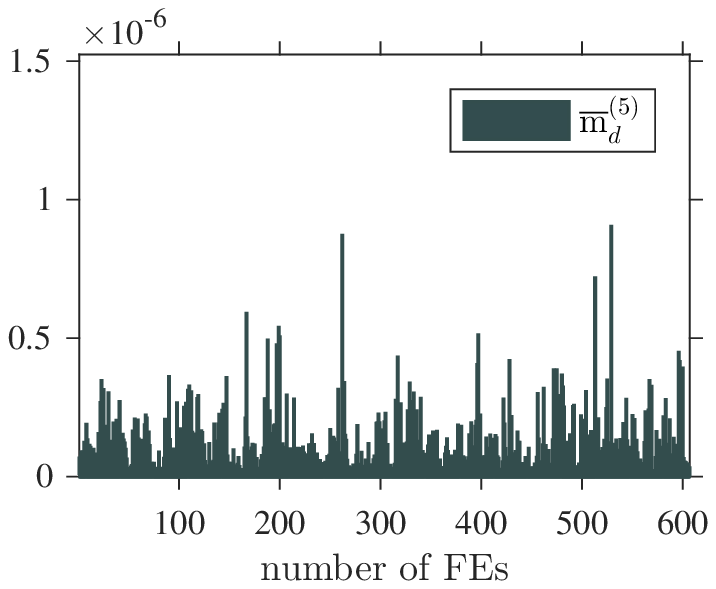}
	\label{fig:incr-pi-shape-2d-t-error-majorant-distr-h}}
	%
	\caption{Ex. \ref{ex:incr-pi-shape-2d-t}. 
	Local distribution of the error $\incred{k}$ and the indicator $\incrmdI{k}$, $k = 2, 3, 4, 5$.}
	\label{fig:incr-pi-shape-2d-t-error-majorant-distr}
\end{figure}

\end{example}

\begin{example}
\label{ex:incr-circ-pi-shape-2d-t}
\rm
In the final example of the current section, we investigate the problem 
with a right-hand side that contains rapidly changing singularity in the source term  
$$f = t\, \sin t \, e^{- 100 \, ((y - 0.8)^2 + (x + 0.8)^2)} 
+ t\, \cos t \, e^{- 100 \, ((y - 0.8)^2 + (x - 0.8)^2)}.$$ 
The domain is defined by a curved $\Pi$-shape domain 
$\Omega := (- 1, 1) \times (0, 1) \backslash B_{\rfrac{1}{2}}(0, 0)$, where 
$B_{\rfrac{1}{2}}(0, 0)$ defines a circle of a radius $\tfrac{1}{2}$ with the centre at 
the point $(0, 0)$. The operator is chosen as
$A = 
\begin{bmatrix}
    1       & 0 \\
    0       & 10 \\
\end{bmatrix}
$, $\vectorb = {\boldsymbol 0}$, $c = 0$, and the initial and boundary conditions 
are zero.  

Because of the nature of FE discretisation and the semicircle-shaped boundary of $\Omega$, 
the spatial mesh does not exactly represents the domain. This implies that the error ($u - v$) 
is not exactly zero in curved part of the boundary. For that reason, the functional majorant is not applicable unless
we project the error on the Dirichlet part of the boundary. However, 
both the majorant and the error indicator can be used in the heuristic sense.

\begin{figure}[!t]
	\centering
	\subfloat[$Q^{(0)}$: initial mesh, 56 EL, 45 ND]{
	\includegraphics[width=6.2cm, trim={1.0cm 2cm 1.0cm 1.2cm}, clip]
	{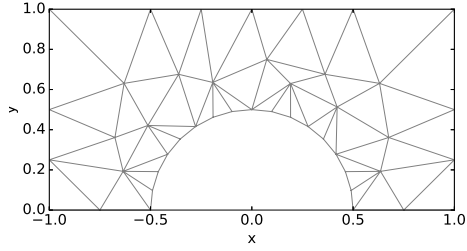}
	\label{eq:incr-pi-shape-2d-t-bulk-10-error-majorant-a}
	}\\[-5pt]
	\subfloat[\qquad $Q^{(10)}$: refinement based on $\ed$ 
	\newline 
	$~$ \qquad 5923 EL, 3095 ND]{
	\includegraphics[width=6cm, trim={1.0cm 2cm 1.0cm 1.2cm}, clip]{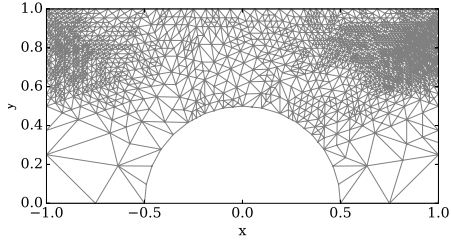}
	\label{eq:incr-pi-shape-2d-t-bulk-10-error-majorant-b}
	}
	\subfloat[\qquad $Q^{(10)}$: refinement based on $\mdI$ 
	\newline
	$~$ \qquad 5828 EL, 3048 ND]{
	\includegraphics[width=6cm, trim={1.0cm 2cm 1.0cm 1.2cm}, clip]{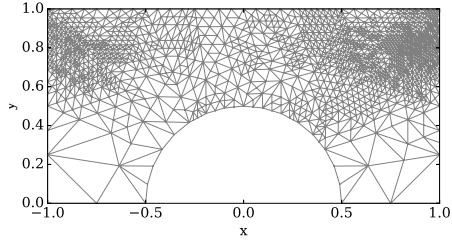}
	\label{eq:incr-pi-shape-2d-t-bulk-10-error-majorant-c}
	}\\[-5pt]
	\subfloat[\qquad $Q^{(13)}$: refinement based on $\ed$ 
	\newline 
	$~$ \qquad 13551 EL, 6881 ND]{
	\includegraphics[width=6cm, trim={1.0cm 2cm 1.0cm 1.2cm}, clip]{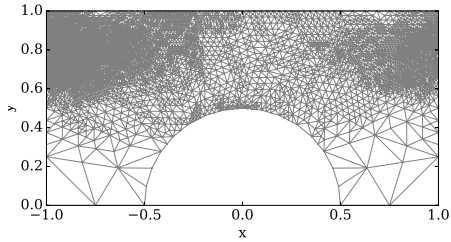}
	\label{eq:incr-pi-shape-2d-t-bulk-10-error-majorant-d}
	}
	\subfloat[\qquad $Q^{(13)}$: refinement based on $\mdI$ 
	\newline
	$~$ \qquad 16664 EL, 8449 ND]{
	\includegraphics[width=6cm, trim={1.0cm 2cm 1.0cm 1.2cm}, clip]{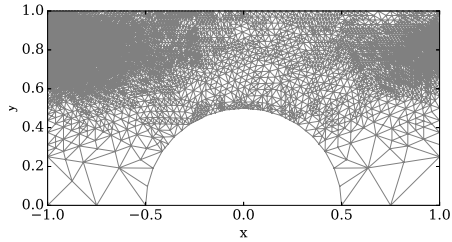}
	\label{eq:incr-pi-shape-2d-t-bulk-10-error-majorant-e}
	}
	\caption{Ex. \ref{ex:incr-circ-pi-shape-2d-t}. 
	Initial mesh (a) and meshes produced during refinement 
	based on the true error (b), (d) and 
	the majorant (c), (e) on the time-slices $Q^{(10)}$ and $Q^{(13)}$,
	using bulk marking $\Marker_{0.1}$.}
	\label{eq:incr-circ-pi-shape-2d-t-bulk-10-error-majorant}
\end{figure}

As it is expected, we see the mesh being refined over time at the singularity points 
of the domain $\Omega$, i.e., $(-0.8, -0.8)$ and $(-0.8, 0.8)$ (Figure 
\ref{eq:incr-circ-pi-shape-2d-t-bulk-10-error-majorant}). Also, 
Table \ref{tab:incr-circ-pi-shape-2d-t-incremented-error-majorant} confirms that 
despite of the singularities in the solution and the non-unitary operator $A$, the efficiency 
of the error estimates stays rather adequate, i.e., $\Ieff \in [1.89, 3.77]$. 

\begin{table}[!ht]
\footnotesize
\centering
\begin{tabular}{c|ccc}
$Q^{(k)}$ & $[e]_{\cup_{i = 1}^k {Q^{(i)}}}$ 
& $\overline{\mathrm M}_{\cup_{i = 1}^k {Q^{(i)}}}$ 
& $\Ieff = \overline{\mathrm M}_{\cup_{i = 1}^k {Q^{(i)}}} / [e]_{\cup_{i = 1}^k {Q^{(i)}}}$\\
\midrule
$Q^{(0)}$ &  $3.45 \cdot 10^{-6}$ & $4.92 \cdot 10^{-5}$ &       3.77 \\
$Q^{(2)}$ & $2.68 \cdot 10^{-5}$ & $1.58 \cdot 10^{-4}$ &       2.43 \\
$Q^{(4)}$ &  $7.36 \cdot 10^{-5}$ & $3.41 \cdot 10^{-4}$ &       2.15 \\
$Q^{(6)}$ &  $1.21 \cdot 10^{-4}$ & $5.45 \cdot 10^{-4}$ &       2.12 \\
$Q^{(8)}$ &  $1.57 \cdot 10^{-4}$ & $6.44 \cdot 10^{-4}$ &       2.02 \\
$Q^{(10)}$ &  $1.84 \cdot 10^{-4}$ & $6.98 \cdot 10^{-4}$ &       1.95 \\
$Q^{(12)}$ &  $2.02 \cdot 10^{-4}$ & $7.39 \cdot 10^{-4}$ &       1.91 \\
$Q^{(14)}$ &  $2.18 \cdot 10^{-4}$ & $7.77 \cdot 10^{-4}$ &       1.89 \\
\end{tabular}
\caption{Ex. \ref{ex:incr-circ-pi-shape-2d-t}. Accumulation of the total error and the majorant.}
\label{tab:incr-circ-pi-shape-2d-t-incremented-error-majorant}
\end{table}

\end{example}

%% file: sections/examples-parabolic-space-time.tex
\section{Space-time approach}
\label{sec:space-time}
 
Due to main drawbacks of the incremental method  (i.e., being time-consuming 
and complicated to parallelise), the space-time FEM approach has been developed. 
The simplest ideas for space-time solvers are based on 
time-parallel integration techniques for ordinary differential equations (ODEs) (an 
overview of the history of this approach can be found in \cite{Gander2015}).
Monograph \cite{Hackbusch1984} introduces a scheme that executes a multigrid 
method for the elliptic problem on each time-step, such that the time is treated as 
an $(d+1)$-th spatial direction in the space-time grid. Later, more space-time 
discretisation methods were suggested, i.e., the so-called parallel time-stepping 
method \cite{Womble1990}, the multigrid waveform relaxation method (space 
parallelism) \cite{VandewallePiessens1992, LubichOstermann1987}, and the full 
space-time multigrid method (using the Fourier mode analysis)
\cite{HortonVandewalle1995}. Recently, due to appearance of powerful multi-core 
computers and advances in parallelised computational methods, application of 
space-time discretisation techniques became popular in various scientific and 
industrial applications (see, e.g., \cite{TakizawaTezduyar2011, Takizawaetall2012, 
TakizawaTezduyar2014, Karabelas2015}, and the references therein). A 
comprehensive review on existing works in space-time techniques can be found in 
\cite{LangerMooreNeumueller2016a}.
Since the majorant is defined as an integral over the whole space-time cylinder $Q$, 
this section contains clarifications on how to apply the 
optimisation algorithm in case of space-time discretisation techniques (see Algorithm 
\ref{alg:majorant-min-space-time}). Numerical results obtained by testing the 
majorant w.r.t. this approach are exposed in 
Ex. \ref{ex:space-time-unit-1d-t} and \ref{ex:space-time-example-2}.

When applying the majorant to the space-time schemes, where time is 
considered as an additional dimension, the optimality condition 
\eqref{eq:majorant-derivative-equal-to-zero} yields the variational formulation on the 
whole $Q$, i.e., 
\begin{equation*}
	\tfrac{\CFriedrichs^2}{\beta \, \underline{\nu}_A} \,
	(\dvrg_x \flux, \dvrg_x \vectorw)_{Q} 
         +  (A^{-1}\flux, \vectorw)_{Q}
	= - \tfrac{\CFriedrichs^2}{\beta \, \underline{\nu}_A} 
	   \big(f - v_t, \dvrg_x \vectorw)_{Q}
	   + (\nabla_x v, \vectorw)_{Q}.
\label{eq:majorant-space-time-variational-formulation}
\end{equation*}
In this case, 
$\flux, \vectorw \in {\rm span} \, \big\{ \, \vectorvarphi_1, ..., \vectorvarphi_M \big\} 
\subset H^{\dvrg_x}(Q)$, i.e., 
$\flux = \sum_{i = 1}^M Y_i \,\vectorvarphi_i$, where $Y \in {\mathds R}^{M}$
is a vector of degrees of freedom (DOFs) approximating $\flux$,
and the test function can be set to $\vectorw = \vectorvarphi_j$, $j = 1, \ldots, M$.
Then, Algorithm \ref{alg:majorant-min-space-time} of minimising the majorant 
follows the steps of Algorithm 3.2 in \cite[Section 3.3.1]{Malietall2014}.
The structure of Algorithm \ref{alg:majorant-min-space-time} is similar to the one 
for the incremental scheme with exception that refinement steps are carried out on 
the whole space-time domain $Q$ instead of the time-slices $Q^{(k)}$. In the current case, the 
adaptivity in time direction is automated since $t$ is treated as $x_{d+1}$-coordinate.

\begin{algorithm}[!t]
\caption{\quad Global minimisation of $\maj{}$ (in the case of a space-time scheme)}
\label{alg:majorant-min-space-time}
\begin{algorithmic} 
\STATE {\bf Input:} $Q$: $v$ 
\COMMENT{approximate solution}
\STATE $\quad \qquad$ $\vectorvarphi_i$, $i = 1, \ldots, M$ 
\COMMENT{$H^{\dvrg_x}(Q)$-conforming basis functions}
\STATE $\quad \qquad$ $L^{\rm iter}_{\rm max}$ 
\COMMENT{number of iterations}
\STATE
\STATE Assemble the matrices $\widetilde{S}$, 
$\widetilde{K} \in {\mathds R}^{M \times M}$ and the vectors 
$\widetilde{z}$, $\widetilde{g} \in {\mathds R}^{M}$ by using 
			\begin{alignat*}{2} \!\!
				\{ \widetilde{S}_{ij} \}_{i, j=1}^M  & = 
				(\dvrg_x \vectorvarphi_i, \dvrg_x \vectorvarphi_j)_Q, \; \quad
				\{ \widetilde{z}_{j} \}_{j=1}^M = 
				\big( f - v_t,  \dvrg_x \vectorvarphi_j \big)_Q,  
				\\
				\{ \widetilde{K}_{ij} \}_{i, j=1}^M & = 
				(A^{-1} \vectorvarphi_i, \vectorvarphi_j)_Q, \qquad \; \quad
				\{ \widetilde{g}_{j} \}_{j=1}^M = 
				\big(\nabla_x v, \vectorvarphi_j \big)_Q.
			\end{alignat*}
\STATE Initialise $\beta$, e.g., $ \beta = 1$.
\vspace{4pt}
\FOR{$l = 1$ {\bf to} $L^{\rm iter}_{\rm max}$}
\vspace{4pt}
\STATE Solve the system \quad 
$
\left( \tfrac{\CFriedrichs^2}{\beta \, \underline{\nu}_A} \widetilde{S} + \widetilde{K} \right)\, {Y} = 
	-\tfrac{\CFriedrichs^2}{\beta \, \underline{\nu}_A} \widetilde{z} + \widetilde{g}.
$
\STATE Approximate the flux $\flux = \sum\limits_{i=1}^M Y_i \vectorvarphi_i$.
\STATE Compute the components of the majorant by using 
$$\mfI := \| \, f + \dvrg_x \flux - v_t  \, \|^2_{Q} \quad \mbox{and}  \quad \mdI := \| \flux - A \nabla_x v \,\|^2_{A^{-1}, Q}.$$
\STATE Compute the optimal $\beta$ by using
$
\beta := 
\left( \tfrac{\CFriedrichs^2 \mfI }{ \underline{\nu}_A \mdI} \right)^{\rfrac{1}{2}}.
$
\ENDFOR
\vspace{5pt}
\STATE Compute the majorant by using
\begin{equation*}
\maj{} (v, \flux; \beta) 
	:= (1 + \beta) \, \| \flux - A \nabla_x v \,\|^2_{A^{-1}, Q} + 
   	    (1 + \tfrac{1}{\beta}) \, \tfrac{\CFriedrichs^2}{\underline{\nu}_A}
		\| \, f + \dvrg_x \flux - v_t  \, \|^2_{Q}.
		\label{eq:majorant-simplified-total-time}
\end{equation*}
\STATE {\bf Output:} $\maj{} (v, \flux; \beta)$ \COMMENT{majorant on $Q$}
\end{algorithmic}
\end{algorithm}

\begin{example}
\label{ex:space-time-unit-1d-t}
\rm
First, let us consider the numerical properties of $\maj{}$ and the corresponding 
indicator on the unit interval $\Omega = (0, 1) \subset \Real$, $T = 1$, $A = I$, 
$\vectorb = {\boldsymbol 0}$, $c = 0$, and  
homogeneous Dirichlet BC. The exact solution is $u = x\,(1 - x)\,(t^2 + t + 1)$ with 
the IC being  
$u_0 = x\,(1 - x)$. The approximation is reconstructed by $\Pone$ FEs, and the flux 
is approximated by $\Ptwo$ FEs. In Figure 20, we confirm the order of convergence 
$O(h^2)$ for $\maj{}$ after several simultaneous refinement iterations (w.r.t. space and time).


\begin{figure}[!ht]
	\centering
	\vskip 10pt
	\includegraphics[width=5.4cm]{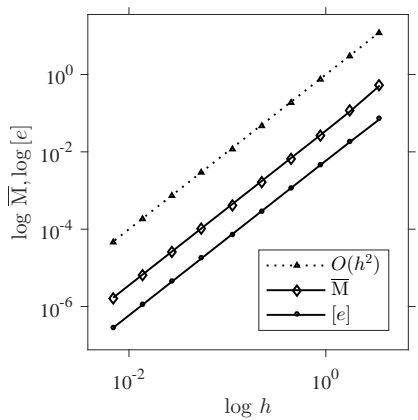}
	\label{fig:example-space-time-1d-t-unit-opt-convergence}
	\caption{Ex. \ref{ex:space-time-unit-1d-t}.
	The optimal convergence of the total error and majorant. }
\end{figure}

\begin{figure}[!ht]
	\centering
	\subfloat[2 REF: 32 EL, \newline $\ed = 9.11 \cdot 10^{-2}, \mdI = 9.16 \cdot 10^{-1}$]{
	\includegraphics[width=5.3cm]{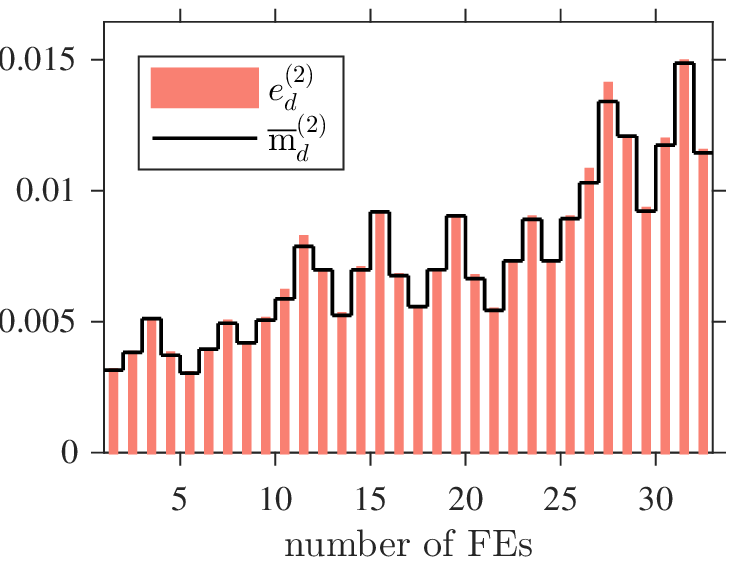}\label{fig:example-space-time-1d-t-unit-e-maj-ref-2}}\!	
	\subfloat[3 REF: 128 EL, \newline $\ed = 2.30 \cdot 10^{-2}, \mdI = 2.30 \cdot 10^{-2}$]{
	\includegraphics[width=5.2cm]{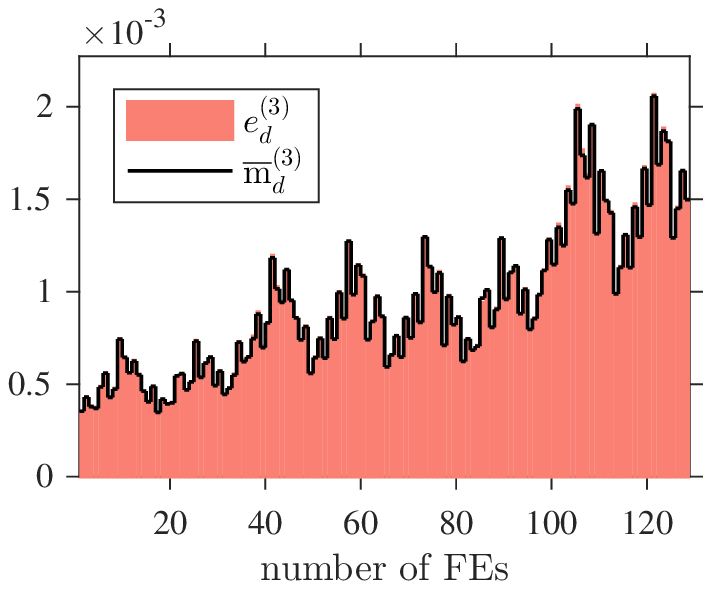}}\!
	\subfloat[4 REF: 512 EL, \newline $\ed = 1.44 \cdot 10^{-3}, \mdI = 1.44 \cdot 10^{-3}$]{
	\includegraphics[width=5.7cm]{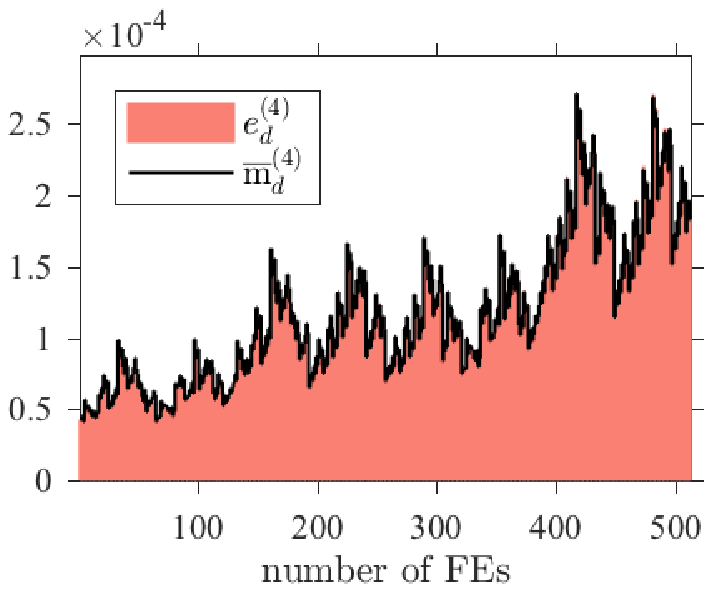} \label{fig:example-space-time-1d-t-unit-e-maj-ref-4}}
	\caption{Ex. \ref{ex:space-time-unit-1d-t}. 
	$\ed$ and $\mdI$ distributions after REF $\# = 2, 3, 4$.}
	\label{fig:example-space-time-1d-t-unit-e-maj-distr}
\end{figure}

Next, we consider the true error and majorant distributions obtained on each refinement 
step starting from the initial mesh $\mathcal{T}_{3 \times 3}$ (see Figure 
\ref{fig:example-space-time-1d-t-unit-e-maj-distr}). In Figures 
\ref{fig:example-space-time-1d-t-unit-e-maj-ref-2}--\ref{fig:example-space-time-1d-t-unit-e-maj-ref-4}, 
we illustrate $\ed$ and $\mdI$ after the refinement steps (REF) 2, 3, 4, 5,
respectively. Here, the elements are ordered by the certain numbering procedure in 
FE code implementation. The graphs confirm
that the indicator indeed manages to mimic the error distribution and to catch 
local jumps of the error very efficiently. Table \ref{tab:example-space-time-unit-1d-t-e-maj-ieff} 
confirms the efficiency of the total majorant. 
\begin{table}[!ht]
\centering
\footnotesize
\begin{tabular}{c|c|ccc}
$\#$ REF & $\#$ EL & $\error$ & $\maj{}$ & $\Ieff$\\
\midrule
1 & 8 & $3.5229 \cdot 10^{-1}$ & $4.0889 \cdot 10^{-1}$ & 1.08 \\
3 & 128 & $2.2969 \cdot 10^{-2}$ & $2.7215 \cdot 10^{-2}$ & 1.09 \\ 
5 & 2048 & $1.4393 \cdot 10^{-3}$ & $1.7209 \cdot 10^{-3}$ & 1.09 \\ 
7 & 131072 & $2.2493 \cdot 10^{-5}$ & $2.6969 \cdot 10^{-5}$ & 1.09 \\
9 & 2097152 & $1.4058 \cdot 10^{-6}$ & $1.6861 \cdot 10^{-6}$ & 1.10 \\
\end{tabular}
\caption{Ex. \ref{ex:space-time-unit-1d-t}. 
Total error, majorant, and efficiency index w.r.t. refinement steps.}
\label{tab:example-space-time-unit-1d-t-e-maj-ieff}
\end{table}

\end{example}

\newpage
\begin{example}
\label{ex:space-time-unit-2d-t}
\rm
Let us consider the same problem discussed in Ex. \ref{ex:incr-unit-2d-t} by 
the space-time discretisation scheme. The used FE spaces are as follows: 
$v \in \Pone$ and $\flux \in \Ptwo$. We consider two meshes obtained after the uniform
refinement on the steps REF 1, 2 (see  
Figures \ref{fig:example-space-time-2d-t-y-RT2-e-m-hist-ref-1}--\ref{fig:example-space-time-2d-t-y-RT2-e-m-hist-ref-2}). Again, the local error and the indicator distributions are shown element-wise. 
Table \ref{tab:example-space-time-unit-2d-t-e-maj-ieff} provides information about
the efficiency index of the majorant on every REF, which improved on average by 
$20\%$ in comparison to the time-stepping approach.

\begin{figure}[!ht]
	\centering
	\subfloat[1 REF: 48 EL \newline $\ed = 4.77 \cdot 10^{-2}, \maj{} = 5.60 \cdot 10^{-2}$]{
	\includegraphics[width=5.7cm]{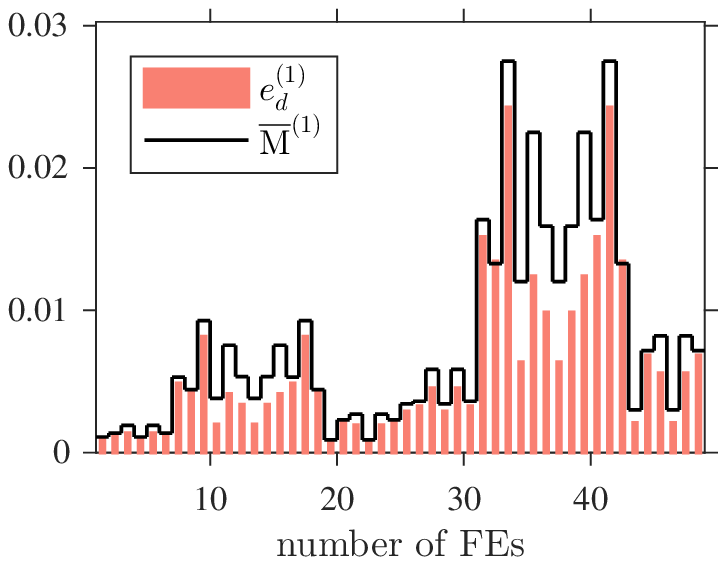}
	\label{fig:example-space-time-2d-t-y-RT2-e-m-hist-ref-1}}
	\subfloat[2 REF: 384 EL \newline $\ed = 1.73 \cdot 10^{-2}, \maj{} = 1.92 \cdot 10^{-2}$]{
	\includegraphics[width=5.4cm]{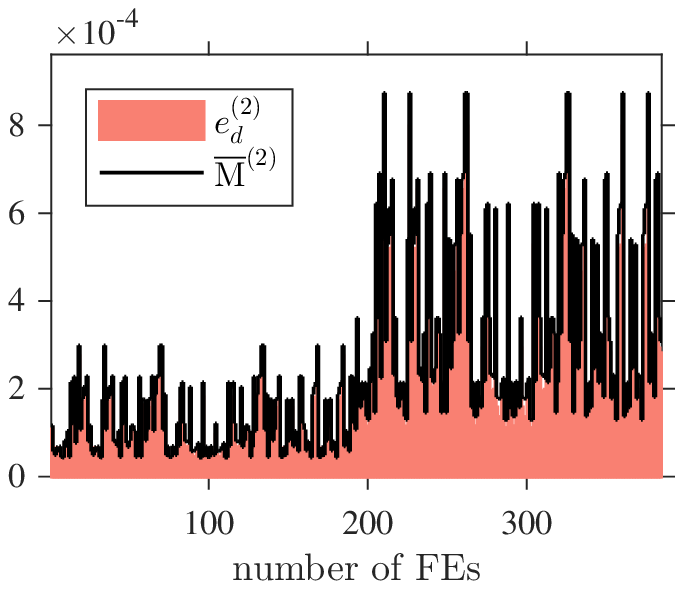}
	\label{fig:example-space-time-2d-t-y-RT2-e-m-hist-ref-2}}\\
	%
	\caption{Ex. \ref{ex:space-time-unit-2d-t}. 
	Distribution of the error and the majorant after refinement steps $\# = 1, 2$.}
\end{figure}

\begin{table}[!ht]
\centering
\footnotesize
\begin{tabular}{c|c|ccc}
$\#$ REF & $\#$ EL & $\error$ & $\maj{}$ & $\Ieff$\\
\midrule
1 & 48 & $4.7719 \cdot 10^{-2}$ & $5.6030 \cdot 10^{-2}$ & 1.08  \\
2 & 384 & $1.7268 \cdot 10^{-2}$ & $1.9175 \cdot 10^{-2}$ & 1.05 \\ 
3 & 3072 & $8.4558 \cdot 10^{-3}$ & $9.4992 \cdot 10^{-3}$ & 1.06 \\ 
4 & 24576 & $6.1427 \cdot 10^{-3}$ & $6.9789 \cdot 10^{-3}$ & 1.07 \\ 
5 & 196608 & $5.5479 \cdot 10^{-3}$ & $6.2617 \cdot 10^{-3}$ & 1.06 \\
\end{tabular}
\caption{Ex. \ref{ex:space-time-unit-2d-t}. 
Total error, the majorant, and the efficiency index w.r.t. to refinement steps.}
\label{tab:example-space-time-unit-2d-t-e-maj-ieff}
\end{table}

\end{example}

\begin{example}
\label{ex:space-time-example-2}
\rm
Finally, we construct an example with dissipating exact solution and 
zero RHS and investigate how the parameter $\sigma$ affects the quality of the 
error indication. This example demonstrates also the automated adaptive strategy 
suggested by the majorant.
The initial condition $u_0 = 6\, \sin \pi x$ determines the solution 
$u = 6\, \sin \pi x \, e^{-\tfrac{\pi^2 t}{\sigma}}$. For different $\sigma$, one 
obtains the solutions presented in Figure 
\ref{fig:example-space-time-example-2-exact-solution}.

\begin{figure}[!ht]
	\centering
	\subfloat[$\sigma = 1$]{
	\includegraphics[width=5.4cm, trim={0cm 0cm 0cm 0cm},clip]{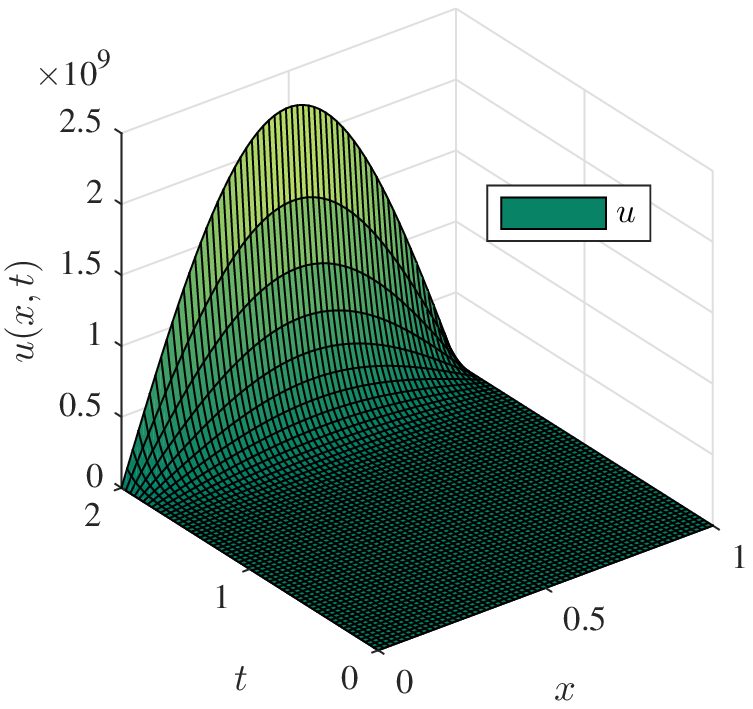}
	\label{fig:example-space-time-example-2-exact-solution-sigma-1}}\qquad
	\subfloat[$\sigma = 10$]{
	\includegraphics[width=5.4cm, trim={0cm 0cm 0cm 0cm},clip]{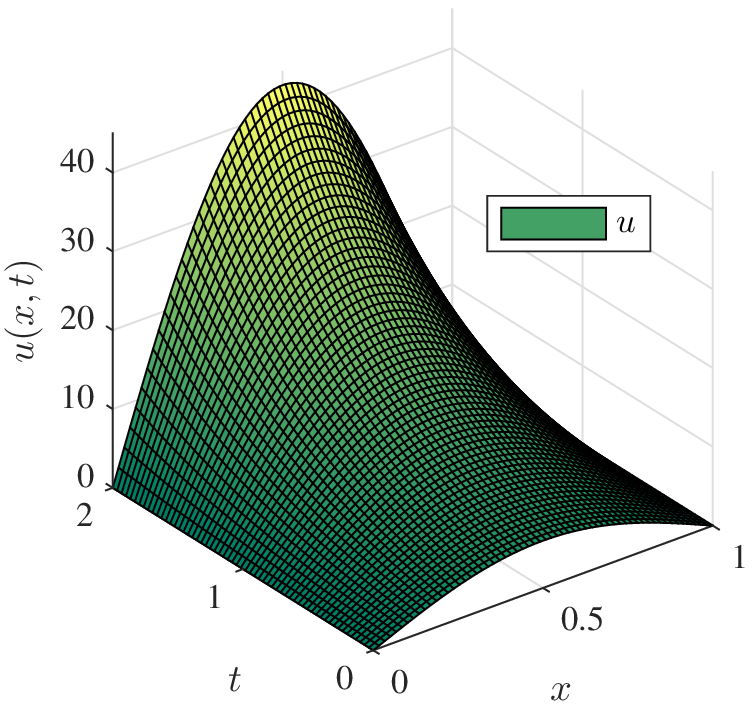}
	\label{fig:example-space-time-example-2-exact-solution-sigma-10}}
	\caption{Ex. \ref{ex:space-time-example-2}. Exact solutions for $\sigma = 1$ and $10$.}
	\label{fig:example-space-time-example-2-exact-solution}
\end{figure}

Assume first that $\sigma = 1$. Let the number of refinements be ${\rm REF} = 12$ 
and the bulk marking procedure have the parameter $\theta =0.3$. Having 
compared the meshes obtained while refining based on the error (Figures 
\ref{fig:example-space-time-example-2-bulk-0.3-meshes-a}--\ref{fig:example-space-time-example-2-bulk-0.3-meshes-d}) and the majorant distribution (Figures \ref{fig:example-space-time-example-2-bulk-0.3-meshes-b}--\ref{fig:example-space-time-example-2-bulk-0.3-meshes-d}), one observes that 
the majorant adequately detects the region of large gradients in the solution. 
Using the true error for the refinement, 
we obtain the following $\L{2}$ and energy norms of errors: $\|e\|^2_Q = 3.56 \cdot 10^{-7}$ 
and $\| \nabla_x e\|^2_Q = 2.02 \cdot 10^{-2}$. Whereas, the refinement based on $\mdI$ 
results into the errors $\|e\|^2_Q = 3.56 \cdot 10^{-7}$ and $\| \nabla_x e\|^2_Q = 2.02 \cdot 10^{-2}$. 
The efficiency index obtained after $12$ refinement steps equals to $1.25$, which is quite adequate 
number taking into account the jump in the exact solution.

\begin{figure}[!ht]
	\centering
	\subfloat[REF $\#$ 4: $\mbox{based on} \;\ed$ \newline 
	796 EL, 440 ND 
	]{
	\includegraphics[width=5.6cm, trim={6.0cm 0cm 6.4cm 1cm},clip]{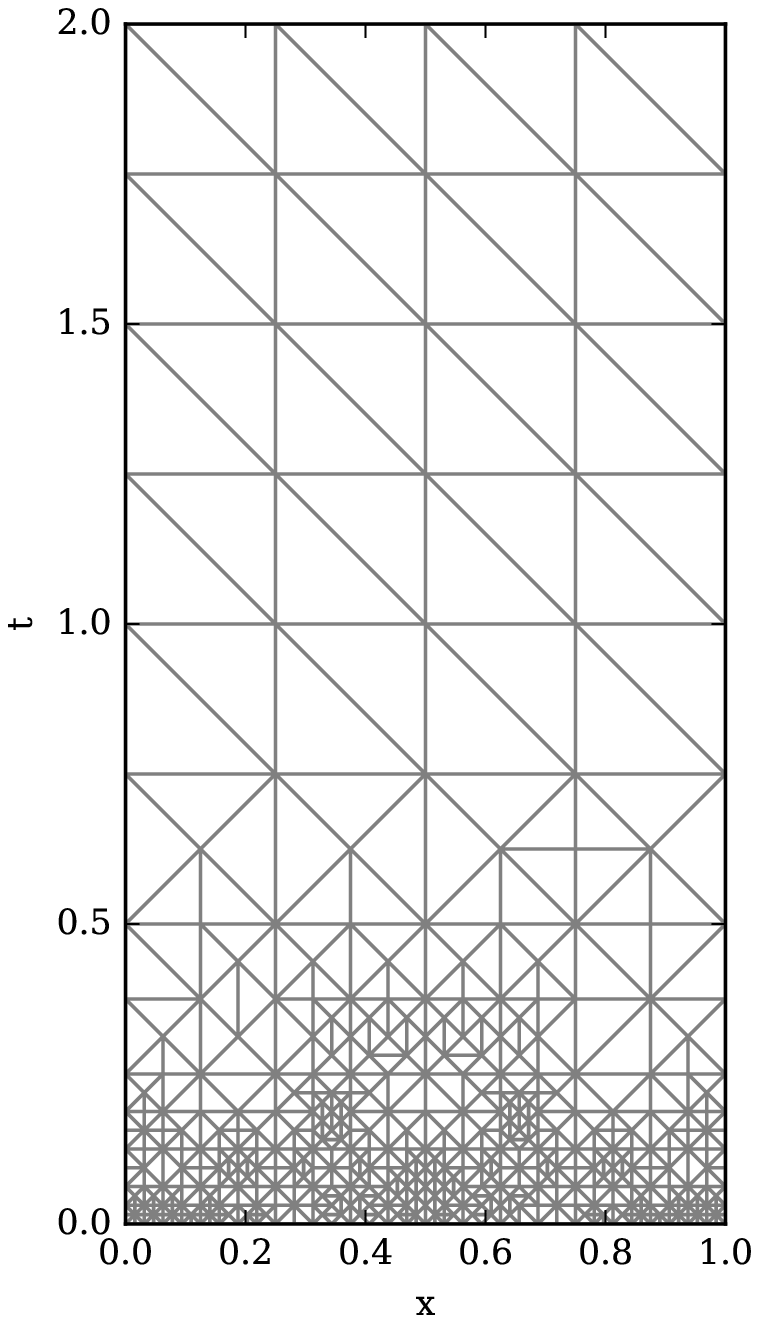}
	\label{fig:example-space-time-example-2-bulk-0.3-meshes-a}}\!
	\subfloat[REF $\#$ 4: $\mbox{based on} \; \mdI$ \newline 
	776 EL, 431 ND 
	]{
	\includegraphics[width=5.6cm,  trim={6.0cm 0 6.4cm 1cm},clip]{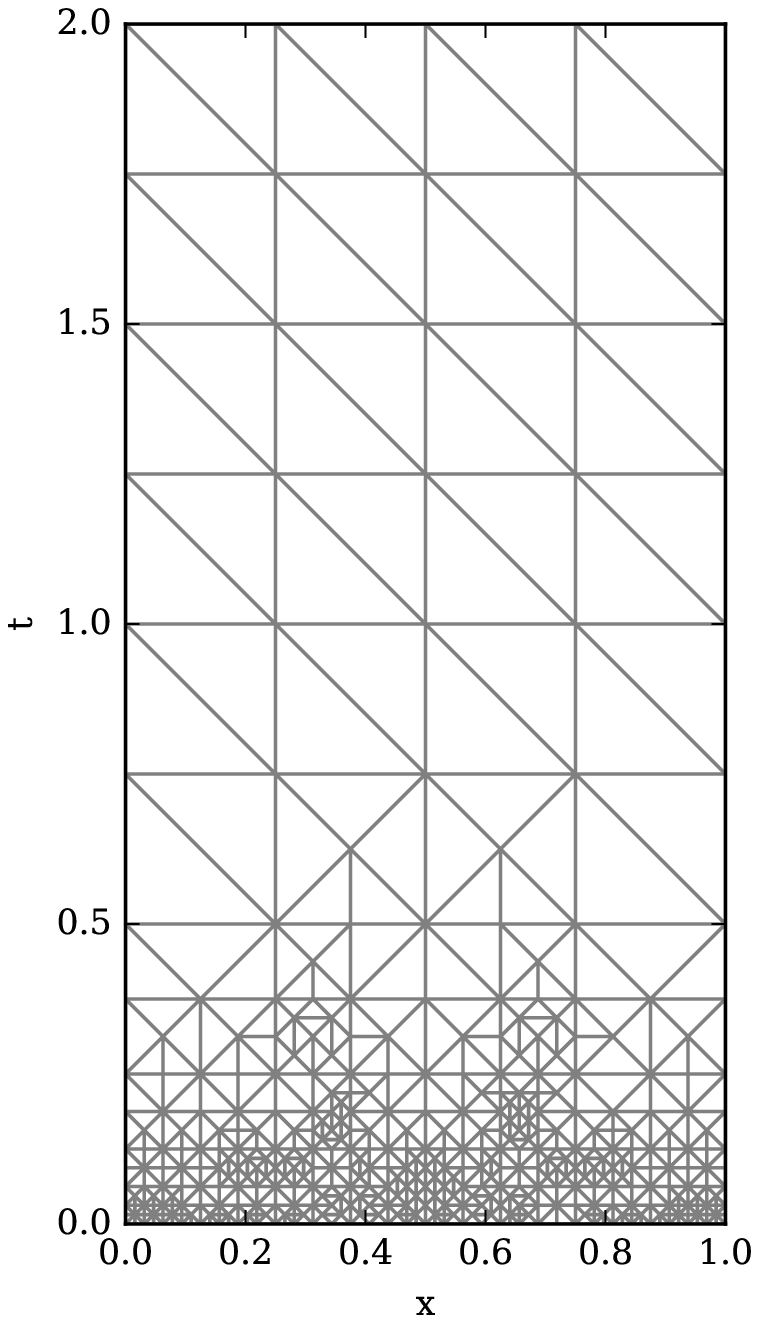}
	\label{fig:example-space-time-example-2-bulk-0.3-meshes-b}}\\[-10pt]
	%
	%
	\subfloat[REF $\#$ 6: $\mbox{based on} \;\ed$ \newline 
	3340 EL, 1739 ND 
	]{
	\includegraphics[width=5.6cm, trim={6.0cm 0 6.4cm 1cm},clip]{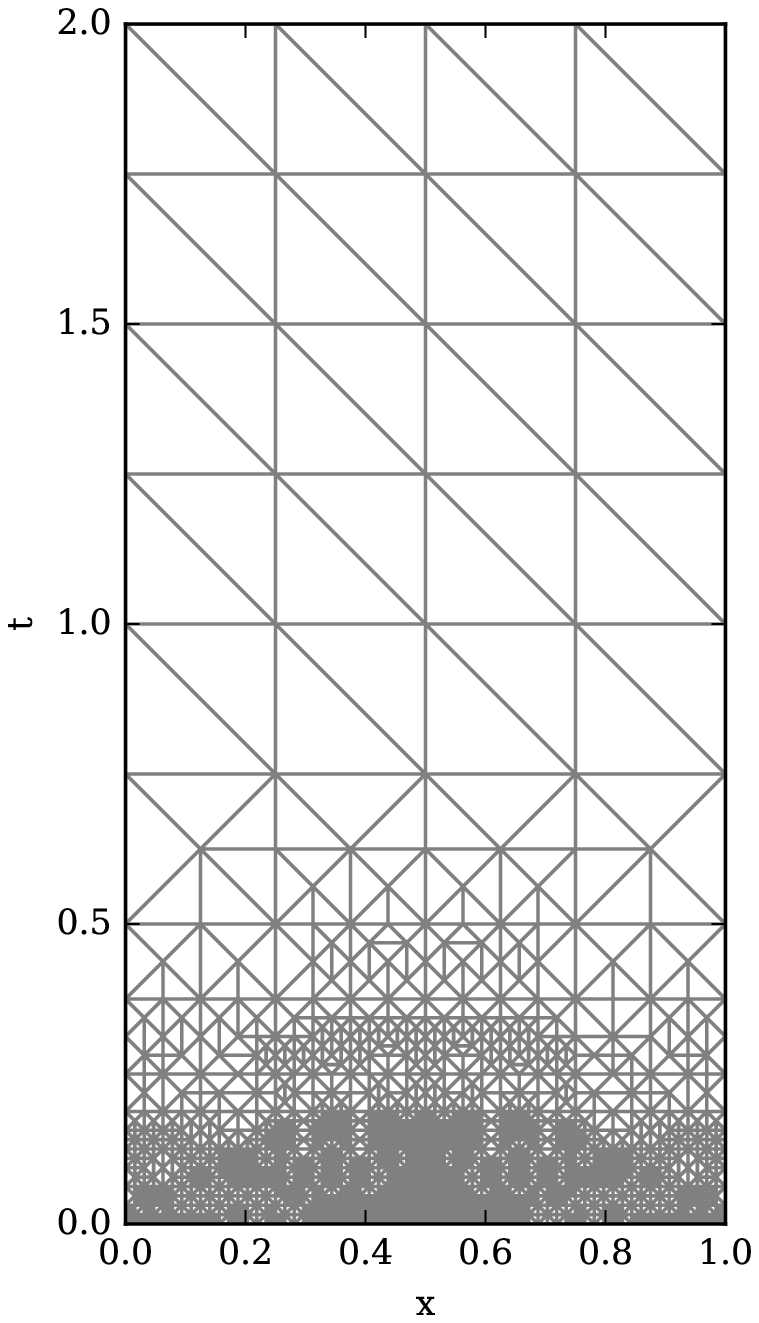}
	\label{fig:example-space-time-example-2-bulk-0.3-meshes-c}}\!
	\subfloat[REF $\#$ 6: $\mbox{based on} \; \mdI$ \newline 
	3374 EL, 1758 ND 
	]{
	\includegraphics[width=5.6cm,  trim={6.0cm 0 6.4cm 1cm},clip]{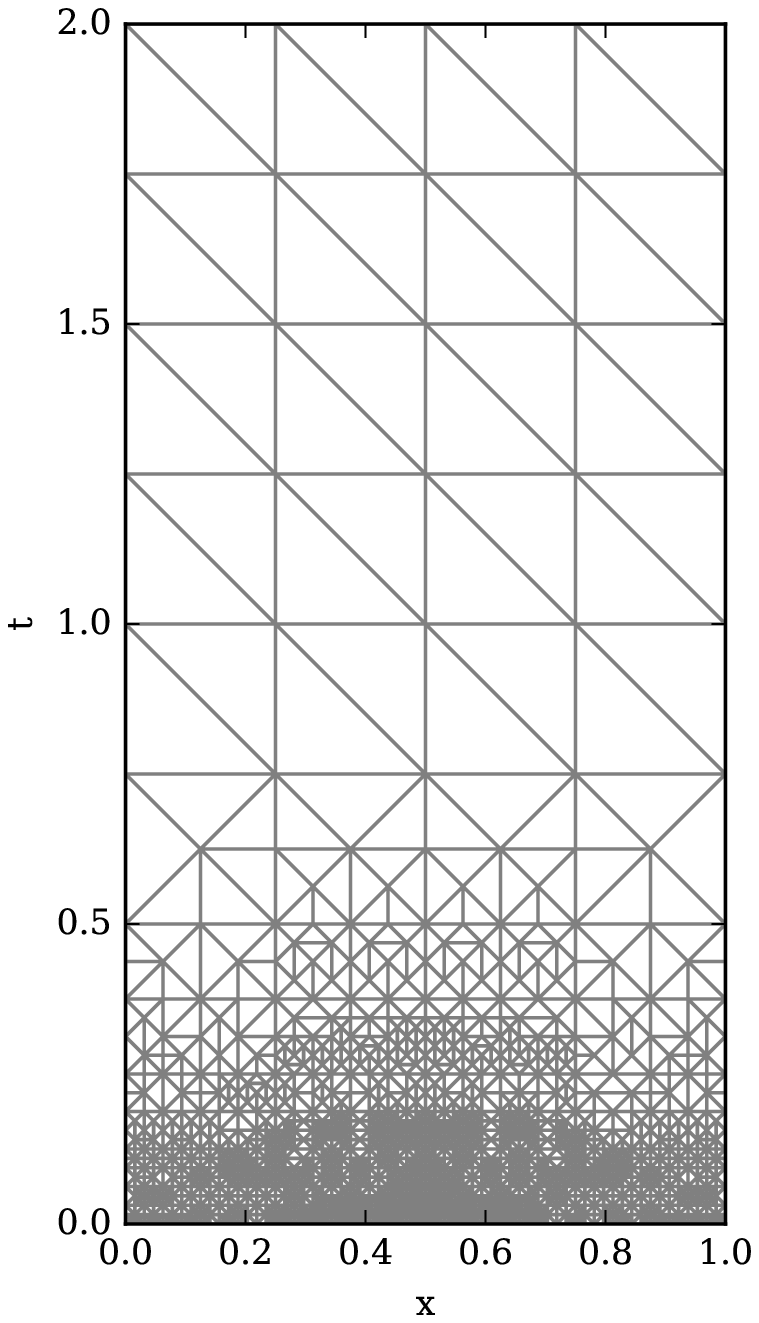}
	\label{fig:example-space-time-example-2-bulk-0.3-meshes-d}}\\
	\caption{Ex. \ref{ex:space-time-example-2} for $\sigma = 1$. 
	Space-time refinement based on the true energy error (a), (c) 
	and on the majorant (b), (d) using the bulk marker $\Marker_{0.3}$.}
	\label{fig:example-space-time-example-2-bulk-0.3-meshes}
\end{figure}

\begin{figure}[!ht]
	\centering
	\subfloat[REF $\#$ 6: $\mbox{based on} \;\ed$ \newline 
	2548 EL, 1321 ND 
	]{
	\includegraphics[width=5.6cm, trim={6.0cm 0cm 6.4cm 1cm},clip]{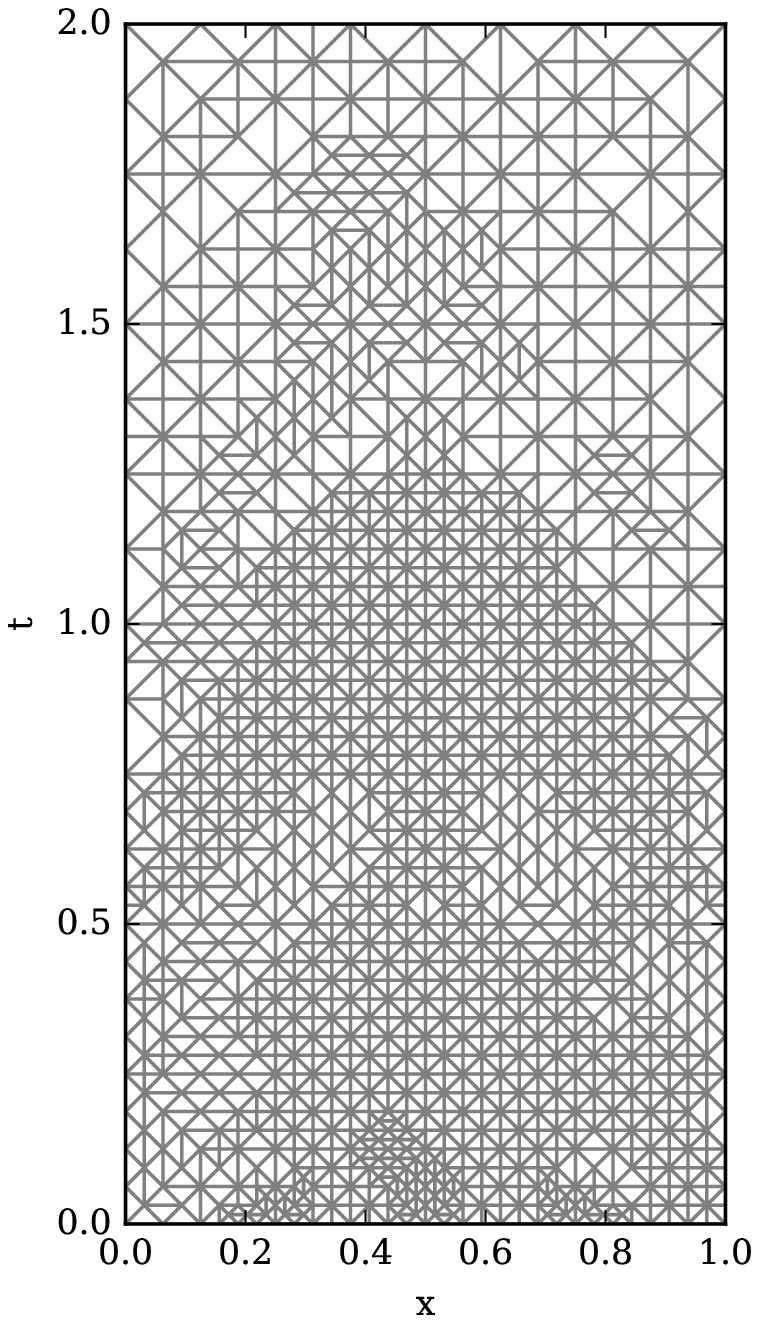}
	\label{fig:example-space-time-example-2-bulk-0.3-meshes-sigma-10-a}}\!
	\subfloat[REF $\#$ 6: $\mbox{based on} \; \mdI$ \newline 
	2660 EL, 1387 ND 
	]{
	\includegraphics[width=5.6cm,  trim={6.0cm 0 6.4cm 1cm},clip]{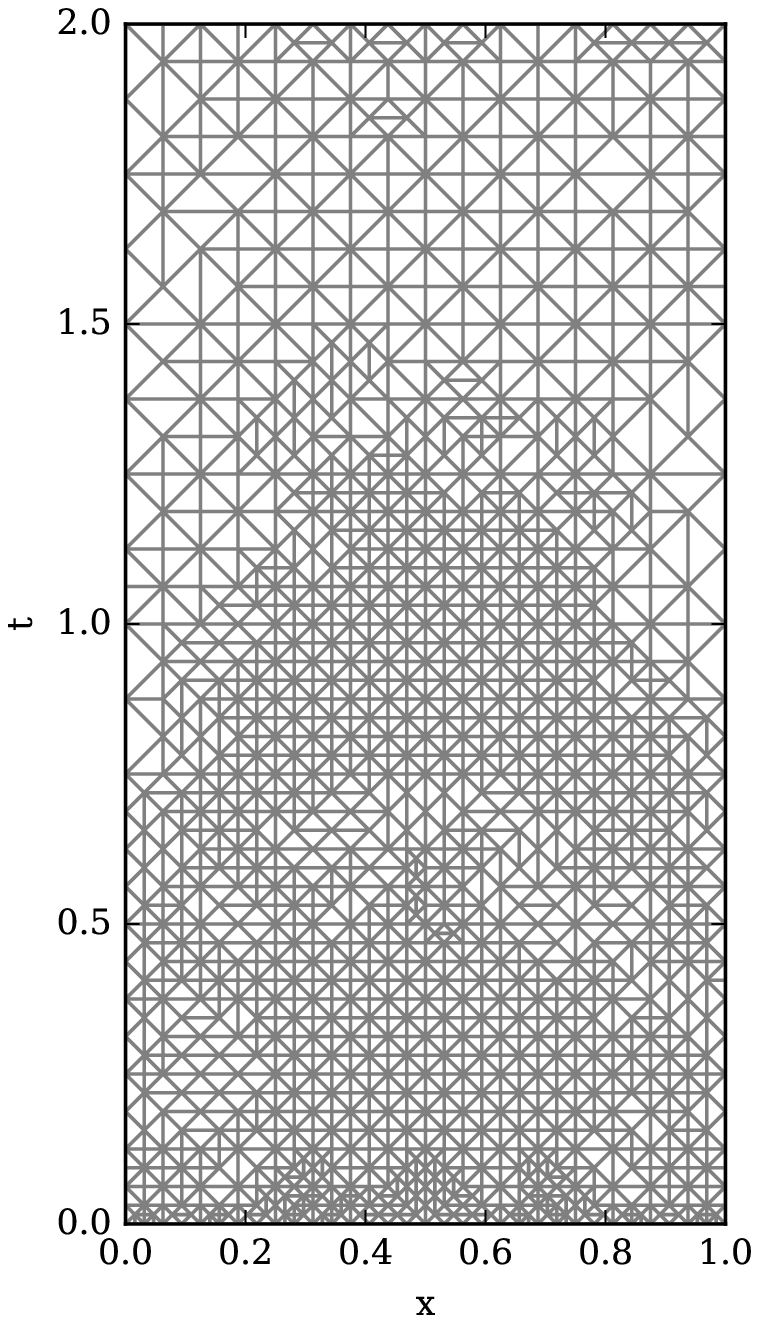}
	\label{fig:example-space-time-example-2-bulk-0.3-meshes-sigma-10-b}}\\[-10pt]
	%
	%
	\subfloat[REF $\#$ 8: $\mbox{based on} \;\ed$ \newline 
	8706 EL, 4437 ND 
	]{
	\includegraphics[width=5.6cm, trim={6.0cm 0 6.4cm 1cm},clip]{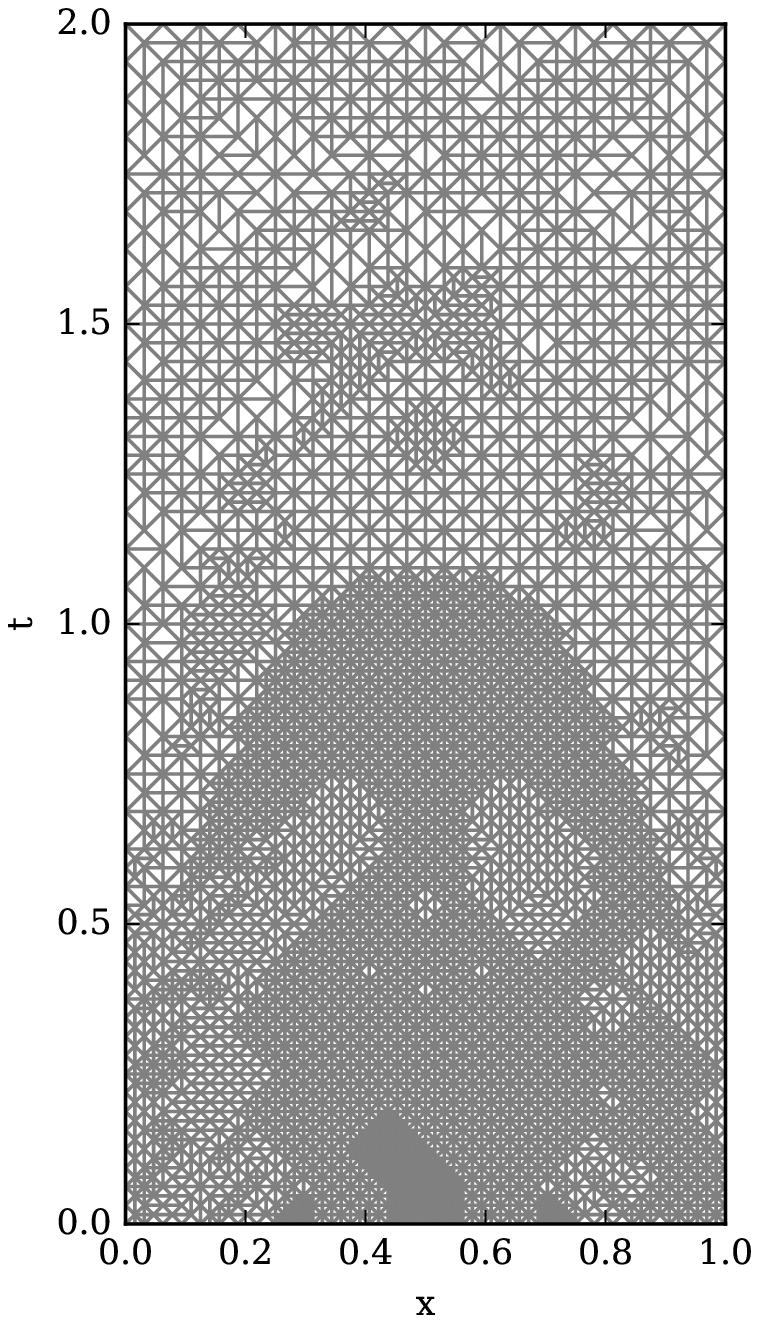}
	\label{fig:example-space-time-example-2-bulk-0.3-meshes-sigma-10-c}}\!
	\subfloat[REF $\#$ 8: $\mbox{based on} \; \mdI$ \newline 
	9274 EL, 4735 ND 
	]{
	\includegraphics[width=5.6cm,  trim={6.0cm 0 6.4cm 1cm},clip]{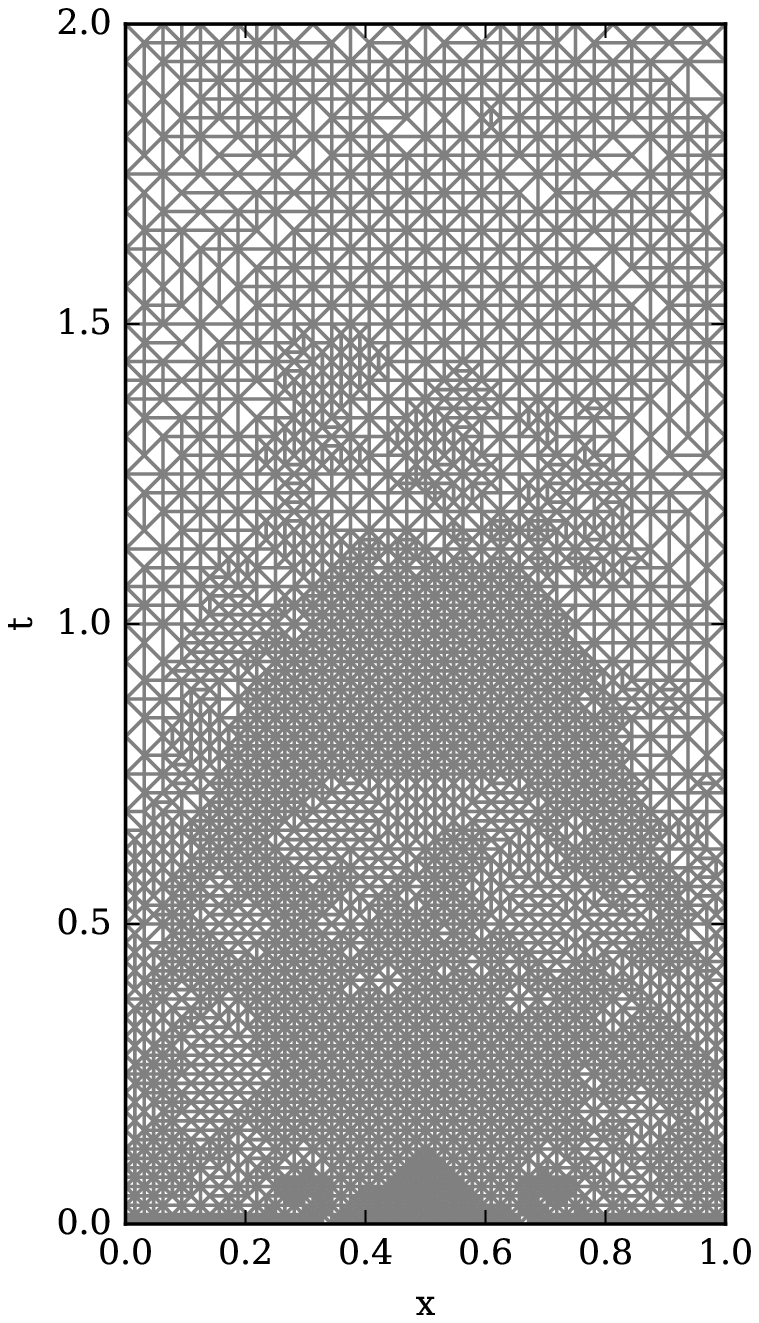}
	\label{fig:example-space-time-example-2-bulk-0.3-meshes-sigma-10-d}}\\
	\caption{Ex. \ref{ex:space-time-example-2} for $\sigma = 10$.
	Space-time refinement based on the true energy error (a), (c) 
	and on the majorant (b), (d) using the bulk marker $\Marker_{0.3}$.}
	\label{fig:example-space-time-example-2-bulk-0.3-meshes-sigma-10}
\end{figure}
 
\begin{table}[!t]
\centering
\footnotesize
\begin{tabular}{c|c|ccc}
$\#$ REF & $\#$ EL & $\error$ & $\maj{}$ & $\Ieff$\\
\midrule
      1 &       71 &  1.50 & 2.29 &       1.52 \\
      3 &      218 &  $8.35 \cdot 10^{-1}$ & 1.10 &       1.31 \\
      5 &      881 &  $3.90 \cdot 10^{-1}$ & $4.97 \cdot 10^{-1}$ &       1.27 \\
      7 &     3332 &  $2.02 \cdot 10^{-1}$ & $2.54 \cdot 10^{-1}$ &       1.26 \\
      9 &    11538 &  $1.10 \cdot 10^{-1}$ & $1.38 \cdot 10^{-1}$ &       1.26 \\
     11 &    38218 &  $6.28 \cdot 10^{-2}$ & $7.87 \cdot 10^{-2}$ &       1.25 \\
\end{tabular}
\caption{Ex. \ref{ex:space-time-example-2} for $\sigma = 1$. 
Total error, majorant, and efficiency index w.r.t. to refinement 
steps with the bulk marking $\theta = 0.3$.}
\label{tab:example-space-time-example-2-e-maj-ieff-sigma-1}
\end{table}

The same comparison of meshes is done for $\sigma = 10$ in
Figure \ref{fig:example-space-time-example-2-bulk-0.3-meshes-sigma-10}.
Table \ref{tab:example-space-time-example-2-e-maj-ieff-sigma-10} illustrates that the 
efficiency index is affected by the value of $\sigma$ but still stays relatively adequate. 
For a lage $\sigma$, the reliability term of the majorant $\mfI$
is more complicated to minimise due to the term $\sigma \, v_t$. 
In the space-time setting, it can be interpreted as a very strong convection term. In 
order to reduce the effect of $\mfI$ on the global majorant, one needs to use more 
sophisticated methods (tailored to the convection-dominated problems), which 
are beyond the focus of this paper. 

\begin{table}[!t]
\centering
\footnotesize
\begin{tabular}{c|c|ccc}
$\#$ REF & $\#$ EL & $\error$ & $\maj{}$ & $\Ieff$\\
\midrule
      1 &       71 &  1.55 & 3.18 &       2.05 \\
      3 &      219 &  $8.08 \cdot 10^{-1}$ & 1.65 &       2.05 \\
      5 &      741 &  $4.28 \cdot 10^{-1}$ & $8.60 \cdot 10^{-1}$ &       2.01 \\
      7 &     2557 & $2.29 \cdot 10^{-1}$ & $4.72 \cdot 10^{-1}$ &       2.06 \\
      9 &     8736 &  $1.23 \cdot 10^{-1}$ & $2.58 \cdot 10^{-1}$ &       2.10 \\
     11 &    29501 &  $6.68 \cdot 10^{-2}$ & $1.39 \cdot 10^{-1}$ &       2.09 \\
\end{tabular}
\caption{Ex. \ref{ex:space-time-example-2} for $\sigma = 10$. 
Total error, majorant, and efficiency index w.r.t. to refinement 
steps with the bulk marking $\theta = 0.3$.}
\label{tab:example-space-time-example-2-e-maj-ieff-sigma-10}
\end{table}

\end{example}

%% file: sections/conclusions.tex
\section{Conclusion}
\label{sec:conclusions}

In this paper, we have discussed guaranteed bounds of the distance to the 
exact solution of the evolutionary reaction-diffusion problem with the Dirichlet BC. 
We have shown that these error estimates are directly computable and rather efficient. 
The functional formulation of the estimates makes them flexible with respect to any 
discretisation method, such as considered in the current work, the time-stepping or the space-time schemes. 
Numerical experiments performed for both approaches confirmed that the 
estimates provide accurate bounds of the overall error and generate efficient 
indicators of the local error distribution. 

\section{Acknowledgments}

This joint work would not be possible without the initial support of the Department of 
Mathematical Information Technology of the University of Jyvaskyla. The authors also 
express gratitude to Prof. Sergey Repin for many helpful 
suggestions and fruitful discussions.

%% file: main-func-estim-paper.bbl
\def\cprime{$'$} \def\cprime{$'$} \def\cprime{$'$}
\begin{thebibliography}{10}

\bibitem{BabushkaRheinboldtSIAM1978}
I.~Babu{\v{s}}ka and W.~C. Rheinboldt.
\newblock Error estimates for adaptive finite element computations.
\newblock {\em SIAM J. Numer. Anal.}, 15(4):736--754, 1978.

\bibitem{BabushkaRheinboldt1978}
I.~Babu\v{s}ka and W.~C. Rheinboldt.
\newblock A-posteriori error estimates for the finite element method.
\newblock {\em Internat. J. Numer. Meth. Engrg.}, 12:1597--1615, 1978.

\bibitem{BeckerRannacher1996}
R.~Becker and R.~Rannacher.
\newblock A feed--back approach to error control in finite element methods:
  {B}asic approach and examples.
\newblock {\em East--West J. Numer. Math.}, 4(4):237--264, 1996.

\bibitem{Braess2001}
D.~Braess.
\newblock {\em Finite elements}.
\newblock Cambridge University Press, Cambridge, second edition, 2001.

\bibitem{BrennerScott1994}
S.~Brenner and R.~L. Scott.
\newblock {\em The mathematical theory of finite element methods}.
\newblock Springer, New York, 1994.

\bibitem{Ciarlet1978}
P.~G. Ciarlet.
\newblock {\em The finite element method for elliptic problems}.
\newblock North-Holland Publishing Co., Amsterdam-New York-Oxford, 1978.
\newblock Studies in Mathematics and its Applications, Vol. 4.

\bibitem{CourantFriedrichsLewy1967}
R.~Courant, K.~Friedrichs, and H.~Lewy.
\newblock On the partial difference equations of mathematical physics.
\newblock {\em IBM J. Res. Develop.}, 11:215--234, 1967.

\bibitem{Deuflhardetall1989}
P.~Deuflhard, P.~Leinen, and H.~Yserentant.
\newblock Concept of an adaptive hierarchical finite element code.
\newblock {\em Impact Computing Sci. Engrg.}, 1(1):3--35, 1989.

\bibitem{Dorfler1996}
W.~D{\"o}rfler.
\newblock A convergent adaptive algorithm for {P}oisson's equation.
\newblock {\em SIAM J. Numer. Anal.}, 33(3):1106--1124, 1996.

\bibitem{Friedman1964}
A.~Friedman.
\newblock {\em Partial differential equations of parabolic type}.
\newblock Prentice-Hall, Inc., Englewood Cliffs, N.J., 1964.

\bibitem{Friedrichs1937}
K.~Friedrichs.
\newblock On certain inequalities and characteristic value problems for
  analytic functions and for functions of two variables.
\newblock {\em Trans. Amer. Math. Soc.}, 41(3):321--364, 1937.

\bibitem{GaevskayaRepin2005}
A.~V. Gaevskaya and S.~I. Repin.
\newblock A posteriori error estimates for approximate solutions of linear
  parabolic problems.
\newblock {\em Springer, Differential Equations}, 41(7):970--983, 2005.

\bibitem{Gander2015}
M.~Gander.
\newblock 50 years of time parallel time integration.
\newblock In {\em Multiple Shooting and Time Domain Decomposition}, volume~16,
  pages 69--114. Springer-Verlag, Berlin, 2015.
\newblock Theory, algorithm, and applications.

\bibitem{GrossmannRoosStynes2007}
C.~Grossmann, H.-G. Roos, and M.~Stynes.
\newblock {\em Numerical treatment of partial differential equations}.
\newblock Universitext. Springer, Berlin, 2007.
\newblock Translated and revised from the 3rd (2005) German edition by Martin
  Stynes.

\bibitem{Hackbusch1984}
W.~Hackbusch.
\newblock Parabolic multigrid methods.
\newblock In {\em Computing methods in applied sciences and engineering, {VI}
  ({V}ersailles, 1983)}, pages 189--197. North-Holland, Amsterdam, 1984.

\bibitem{HortonVandewalle1995}
G.~Horton and S.~Vandewalle.
\newblock A space-time multigrid method for parabolic partial differential
  equations.
\newblock {\em SIAM J. Sci. Comput.}, 16(4):848--864, 1995.

\bibitem{Johnson2009}
C.~Johnson.
\newblock {\em Numerical solution of partial differential equations by the
  finite element method}.
\newblock Dover Publications Inc., Mineola, NY, 2009.
\newblock Reprint of the 1987 edition.

\bibitem{Karabelas2015}
E.~Karabelas.
\newblock {\em Space-time discontinuous Galerkin methods for cardic
  electro-mechanics}.
\newblock PhD thesis, Technische Universitat Graz, 2015.

\bibitem{Ladyzhenskaya1985}
O.~A. Ladyzhenskaya.
\newblock {\em The boundary value problems of mathematical physics}.
\newblock Springer, New York, 1985.

\bibitem{Ladyzhenskayaetall1967}
O.~A. Ladyzhenskaya, V.~A. Solonnikov, and N.N. Uraltseva.
\newblock {\em Linear and quasilinear equations of parabolic type}.
\newblock Nauka, Moscow, 1967.

\bibitem{Lang2001}
J.~Lang.
\newblock {\em Adaptive multilevel solution of nonlinear parabolic {PDE}
  systems}, volume~16 of {\em Lecture Notes in Computational Science and
  Engineering}.
\newblock Springer-Verlag, Berlin, 2001.
\newblock Theory, algorithm, and applications.

\bibitem{LangerMooreNeumueller2016a}
U.~Langer, S.~Moore, and M.~Neum\"uller.
\newblock Space-time isogeometric analysis of parabolic evolution equations.
\newblock {\em Comput. Methods Appl. Mech. Engrg.}, 306:342--363, 2016.

\bibitem{LoggMardalWells2012}
A.~Logg, K.-A. Mardal, and G.~N. Wells, editors.
\newblock {\em Automated solution of differential equations by the finite
  element method}, volume~84 of {\em Lecture Notes in Computational Science and
  Engineering}.
\newblock Springer, Heidelberg, 2012.
\newblock The FEniCS book.

\bibitem{LubichOstermann1987}
Ch. Lubich and A.~Ostermann.
\newblock Multigrid dynamic iteration for parabolic equations.
\newblock {\em BIT}, 27(2):216--234, 1987.

\bibitem{Malietall2014}
O.~Mali, P.~Neittaanm{\"a}ki, and S.~Repin.
\newblock {\em Accuracy verification methods}, volume~32 of {\em Computational
  Methods in Applied Sciences}.
\newblock Springer, Dordrecht, 2014.
\newblock Theory and algorithms.

\bibitem{MatculevichNeitaanmakiRepin2015}
S.~Matculevich, P.~Neittaanm{\"a}ki, and S.~Repin.
\newblock A posteriori error estimates for time-dependent reaction-diffusion
  problems based on the {P}ayne--{W}einberger inequality.
\newblock {\em AIMS}, 35(6):2659--2677, 2015.

\bibitem{MatculevichRepin2014}
S.~Matculevich and S.~Repin.
\newblock Computable estimates of the distance to the exact solution of the
  evolutionary reaction-diffusion equation.
\newblock {\em Appl. Math. and Comput.}, 247:329--347, 2014.

\bibitem{MatculevichRepinDiffUr2016}
S.~Matculevich and S.~Repin.
\newblock Estimates for the difference between exact and approximate solutions
  of parabolic equations on the basis of {P}oincar\'{e} inequalities for traces
  of functions on the boundary.
\newblock {\em Differential Equations}, 52(10):1355--1365, 2016.

\bibitem{MatculevichRepin2016}
S.~Matculevich and S.~Repin.
\newblock Explicit constants in poincar\'{e}-type inequalities for simplicial
  domains.
\newblock {\em Comput. Methods Appl. Math.}, 16(2):277--298, 2016.

\bibitem{NeittaanmakiRepin2004}
P.~Neittaanm{\"a}ki and S.~Repin.
\newblock {\em Reliable methods for computer simulation}, volume~33 of {\em
  Studies in Mathematics and its Applications}.
\newblock Elsevier Science B.V., Amsterdam, 2004.
\newblock Error control and a posteriori estimates.

\bibitem{Repin1997}
S.~Repin.
\newblock A posteriori estimates for approximate solutions of variational
  problems with strongly convex functionals.
\newblock {\em Problems of Mathematical Analysis}, 17:199--226, 1997.

\bibitem{Repin2000}
S.~Repin.
\newblock A posteriori error estimation for variational problems with uniformly
  convex functionals.
\newblock {\em Math. Comput.}, 69(230):481--500, 2000.

\bibitem{RepinDeGruyter2008}
S.~Repin.
\newblock {\em A posteriori estimates for partial differential equations},
  volume~4 of {\em Radon Series on Computational and Applied Mathematics}.
\newblock Walter de Gruyter GmbH \& Co. KG, Berlin, 2008.

\bibitem{Repin1999}
S.~I. Repin.
\newblock A unified approach to a posteriori error estimation based on duality
  error majorants.
\newblock {\em Math. Comput. Simulation}, 50(1-4):305--321, 1999.
\newblock Modelling '98 (Prague).

\bibitem{Repin2002}
S.~I. Repin.
\newblock Estimates of deviations from exact solutions of initial-boundary
  value problem for the heat equation.
\newblock {\em Rend. Mat. Acc. Lincei}, 13(9):121--133, 2002.

\bibitem{RepinTomar2010}
S.~I. Repin and S.~K. Tomar.
\newblock A posteriori error estimates for approximations of evolutionary
  convection-diffusion problems.
\newblock {\em J. Math. Sci. (N. Y.)}, 170(4):554--566, 2010.
\newblock Problems in mathematical analysis. No. 50.

\bibitem{StrangFix1973}
G.~Strang and G.~Fix.
\newblock {\em An analysis of the finite element method}.
\newblock Prentice Hall, Englewood Cliffs, 1973.

\bibitem{Stroud1974}
A.~H. Stroud.
\newblock {\em Numerical quadrature and solution of ordinary differential
  equations}.
\newblock Springer-Verlag, New York-Heidelberg, 1974.
\newblock A textbook for a beginning course in numerical analysis, Applied
  Mathematical Sciences, Vol. 10.

\bibitem{StroudSecrest1966}
A.~H. Stroud and Don Secrest.
\newblock {\em Gaussian quadrature formulas}.
\newblock Prentice-Hall, Inc., Englewood Cliffs, N.J., 1966.

\bibitem{Takizawaetall2012}
K.~Takizawa, K.~Schjodt, A.~Puntel, N.~Kostov, and T.~E. Tezduyar.
\newblock Patient-specific computer modeling of blood flow in cerebral arteries
  with aneurysm and stent.
\newblock {\em Comput. Mech.}, 50(6):675--686, 2012.

\bibitem{TakizawaTezduyar2011}
K.~Takizawa and T.~E. Tezduyar.
\newblock Multiscale space-time fluid-structure interaction techniques.
\newblock {\em Comput. Mech.}, 48(3):247--267, 2011.

\bibitem{TakizawaTezduyar2014}
K.~Takizawa and T.~E. Tezduyar.
\newblock Space-time computation techniques with continuous representation in
  time ({ST}-{C}).
\newblock {\em Comput. Mech.}, 53(1):91--99, 2014.

\bibitem{Thomee2006}
V.~Thom{\'e}e.
\newblock {\em Galerkin finite element methods for parabolic problems},
  volume~25 of {\em Springer Series in Computational Mathematics}.
\newblock Springer-Verlag, Berlin, second edition, 2006.

\bibitem{VandewallePiessens1992}
S.~Vandewalle and R.~Piessens.
\newblock Efficient parallel algorithms for solving initial-boundary value and
  time-periodic parabolic partial differential equations.
\newblock {\em SIAM J. Sci. Statist. Comput.}, 13(6):1330--1346, 1992.

\bibitem{Wloka1987}
J.~Wloka.
\newblock {\em Partial Differential Equations}.
\newblock Cambridge University Press, 1987.

\bibitem{Womble1990}
D.~E. Womble.
\newblock A time-stepping algorithm for parallel computers.
\newblock {\em SIAM J. Sci. Statist. Comput.}, 11(5):824--837, 1990.

\bibitem{Zeidler1990A}
E.~Zeidler.
\newblock {\em Nonlinear functional analysis and its applications. {II}/{A}}.
\newblock Springer-Verlag, New York, 1990.

\bibitem{Zeidler1990B}
E.~Zeidler.
\newblock {\em Nonlinear functional analysis and its applications. {II}/{B}}.
\newblock Springer-Verlag, New York, 1990.

\bibitem{ZienkiewiczZhu1987}
O.~C. Zienkiewicz and J.~Z. Zhu.
\newblock A simple error estimator and adaptive procedure for practical
  engineering analysis.
\newblock {\em Internat. J. Numer. Meth. Engrg.}, 24(2):337--357, 1987.

\bibitem{ZienkiewiczZhu1988}
O.~C. Zienkiewicz and J.~Z. Zhu.
\newblock Adaptive techniques in the finite element method.
\newblock {\em Commun. Appl. Numer. Methods}, 4:197--204, 1988.

\end{thebibliography}
